%% file: Newman_OptimalTensorAlgebras_ArXiv.tex
\documentclass[hidelinks,onefignum,onetabnum]{siamart220329}

\input{ex_shared}

\ifpdf
\hypersetup{
  pdftitle={Optimal Tensor Algebras},
  pdfauthor={E. Newman and K. Keegan}
}
\fi

\begin{document}

\maketitle

\begin{abstract}
Recent advances in {matrix-mimetic} tensor frameworks have made it possible to preserve linear algebraic properties for multilinear data analysis and, as a result, to obtain optimal representations of multiway data. Matrix mimeticity arises from interpreting tensors as operators that can be multiplied, factorized, and analyzed analogous to matrices. Underlying the tensor operation is an algebraic framework parameterized by an invertible linear transformation. The choice of linear mapping is crucial to representation quality and, in practice, is made heuristically based on expected correlations in the data. However, in many cases, these correlations are unknown and common heuristics lead to suboptimal performance. In this work, we simultaneously learn optimal linear mappings and corresponding tensor representations without relying on prior knowledge of the data. Our new framework explicitly captures the coupling between the transformation and representation using variable projection.  We preserve the invertibility of the linear mapping by learning orthogonal transformations with Riemannian optimization. We provide original theory of uniqueness of the transformation and convergence analysis of our variable-projection-based algorithm. We demonstrate the generality of our framework through numerical experiments on a wide range of applications, including financial index tracking, image compression, and reduced order modeling.  We have published all the code related to this work at \url{https://github.com/elizabethnewman/star-M-opt}.
\end{abstract}




\section{Introduction}
\label{sec:introduction}

In recent years, multilinear or \emph{tensor}-based analogs to the singular value decomposition (SVD) have revolutionized multiway data analysis.  Tensor decompositions seek to represent data and operators using \emph{interpretable}, \emph{compressible}, \emph{computationally efficient} strategies that avoid vectorization, respect high-dimensional correlations (e.g., spatio-temporal),  and exploit structure cleverly.  The right factorization for a particular application is often adapted from one of three classical strategies: Canonical Polyadic/Parallel Factor~\cite{Hitchcock1927, Harshman1970, CarrollChang1970}, Tucker~\cite{Tucker1966, LathauwerMoorVandewalle2000}, and Tensor-Train~\cite{Oseledets2011}.  Despite the popularity of these decompositions, all suffer from a so-called ``curse of multidimensionality;'' that is, linear algebra concepts break in high dimensions, including the optimality of low-rank approximations (Eckart-Young Theorem~\cite{Eckart1936}).

In this paper, we focus on  a  \emph{matrix-mimetic} tensor framework that preserves well-known concepts from matrix algebra (e.g., transpose, orthogonality, ...).  The key to the tensor algebraic framework is a well-defined tensor-tensor product that looks and feels like matrix-matrix multiplication. The tensor operation was first introduced as the specific $t$-product in~\cite{KILMER2011641, KilmerMartinPerrone2008, KilmerEtAl2013:thirdOrderTensorOperators} and later extended to a family of tensor-tensor products, called the $\starM$-product ($\starM$ is pronounced ``star-M'')~\cite{KernfeldKilmer2015}.  Recent work proved the tensor SVD under the $\starM$-product ($t$-SVDM) satisfied an Eckart-Young Theorem and produced provably superior representations compared to the matrix SVD and the classical tensor decompositions~\cite{Kilmer2021:pnas}.  The unique optimality properties of $\starM$-framework have had success in engineering applications~\cite{10125009, 9313596, 8421595}, inverse problems~\cite{Newman2020, ma2020randomized}, machine learning~\cite{Newman2024:TNN}, and more.

Each $\starM$-product is parameterized by an invertible linear transformation $\bfM$, which in turn induces a {tensor algebra}. The choice of algebra can significantly impact the quality of representations. Typically, the algebra is chosen heuristically based on a priori knowledge of multilinear correlations in the data (e.g., periodic temporal behavior~\cite{Keegan2022}).  However, often no such prior knowledge exists. The goal of this paper is to \emph{learn optimal linear transformations} that harness implicit correlations in multiway data and lead to compressible and accurate $\starM$-based representations. 

\subsection{Our Contributions}

Our main contribution is $\starM$-optimization, a new algorithm to learn a linear transformation and desired tensor representation simultaneously. The core of $\starM$-optimization leverages matrix-mimetic optimality properties through variable projection, explicitly coupling the representation and the linear mapping.  In some cases, the variable projection formulation necessitates differentiation through tensor factorizations; hence, we introduce new formulas for differentiating through the $t$-SVDM (\Cref{sec:tsvdmDerivatives}). To analyze the uniqueness of optimal transformations, we develop original theory regarding invariants of the $\starM$-product (\Cref{sec:invariantProduct}) and two prototype problems (\Cref{sec:uniqueInvariancePrototypeShort}).  We then provide theoretical guarantee of convergence $\starM$-optimization for a nonconvex objective function (\Cref{sec:analysis}). Our experiments demonstrate the breadth of the $\starM$-optimization framework and provide concrete intuition about the role of the linear transformations in forming quality tensor representations.  To ensure reproducibility of the results, all code and reproducible experiments are available in a public repository at \url{https://github.com/elizabethnewman/star-M-opt}.

\subsection{Related Work} 
Matrix completion seeks to recover a matrix of the lowest possible rank that exactly matches partially-observed data, which is an NP-Hard problem. Practical algorithms form a convex relaxation of this problem by minimizing the matrix nuclear norm (i.e., the sum of the singular values) iteratively~\cite{Recht2010:nuclearNorm}.  In~\cite{ZhangEtAl2014:tensorNuclearNorm}, the authors introduced a new tensor nuclear norm (TNN) based on the $t$-product and obtained quality low-tensor-rank approximations of partially-observed video data. Subsequent works extended tensor completion to the more general $\starM$-product~\cite{LuPengWei2019, SongEtAl2020:tensorNuclearNorm}.  In~\cite{KongLuLin2021:TensorQRank}, the authors proposed learning the $\starM$-product as a subproblem of tensor completion and demonstrated success on a variety of imaging datasets.  Due to the nonconvexity of the subproblem, the authors developed estimates via variance maximization and approximate Riemannian optimization to tractably update the linear transformation. Additional work in~\cite{WuEtAl2022} introduced a Schatten-$p$ quasi-norm regularization term to learn transformations that promote low-rank structure for TNN minimization. Further extensions have included learning nonlinear tensor-tensor products via neural networks for tensor completion~\cite{Luo_2022}. 

Learning the tensor-tensor product has thus far been restricted to tensor completion problems and the proposed algorithms have not fully made use of the optimality properties of the $\starM$-framework. Our paper takes a notably different perspective that places learning the transformation as the central goal, which enables a more general problem formulation, a wider range of potential applications, and, importantly, a strong algebraic foundation.  

\subsection{Outline of the Paper}

The paper is organized as follows.  In~\Cref{sec:background}, we introduce notation and foundational concepts for the $\starM$-framework. In~\Cref{sec:starMOpt}, we present our new $\starM$-optimization framework for learning tensor-tensor products and provide theoretical convergence guarantees.  In~\Cref{sec:numerical}, we demonstrate the broad applicability of our framework through numerical experiments on financial index tracking, image compression, and reduced order modeling. In~\Cref{sec:conclusions}, we conclude and discuss several future directions.

\section{Background and Notation}
\label{sec:background}

We first introduce the tensor mechanics with standard notation in~\Cref{sec:notation}, the $\starM$-product in~\Cref{sec:mimeticDefinitions}, and the $t$-SVDM and Eckart-Young Theorem for tensors in~\Cref{sec:tsvdm}. 
In~\Cref{sec:underlyingAlgebra}, we connect the tensor mechanics to the algebraic structure that underlies the $\starM$-product. We conclude with new theory about invariants of the $\starM$-product in~\Cref{sec:invariantProduct}. 

\subsection{Notation}
\label{sec:notation}
Tensors are multidimensional arrays of data. The \emph{order} of the tensor refers to the number of indices required to describe a single entry. Scalars (lowercase letters $a$) are order-$0$ tensors, vectors (bold lowercase $\bfa$) are order-$1$, matrices (bold uppercase $\bfA$) are order-$2$, and tensors (bold uppercase script $\TA$) are order-$3$ or higher.    
While we focus on real-valued, order-$3$ tensors, the background concepts extend to complex-valued, higher-order tensors~\cite{Kilmer2021:pnas, Keegan2022}. 

\begin{remark}
Throughout the paper, we will use {\sc Matlab} indexing notation; for example, $\bfA_{:,j}$ or $\bfA(:,j)$ denotes that $j$-th column of a matrix.  
\end{remark}

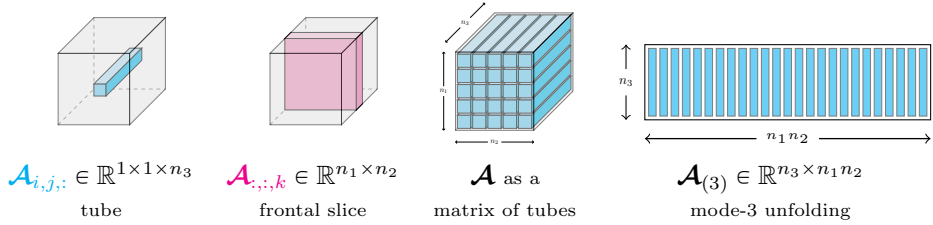
\begin{figure}
\centering
\def\w{0.15}
    \begin{tabular}{ccccc}
		\includegraphics[width=\w\linewidth]{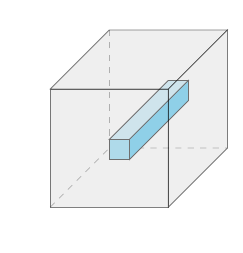} 
	& 	\includegraphics[width=\w\linewidth]{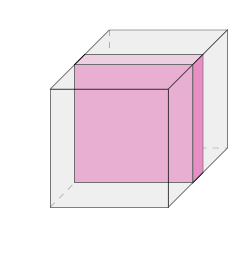}
 	& 	\includegraphics[width=\w\linewidth]{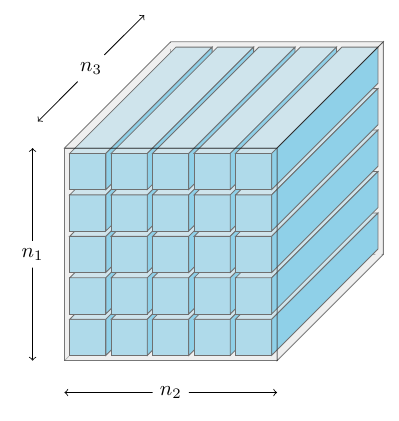}
	&	\begin{tikzpicture}
			\def\w{0.15}
			\def\s{0.05}
			\def\bb{0.25}
			\def\b{0.5}
	
			\draw (0,0) rectangle (25*\w+\s,1);

			\foreach \i in {0,...,24}{
				\draw[gray, fill=cyan!50] (\i*\w+\s,\s) rectangle (\i*\w+\w,1-\s);
			}

			\node (n1n2) at (25*\w/2+\s/2,-\bb) {\tiny $n_1n_2$};
			\draw[<-, line width=0.5pt]  (0,-\bb) -- (n1n2);
			\draw[->,  line width=0.5pt] (n1n2) -- (25*\w+\s,-\bb);
	
			\node (n3) at (-\bb,0.5) {\scalebox{0.5}{$n_3$}};
			\draw[<-]  (-\bb,0) -- (n3);
			\draw[->] (n3) -- (-\bb,1);
		
			\end{tikzpicture}\\

		$\cyan{\TA_{i,j,:}} \in \Rbb^{1\times 1\times n_3}$ 
	& 	$\magenta{\TA_{:,:,k}} \in \Rbb^{n_1\times n_2}$
         & 	$\TA$ \scriptsize as a
         & 	$\TA_{(3)}\in \Rbb^{n_3\times n_1n_2}$\\
		\scriptsize tube 
	&	\scriptsize  frontal slice
        &	\scriptsize matrix of tubes
        &	\scriptsize mode-$3$ unfolding
	\end{tabular}
	
	\caption{Visualization of a third-order tensor  $\TA\in \Rbb^{n_1\times n_2\times n_3}$ as a matrix of tubes and various partitions and unfoldings.}
	\label{fig:tensorNotation}
\end{figure}

Like rows and columns of matrices, we can partition tensors along various dimensions (\Cref{fig:tensorNotation}). 
A \emph{tube} is a vector lying along the third-dimension\footnote{Because tubes are one-dimensional arrays, we will denote them with bold lowercase letters $\bfa$. 
Whether the notation corresponds to a column vector or a tube will be clear from context.}. Given a tensor $\TA \in \Rbb^{n_1\times n_2\times n_3}$, the $(i,j)$-tube  is $\TA_{i,j,:} \in \Rbb^{1\times 1\times n_3}$. We can interpret a tensor as a \emph{matrix of tubes}, which will be a useful perspective in this paper.  Slices are matrices oriented along different axes of the tensor.  The $k$-th frontal slice is a matrix in the standard notation $\TA_{:,:,k} \in \Rbb^{n_1\times n_2}$.  
Using frontal slices, we define the tensor Frobenius norm as $\|\TA\|_F^2 = \sum_{k=1}^{n_3} \|\TA_{:,:,k}\|_F^2$. 

We now present two methods to matricize and act upon tensors.  

\begin{definition}[mode-$3$ unfolding]\label{def:mode3unfold}
Given a tensor $\TA\in \Rbb^{n_1\times n_2\times n_3}$, the mode-$3$ unfolding is a matrix $\TA_{(3)} \in \Rbb^{n_3\times n_1n_2}$ whose columns correspond to tubes of the tensor.  Explicitly, for $i=1,\dots, n_1$ and $j=1,\dots, n_2$,
    \begin{align}
    (\TA_{(3)})_{:,K(i,j)} = \myVec(\TA_{i,j,:}) \qquad \text{where} \qquad K(i,j) = i + (j - 1)n_1.
    \end{align}
Here, $\myVec: \Rbb^{1\times 1\times n_3} \to \Rbb^{n_3}$ reshapes tubes into column vectors. See~\Cref{fig:tensorNotation}.
\end{definition}

We can act upon matricized tensors as follows.

\begin{definition}[mode-$3$ product]\label{def:mode3product}
Given a third-order tensor $\TA\in \Rbb^{n_1\times n_2\times n_3}$ and matrix $\bfM\in \Rbb^{p\times n_3}$,  the mode-$3$ product applies $\bfM$ along each of the tubes; i.e.,
    \begin{align}
    \TA \times_3 \bfM = \myFold_3(\bfM \TA_{(3)})
    \end{align}
where $\myFold_3(\cdot)$ reverses the mode-$3$ unfolding.  
\end{definition}

\begin{remark}
\Cref{def:mode3unfold,def:mode3product} can be generalized to any mode of a tensor; see~\cite{KoldaBader2009} for details.  
\end{remark}

\subsection{Matrix-Mimetic Definitions}
\label{sec:mimeticDefinitions}

In this paper, we consider tensors to be \emph{$t$-linear operators}. By viewing tensors as matrices with tubal entries, tensor-tensor products are analogous matrix-matrix products, except that the individual multiplied entries are tubes rather than scalars. 
 
Hence, we first define tubal multiplication. 
\begin{definition}[$\starM$-tubal multiplication]\label{def:starMtubal}
Given tubes $\bfa, \bfb\in \Rbb^{1\times 1\times n_3}$ and an invertible transformation matrix $\bfM\in \Rbb^{n_3\times n_3}$, the $\starM$-product for tubes is 
    \begin{align}
        \bfa \starM \bfb = (\widehat{\bfa} \odot \widehat{\bfb})\times_3\bfM^{-1}.
    \end{align}
where $\widehat{\bfc} = \bfc\times_3 \bfM$ and $\odot$ denotes the Hadamard pointwise product.
We say $\bfc$ is in the \emph{spatial domain} and $\widehat{\bfc}$ is in the \emph{transform domain}. 
\end{definition}

\begin{figure}
\centering
    \begin{tabular}{c}
    \includegraphics[width=0.5\linewidth]{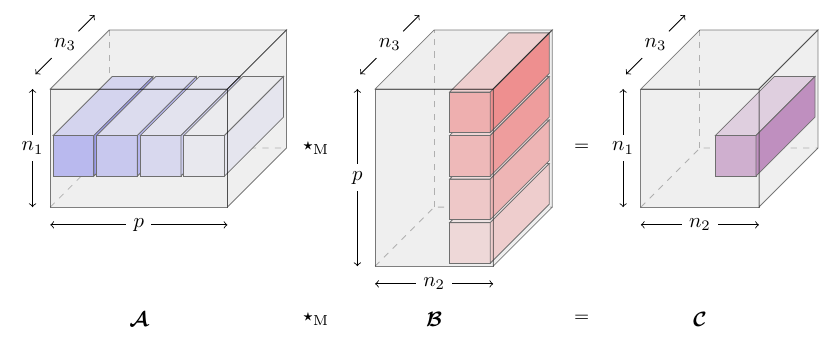} \\
    \includegraphics[width=0.7\linewidth]{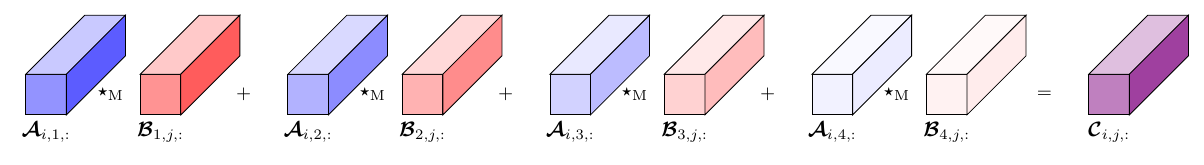}
    \end{tabular}

    \caption{Illustration of $\starM$-product (top) and computation of the $(i,j)$-tube (bottom).}
    \label{fig:starMProduct}
\end{figure}

Building from tubal multiplication, we define the tensor-tensor product as follows:

\begin{definition}[$\starM$-product]\label{defn:starMProduct}
 Given tensors $\TA \in \Rbb^{n_1\times p\times n_3}$ and $\TB\in \Rbb^{p\times n_2\times n_3}$ and an invertible matrix $\bfM\in \Rbb^{n_3\times n_3}$, the $\starM$-product is defined as
        \begin{align}
    	 	(\TA \starM \TB)_{i,j,:} = \sum_{k=1}^p \TA_{i,k,:} \starM \TB_{k,j,:} \qquad \text{for $i=1,\dots, n_1$ and $j=1,\dots, n_2$.}
        \end{align}
In practice, we parallelize the $\starM$-product over the third dimension via
	\begin{align}
		\TA \starM \TB = (\widehat{\TA} \smalltriangleup \widehat{\TB}) \times_3 \bfM^{-1}
	\end{align}
where $\widehat{\TC} = \TC\times_3 \bfM$ and $\smalltriangleup$ denotes the facewise product 
	\begin{align}\label{eq:facewise}
		(\widehat{\TA} \smalltriangleup \widehat{\TB})_{:,:,k} = \widehat{\TA}_{:,:,k} \widehat{\TB}_{:,:,k} \qquad \text{for $k=1,\dots,n_3$}. 
	\end{align}
\end{definition}
\Cref{fig:starMProduct} illustrates \Cref{defn:starMProduct}.  
We mention two crucial observations about the $\starM$-product.  First,~\Cref{defn:starMProduct} is exactly the definition of matrix-matrix multiplication if we were to remove the third dimension. This is an example of matrix-mimeticity unique to our tensor framework.  Second, the $\starM$-product is actually a family of tensor-tensor products, and each choice of $\bfM$ yields a different tensor-tensor product.  

From the foundation of the $\starM$-product, we build matrix-mimetic definitions of familiar linear algebra concepts.  Unless otherwise stated, all of the following definitions assume we are given an invertible, real-valued transform matrix $\bfM\in \Rbb^{n_3\times n_3}$.  Here, we provide the only the definitions needed to understand the $t$-SVDM (\Cref{sec:tsvdm}). For additional matrix-mimetic definitions, we recommend~\cite{KilmerEtAl2013:thirdOrderTensorOperators}.  

\begin{definition}[$\starM$-transpose\footnote{This definition changes if $\bfM$ is complex-valued.  See~\cite{Kilmer2021:pnas} for details.}]\label{def:starMtranspose} 
Given a tensor $\TA \in \Rbb^{n_1\times n_2\times n_3}$, the $\starM$-transpose $\TA^\top \in \Rbb^{n_2\times n_1\times n_3}$ is given by transposing the frontal slices of $\TA$; that is, 
\begin{align}
    (\TA^\top)_{:,:,k} = (\TA_{:,:,k})^\top \quad \text{for $k=1,\dots, n_3$.}
\end{align} 
\end{definition}

Just as the identity matrix is diagonal with ones along the diagonal, we present equivalent notions for tensors and ultimately define the identity tensor.

\begin{definition}[f-diagonal]
A tensor $\TD\in \Rbb^{m\times m\times n_3}$ is facewise diagonal or f-diagonal if every frontal slice is a diagonal matrix.  As a result, the only  non-zero tubes appear along its main diagonal; that is, for $i,j=1,\dots,m$,
	\begin{align}
	\TD_{i,j,:} = \begin{cases}
				\bfd_{i} & i = j\\
				\bf0 & i\not=j.
				\end{cases} 
			\qquad \text{for any $\bfd_i\in \Rbb^{1\times 1\times n_3}$.}
	\end{align}
\end{definition}

\begin{definition}[$\starM$-identity tube]\label{def:starMidentitytube} 
The $\starM$-identity tube, $\bfe \in \Rbb^{1 \times 1\times n_3}$, is $ \bfe = \mathbf{1} \times_3 \mathbf{M}^{-1}$
where $\mathbf{1} \in \Rbb^{1\times 1\times n_3}$ is the constant tube of all ones.
\end{definition}

\begin{definition}[$\starM$-identity]\label{def:starMidentity} 
The $\starM$-identity tensor $\TI \in \Rbb^{m\times m \times n_{3}}$ is an f-diagonal tensor with identity tubes along its main diagonal; that is, $\TI_{i,i,:} = \bfe$ for $i=1,\dots, m$ and all other tubes are equal to zero. 
\end{definition}

Combining the transpose operation and the identity tensor, we define orthogonality under the $\starM$-product as follows:
\begin{definition}[$\starM$-orthogonal]\label{def:starMorthogonality}  
A tensor $\TQ\in \Rbb^{m\times m \times n_3}$ is $\starM$-orthogonal if $\TQ^\top \starM \TQ = \TQ \starM \TQ^\top = \TI$.  
\end{definition}

\begin{remark}
We will often omit the ``$\starM$-'' prefix when describing tensor properties. 
\end{remark}

\subsection{$t$-SVDM and Eckart-Young-Like Theorem}
\label{sec:tsvdm}

Tensor singular value decompositions are central to multiway data analysis, especially for dimension reduction and feature extraction. Unique among tensor SVDs, the $t$-SVDM satisfies an Eckart-Young-like property, yielding \emph{provably optimal} low-rank tensor approximations.  
\begin{definition}[$t$-SVDM]\label{def:tsvdm}
Any tensor $\TA\in \Rbb^{n_1\times n_2\times n_3}$ can be decomposed as
	\begin{align}
	\TA = \TU \starM \TS \starM \TV^\top = \sum_{i=1}^r \TU_{:,i,:} \starM \TS_{i,i,:} \starM \TV_{:,i,:}^\top
	\end{align}
where $\TU\in \Rbb^{n_1\times n_1 \times n_3}$ and $\TV\in \Rbb^{n_2\times n_2\times n_3}$ are $\starM$-orthogonal and $\TS\in \Rbb^{n_1\times n_2\times n_3}$ is f-diagonal with ordered singular tubes 
	\begin{align}
	\|\TS_{1,1:}\|_F \ge \|\TS_{2,2,:}\|_F \ge \cdots \ge \|\TS_{r,r,:}\|_F \ge 0.
	\end{align}
\end{definition}

\begin{figure}
\centering
\includegraphics[width=0.9\linewidth]{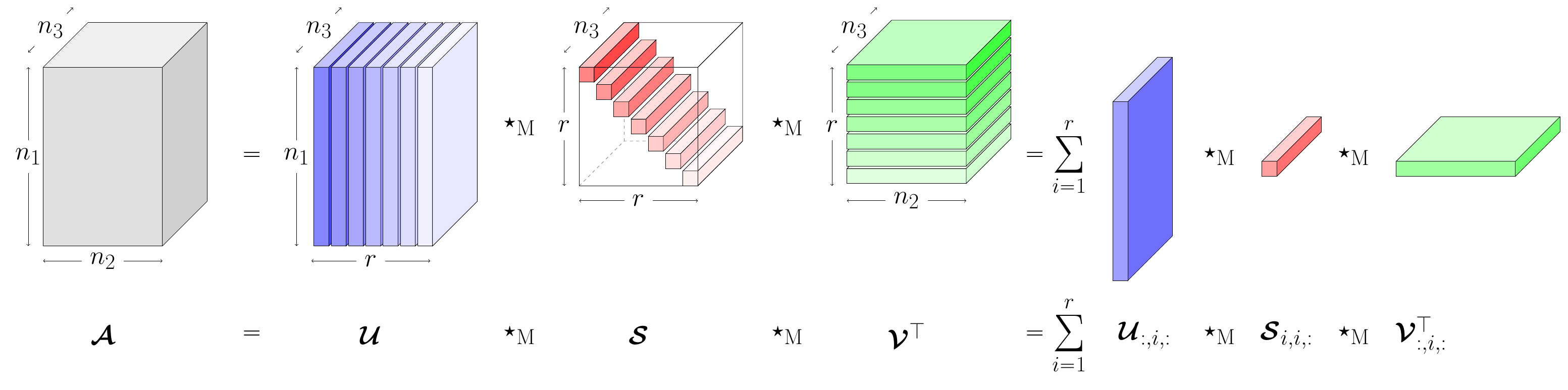}
\caption{Illustration of $t$-SVDM of $\TA$ with $\Trank(\TA,\bfM) = r$.}
\label{fig:tsvdm}
\end{figure}

We illustrate the $t$-SVDM in~\Cref{fig:tsvdm}.  
The notion of the rank of a tensor consistent with the notion of the rank of a matrix arises naturally from the $t$-SVDM. 

\begin{definition}[$t$-rank]\label{def:trank} 
Suppose we have the $t$-SVDM of a tensor $\TA = \TU \starM \TS \starM \TV^\top$.  
The $t$-rank of $\TA$ is the number of nonzero singular tubes; that is, 
	\begin{align}
	\Trank(\TA, \bfM) = \#\{\TS_{i,i,:} \mid  \TS_{i,i,:} \not= {\bf0} \}.
	\end{align}
where $\#$ denotes the cardinality of a set.
\end{definition}  
Note that the $t$-rank depends on $\bfM$, and it is possible to obtain different ranks for different transformations; e.g., see~\cite[Example 3.1.8]{Newman2019:thesis}. 

From the $t$-SVDM, we state the cornerstone of the $\starM$-framework.
\begin{mytheo}{$\starM$-Eckart-Young~\cite{Kilmer2021:pnas}}{eckartYoung}
Let $\TA\in \Rbb^{n_1\times n_2\times n_3}$  and let $\bfM$ be a nonzero multiple of an orthogonal matrix.  
Suppose $\TA$ has $t$-rank-$r$ and a $t$-SVDM $\TA = \TU \starM \TS \starM \TV^\top$. 
Given $k \le r$, an optimal $t$-rank-$k$ approximation to $\TA$ is given by the truncated $t$-SVDM; i.e., 
	\begin{align}
	\TA_k := \TU_{:,1:k,:} \starM \TS_{1:k,1:k,:} \starM \TV_{:,1:k,:}^\top \in \argmin_{\TX \in \Xcal} \|\TA - \TX\|_F
	\end{align}
where $\Xcal = \{\TX \in \Rbb^{n_1\times n_2\times n_3} \mid \Trank(\TX, \bfM) \le k\}$.
\end{mytheo}
A similar theorem can be derived for complex-valued tensors and transformations. We can further prove that the $t$-rank-$k$ approximation is superior to the matrix rank-$k$ approximation for appropriate data representations.  Details can be found in~\cite{Kilmer2021:pnas}.

This Eckart-Young-like result is unique to the $\starM$-framework -- other tensor representations can only prove truncated approximations to be quasi-optimal under their respective decompositions~\cite{KoldaBader2009}.  Learning the transformation $\bfM$ will rely heavily on the optimality guarantees ensured by matrix mimeticity.

\subsection{The Underlying Tensor Algebra}
\label{sec:underlyingAlgebra}

The mechanics of the $\starM$-product are powered by an underlying algebraic structure of tubal multiplication.  The term ``algebra'' refers to square $n_3\times n_3$ matrices over $\Rbb$ under the familiar bilinear operation of matrix multiplication.  As presented in~\cite{KernfeldKilmer2015}, the connection between tensor-tensor products and matrix algebra can be seen when we express tubal multiplication (\Cref{def:starMtubal}) equivalently in vectorized form as\footnote{It is sometimes convenient to consider the action of $\bfb$ on $\bfa$ by $\myVec(\bfa)^\top \bfR_{\bfM}[\bfb]^\top$.}
	\begin{align}\label{eqn:starMtubalVectorized}
	 \bfa \starM \bfb \equiv  \bfR_{\bfM}[\bfa] \myVec(\bfb) \qquad \text{where} \qquad \bfR_{\bfM}[\bfa] := \bfM^{-1} \diag(\bfM \myVec(\bfa)) \bfM.
	\end{align}
In words, tubal multiplication of $\bfa$ acting on $\bfb$ is equivalent to the action of the structured matrix $\bfR_{\bfM}[\bfa]$, which is parameterized by the tube $\bfa$.  The set of all matrices with the structure $\bfR_{\bfM}[\bfa]$ forms a matrix subalgebra, and the term ``tensor algebra'' and variants refers to this structure.  

Importantly, the tensor algebra is completely defined by the choice of $\bfM$.   For example, if $\bfM = \bfF$,  the (unnormalized) discrete Fourier transform matrix, then $\bfR_{\bfF}[\bfa]$ is a circulant matrix whose entries come from $\bfa$ and the resulting algebra is the \emph{algebra of circulants}~\cite{GleichEtAl2013:algebraOfCirculants}.  Comparatively, if $\bfM = \bfI$, the identity matrix, then $\bfR_{\bfI}[\bfa]$ is a diagonal matrix whose entries come from $\bfa$ and the resulting bilinear operation is given by the Hadamard pointwise product. For concreteness, we provide a few expressible tensor algebras in~\Cref{sec:underlyingStructure}.

\subsection{Invariants of Tensor-Tensor Products}
\label{sec:invariantProduct}

Before we learn data-dependent $\starM$-algebras, we present original theory about invariants of the $\starM$-product under modifications of the transformation matrix. We first show how the $\starM$-product behaves under negation of the transformation (\Cref{lem:starMInvariance}). We then present our main theorem that proves the $\starM$-product under left-matrix multiplication of the transformation is invariant only under permutations (\Cref{thm:starMInvariance}). This theory will be crucial for subsequent analysis of our algorithm. 

\begin{mylemma}{$\starM$-product under negation}{starMInvariance}
Given $\TY_1 \in \Rbb^{n_1\times p \times n_3}$ and $\TY_2\in \Rbb^{p\times n_2\times n_3}$, we have $\TY_1 \star_{-\bfM} \TY_2 = -\TY_1 \starM \TY_2$. 
\end{mylemma}

\begin{proof}
By~\Cref{defn:starMProduct}, $\bfM$ appears three times in the $\starM$-product. Thus,
	\begin{align}
	\TY_1 \star_{-\bfM} \TY_2
	= [(\TY_1\times_3 (-\bfM)) \smalltriangleup (\TY_2 \times_3 (-\bfM))] \times_3 (-\bfM)^{-1} 
	= -\TY_1 \starM \TY_2.
	\end{align}
By induction, the formula generalizes to multiplying $q$ tensors by
	\begin{align}
	\TY_1 \star_{-\bfM} \TY_2 \star_{-\bfM} \cdots \star_{-\bfM} \TY_q = (-1)^{q - 1} [\TY_1 \starM \TY_2 \starM \cdots \starM \TY_q].
	\end{align}
\end{proof}

\begin{mytheo}{$\starM$-product invariance}{starMInvariance}
The $\starM$-product and $\star_{\bfP\bfM}$-product are equivalent if and only if $\bfP$ is a permutation matrix. 
\end{mytheo}

\begin{proof}  We will show the bidirectional proof holds for tubal products. 
Tensor-tensor products will follow from~\Cref{defn:starMProduct}.   

$(\Longrightarrow)$  
Let $\bfa, \bfb\in \Rbb^{1\times 1\times n_3}$ be arbitrary. By definition, if $\widehat{\bfc} = \bfc \times_3 \bfM$, we have
	\begin{align}\label{eq:productOfMatricesMprod}
	\bfa \star_{\bfP \bfM} \bfb = (\widehat{\bfa} \star_{\bfP} \widehat{\bfb} ) \times_3 \bfM^{-1}.
	\end{align}
This will equal $\bfa \starM \bfb$ if the $\star_{\bfP}$-product is equivalent to the Hadamard pointwise product.  Thus, we need to find $\bfP$ such that for all tubes $\bfy, \bfz\in \Rbb^{1\times 1\times n_3}$, we have
	\begin{align}\label{eq:pullOutM}
		\bfy \star_{\bfP} \bfz = \bfy \odot \bfz \qquad  \xrightarrow{\texttt{vectorize}} 
		 \qquad \bfR_{\bfP}[\bfy] \bfz = \diag(\bfy) \bfz
	\end{align}
where the equivalent vectorized form redefines $\bfy \gets \myVec(\bfy)$ and $\bfz \gets \myVec(\bfz)$ and $\bfR_{\bfP}[\bfy]$ is defined in~\eqref{eqn:starMtubalVectorized}. The equality must hold for all vectors $\bfz\in \Rbb^{n_3\times 1}$, including the $j$-th standard basis vector  $\bfe_j\in \Rbb^{n_3\times 1}$. 
Thus, we have
	\begin{align}
	\begin{split}
	  \bfR_{\bfP}[\bfy] \bfe_j &= \diag(\bfy) \bfe_j\\
	\Longrightarrow \qquad  \bfP^{-1} \diag(\bfP\bfe_j) \bfP \bfy &=\bfe_j \bfe_j^\top \bfy\\
	\Longrightarrow \qquad  \phantom{\bfP^{-1}\bfP \bfy} \diag(\bfP\bfe_j) &= (\bfP\bfe_j) (\bfP^{-\top} \bfe_j)^\top.
	 \end{split}
	\end{align}
The left-hand side is a nonzero diagonal matrix and the right-hand side is a rank-one matrix.  To ensure the off-diagonal entries of the right-hand side are zero, we must have $\bfP\bfe_j$ and $\bfP^{-\top} \bfe_j$ be (nonzero) multiples of the same standard basis vector; that is, $\bfP\bfe_j = c_k\bfe_k$ and $\bfP^{-\top} \bfe_j = d_{k} \bfe_{k}$ for some $k\in \{1,\dots, n_3\}$ with $c_k,d_k \not=0$.  As a result, the $k$-th diagonal entry is the only nonzero entry. Equating the $k$-th diagonal entries on the left- and right-hand sides, we get $c_k = c_k d_k$. It follows that $d_k = 1$. Moreover, by invertibility, we have $\bfP \bfP^{-1} = \sum_{k=1}^{n_3} c_k d_k \bfe_k\bfe_k^\top = \bfI_{n_3}$.  Thus, $c_kd_k = 1$ for $k =1,\dots,n_3$. Because $d_k = 1$, we also have $c_k = 1$. In sum, $\bfP$ is invertible with standard basis vectors columns, and thus is a permutation matrix.

$(\Longleftarrow)$ Let $\bfP\in \Rbb^{n_3\times n_3}$ be a permutation matrix. 
For any vector $\bfx\in \Rbb^{n_3}$, we have 
	\begin{align}
	\bfR_{\bfI}[\bfx] = \diag(\myVec(\bfx)) = \bfP^\top \diag(\bfP \myVec(\bfx)) \bfP = \bfR_{\bfP}[\bfx]. 
	\end{align}
Because the matrices are equal, the $\star_{\bfP}$-product is equivalent to the pointwise product.  
Using~\eqref{eq:pullOutM} for general $\bfM$, we get $\bfR_{\bfM}[\bfx] = \bfR_{\bfP\bfM}[\bfx]$ for all $\bfx\in \Rbb^{1\times 1\times n_3}$.

\end{proof}

	
\section{$\starM$-Optimization}
\label{sec:starMOpt}

We introduce $\starM$-optimization, a framework to learn an optimal transformation and a desired tensor representation (e.g., low-$t$-rank) simultaneously. In~\Cref{sec:varpro}, we exploit the natural coupling between transformation and representation using variable projection, a bilevel optimization strategy that eliminates one variable using the unique optimality properties of the $\starM$-framework. Building from the theory in~\Cref{sec:invariantProduct}, we describe the uniqueness of the obtained solutions for two prototype problems. To preserve $\starM$-algebraic integrity, we present the main $\starM$-optimization algorithm  in~\Cref{sec:optimizationOnManifold} that learns orthogonal transformations via Riemannian optimization.  We derive formulas to differentiate through $\starM$-representations in~\Cref{sec:starMDerivatives}. We conclude with a discussion of the computational cost~(\Cref{sec:cost}) and convergence analysis (\Cref{sec:analysis}). 

\subsection{Variable Projection Formulation}
\label{sec:varpro}

To learn the optimal transformation, we set up the following optimization problem
    \begin{align}\label{eq:optFull}
    \min_{\bfM \in \Mcal, \TX \in \Xcal} \Phi(\bfM, \TX) 
    \end{align}
where the variables are an invertible transformation matrix $\bfM\in \Mcal \subset \Rbb^{n_3\times n_3}$  and a desired representation $\TX\in \Xcal$. Here, $\Mcal$ and $\Xcal$ are the feasible sets  for $\bfM$ and $\TX$, respectively.  The objective function $\Phi: \Mcal\times \Xcal \to \Rbb$ measures the quality of the representation and can include additional regularization and penalty terms.

To capture the coupling between representation and transformation, we form a bi-level optimization problem using variable projection~\cite{GolubPereyra1973:VarPro, OLearyRust2013:VarPro, Kaufman1975:VarPro, vanLeeuwenAravkin2021:nonsmoothVarPro}, given by
    \begin{subequations}\label{eq:bilevelOpt}
	            \begin{alignat}{3}
            & \min_{\bfM\in \Mcal} &&\quad\overline{\Phi}(\bfM) &&\equiv \Phi(\bfM, \TX(\bfM)) \label{eq:outerOpt} \\
            & \subjectto_{\phantom{\bfM\in \Mcal}} &&\quad  \TX(\bfM) &&\in \argmin_{\TX \in \Xcal} \Phi(\bfM, \TX), \label{eq:innerOpt}
            \end{alignat}
    \end{subequations}
Adopting notation from~\cite{vanLeeuwenAravkin2021:nonsmoothVarPro} and terminology from~\cite{Newman2021:GNvpro}, we call $\Phi$ the \emph{full} objective function and $\overline{\Phi}$ the \emph{reduced} objective function. Because of the optimality guarantees of the $\starM$-framework, we are able to solve the inner optimization problem~\eqref{eq:innerOpt}  to obtain an optimal tensor representation for a given transformation, $\TX(\bfM)$, which explicitly exploits the coupling.  The generality of~\eqref{eq:bilevelOpt} expands the applicability of learnable tensor algebras to problems beyond tensor completion. We present two general prototype problems that are amenable to our variable projection formulation.

\begin{myproto}{$t$-linear regression $\Phi_{\rm reg}$}{regression}
Using a model tensor $\TA\in \Rbb^{n_1\times p\times n_3}$ and observations $\TB\in \Rbb^{n_1\times n_2\times n_3}$, we form a $t$-linear regression problem via
	\begin{align}\label{eq:tLinearRegressionFull}
		\min_{\bfM \in \Mcal, \TX \in\Xcal}  \Phi_{\rm reg}(\bfM, \TX) \equiv \tfrac{1}{2}\|\TA \starM \TX - \TB\|_F^2.
	\end{align}  
where $ \Xcal = \Rbb^{p\times n_2\times n_3}$.  This formulation admits an analytic solution to~\eqref{eq:innerOpt}, obtained by solving the $\starM$-analog of the normal equations
	\begin{align}
	(\TA^\top \starM \TA) \starM \TX_{\rm reg}(\bfM) = \TA^\top \starM \TB.
	\end{align}
\end{myproto}

\begin{myproto}{Optimal low-$t$-rank approximations $\Phi_{\rm low}$}{lowRank}
Given $\TA \in \Rbb^{n_1\times n_2\times n_3}$, we find an optimal $t$-rank-$k$ approximation by solving
	\begin{align}\label{eq:LowRank}
	\min_{\bfM\in \Mcal, \TX\in \Xcal} \Phi_{\rm low}(\bfM, \TX) \equiv \tfrac{1}{2}\|\TA - \TX\|_F^2
	\end{align}
where $\Xcal = \{\TX \in \Rbb^{n_1\times n_2\times n_3}\mid \text{$t$-rank}(\TX, \bfM) \le k\}$.    Despite the nonconvexity of~\eqref{eq:LowRank}~\cite{Recht2010:nuclearNorm}, the truncated $t$-SVDM solves~\eqref{eq:innerOpt} (\Cref{thm:eckartYoung}); i.e., 
	\begin{align}
	\TX_{\rm low}(\bfM) = \TA_k \equiv \TU_{:,1:k,:} \starM \TS_{1:k,1:k,:} \starM \TV_{:,1:k,:}^\top.
	\end{align}

\end{myproto}

While we do not present an exhaustive list of all the possible variations, the prototype problems are central to a wide range of applications, including multiframe blind deconvolution~\cite{FAN20122112}, which sets up a least-squares-like system, and tensor nuclear norm minimization, which relies on the $t$-SVDM~\cite{ZhangEtAl2014:tensorNuclearNorm, KongLuLin2021:TensorQRank, Recht2010:nuclearNorm}.

\subsection{Uniqueness and Invariance of the Prototype Problems}
\label{sec:uniqueInvariancePrototypeShort}
Each prototype problem offers uniqueness properties of the representation $\TX(\bfM)$ and invariance to modifications of $\bfM$.   We briefly mention the properties here and leave additional details in~\Cref{sec:objectiveInvariance} and~\Cref{sec:uniqueInvariancePrototype}. 

The uniqueness of the $t$-linear regression solution follows from similar properties of the uniqueness of least squares solution in the matrix case. Specifically, if the data matrix $\bfA$ has full column rank, the solution will be unique; a similar analogy holds for the tensor case. The reduced $t$-linear regression objective function is invariant to both permutation and negation of the rows of $\bfM$. Thus, if $\bfM^*$ is an $n_3\times n_3$ optimal transformation, there are at least $2^{n_3}(n_3!)$ equally optimal transformations.

The $t$-SVDM is unique up to $\starM$-analogies of equivalent matrix SVD properties, specifically that orthogonal transformations of singular basis vectors corresponding to repeated singular values)~\cite{GoluVanl96, Demmel1997}.  The reduced low-$t$-rank objective function is invariant to permutation of the rows and negation of $\bfM$. Thus, if $\bfM^*$ is an $n_3\times n_3$ optimal transformation, there are at least $2(n_3!)$ equally optimal transformations. 

\subsection{$\starM$-Optimization Algorithm: Riemannian Gradient Descent}
\label{sec:optimizationOnManifold}
\label{alg:starMopt}

The algebraic structure in~\eqref{eqn:starMtubalVectorized} shows that $\bfM$ and its inverse simultaneously diagonalize the tubal product and affect the ``eigenvalues'' in the transformed space. This makes the optimization problem for $\bfM$ nonconvex and thereby challenging to solve. We will use a gradient-based, iterative algorithm to solve for optimal transformations. To ensure invertibility of $\bfM$, we will let the feasible set be the orthogonal group
\begin{align}\label{eq:MSet}
	\Mcal = \Ocal_{n_3} := \left\{\bfQ \in \Rbb^{n_3\times n_3} \mid \bfQ^\top \bfQ = \bfI_{n_3}\right\}.
\end{align}
\begin{remark}
We will use $\Mcal$ and $\Ocal_{n_3}$ interchangeably for the remainder of the paper. 
\end{remark}

We provide a brief overview of $\starM$-optimization; for more intuition and details, see~\Cref{sec:riemannian_optimization_background} and the excellent resources~\cite{Absil2008,  boumal2023intromanifolds, EdelmanAriasSmith1998}. To learn the optimal transformation under orthogonality constraints, we leverage techniques from Riemannian optimization, which consists of three steps. First, we compute the  Euclidean gradient
	\begin{align}\label{eq:negativeGradient}
	\bfG = \nabla \overline{\Phi}(\bfM).
	\end{align}
Second, we compute the Riemannian gradient by projecting onto a tangent space via 
	\begin{align}\label{eqn:projection}
	\myGrad \overline{\Phi}(\bfM) := \bfM\bfOmega  \qquad \text{where} \qquad \bfOmega = \frac{\bfM^\top \bfG - \bfG^\top \bfM}{2}.
	\end{align}
Third, we return to the manifold via a retraction, specifically the exponential mapping
	\begin{align}\label{eq:exponentialMapping2}
	\bfM \gets \retraction_{\bfM}(-\alpha \myGrad \overline{\Phi}(\bfM))  \qquad \text{with} \qquad \retraction_{\bfM}(\bfM\bfOmega) = \bfM \exp(\bfOmega)
	\end{align}
where $\exp: \Rbb^{n_3\times n_3} \to \Rbb^{n_3\times n_3}$ is a matrix function~\cite{Higham2008:matrixFunctions}. The step size $\alpha>0$ is either fixed or chosen from a backtracking line search~\cite[Alg. 4.2]{boumal2023intromanifolds}. We use standard stopping criteria based on the norm of the Riemannian or Euclidean gradient being below a user-defined tolerance. We set the tolerance to be $10^{-10}$ in our experiments and include a maximum number of iterations as an additional stopping condition. 
	
\subsection{Computing $\starM$-Derivatives}
\label{sec:starMDerivatives}

As the name implies, Riemannian gradient descent relies on computing Euclidean derivatives of the objective function.    
We describe the key gradients for our prototype problems, which rely on the orthogonality of the transformations. 
Derivatives for a general invertible $\bfM$ can be found in~\cite{Newman2019:thesis, NewmanTNN2018}. 
For additional details, including the derivative of the $\starM$-product, see~\Cref{sec:fundamental_derivatives}. 

\subsubsection{Derivatives for $t$-Linear Regression (\Cref{proto:lowRank})}  
\label{sec:regressionDerivatives}

Recall, $\TX_{\rm reg}(\bfM)$ solves the inner optimization~\eqref{eq:innerOpt} using the full $t$-linear regression objective function $\Phi_{\rm reg}$ given by~\eqref{eq:tLinearRegressionFull}.  The derivative of the corresponding reduced objective function, $\overline{\Phi}_{\rm reg}$, is given by
	\begin{align}\label{eq:tLinearRegressionGrad}
	\begin{split}
	\nabla \overline{\Phi}_{\rm reg}(\bfM) 
		&=\nabla_{\bfM} \Phi_{\rm reg}(\bfM, \TX_{\rm reg}(\bfM))\\
		&=\nabla\widetilde{\TB}(\bfM) \TR(\bfM) 
		+ {\nabla \TX_{\rm reg}(\bfM) }\nabla_{\TX} \Phi_{\rm reg}(\bfM, \TX_{\rm reg}(\bfM))
	\end{split}
	\end{align}
where $\widetilde{\TB}(\bfM) = \TA\starM \TX_{\rm reg}(\bfM)$ is the approximation and $\TR(\bfM) = \widetilde{\TB}(\bfM) - \TB$ is the residual. We compute the first term in~\eqref{eq:tLinearRegressionGrad} using~\Cref{lem:starMGrad} while holding $\TX_{\rm reg}(\bfM)$ constant. For the second term, because $\TX_{\rm reg}(\bfM)$ is a minimizer of the inner optimization problem, it satisfies the first-order optimality conditions
	\begin{align}\label{eq:optimalityCondition}
	\nabla_{\TX} \Phi_{\rm reg}(\bfM, \TX_{\rm reg}(\bfM))= {\bf 0}. 
	\end{align}
Thus, the second term in~\eqref{eq:tLinearRegressionGrad} can be eliminated, which, in addition to simplifying the derivative formula, avoids non-trivial differentiation through $\TX_{\rm reg}(\bfM)$.   We emphasize that~\eqref{eq:optimalityCondition} is unique to the $\starM$-framework.  Such optimality conditions are not guaranteed to hold for other tensor decomposition techniques.

\subsubsection{$t$-SVDM Derivatives (\Cref{proto:regression})} 
\label{sec:tsvdmDerivatives}
We differentiate through the reduced low-$t$-rank objective function, $\overline{\Phi}_{\rm low}$, with respect to $\bfM$.    Similar to the derivation in~\Cref{sec:regressionDerivatives}, we have
	\begin{align}
	\nabla \overline{\Phi}_{\rm low}(\bfM) ={\nabla \TX_{\rm low}(\bfM)}[\nabla_{\TX} \Phi_{\rm low}(\bfM, \TX_{\rm low}(\bfM))]
	\end{align}
where $\nabla_{\TX} \Phi_{\rm low}(\bfM, \TX_{\rm low}(\bfM)) = -(\TA - \TX_{\rm low}(\bfM))$. Unlike the $t$-linear regression case, the solution $\TX_{\rm low}(\bfM)$ does not satisfy first-order optimality conditions of $\Phi_{\rm low}$ and may not vanish. Hence,  we must differentiate through $\TX_{\rm low}(\bfM)$.

The optimal $t$-rank-$k$ approximation to $\TA$ is the truncated $t$-SVDM given by $\TX_{\rm low}(\bfM) = \TU_{:,1:k,:}\starM \TS_{1:k,1:k,:} \starM \TV_{:,1:k,:}^\top$.  Here, the factors implicitly depend on $\bfM$ and hence we differentiate through each factor.  Notationally, we consider the gradient $\nabla \TX_{\rm low}(\bfM)[\cdot]: \Rbb^{n_1\times n_2\times n_3} \to \Rbb^{n_3\times n_3}$ to be an operator defined as follows:
	\begin{equation}\label{eq:tsvdmFactorDerivatives}
	\begin{alignedat}{3}
	\nabla \TX_{\rm low}(\bfM)[\TR] &=\qquad &&\nabla \TU_{:,1:k,:}(\bfM)[\TR \starM (\TS_{1:k,1:k,:} \starM \TV_{:,1:k,:}^\top )^\top]  \\
	& \qquad + &&\nabla \TS_{1:k,1:k,:} (\bfM) [\TU_{:,1:k,:}^\top \starM \TR \starM \TV_{:,1:k,:}] \\
	& \qquad + &&\nabla \TV_{:,1:k,:}(\bfM)[(\TU_{:,1:k,:} \starM \TS_{1:k,1:k,:})^\top \starM \TR]
	\end{alignedat}
	\end{equation}
where $\TR = \nabla_{\TX} \Phi_{\rm low}(\bfM, \TX_{\rm low}(\bfM))$. The formulas for $t$-SVDM derivatives from similar formulas for the matrix SVD; see~\Cref{app:svdDerivation} for details.

\subsection{Computational Cost} 
\label{sec:cost}
For the prototype problems, the dominant cost at each $\starM$-optimization iteration is a least squares solve or a matrix SVD per frontal slice, which, if $n_1=n_2=n$, costs $\Ocal(n^3 n_3)$ floating point operations.  This is computationally demanding for large-scale problems. These operations can be performed in parallel, decreasing communication cost.
There are additional costs for applying the transformation along the third dimension and performing a line search. The former is the dominant communication cost in the algorithm. Additionally, the learned orthogonal transformation is often dense, requiring a storage cost of $\Ocal(n^3)$ floating point numbers. In practice, we consider $\starM$-optimization to be a one-time, offline cost, and focus on the benefits of using a learned transformation in terms of representation quality, compressibility, and transferability in subsequent data analysis pipelines.

\subsection{Convergence Analysis of $\starM$-Optimization}
\label{sec:analysis}

As described in~\cite[p. 57-58]{boumal2023intromanifolds},  Riemannian gradient descent will converge if two standard assumptions are satisfied: the objective function is bounded below and the algorithm ensures sufficient decreases at each iteration.  Both prototype problems described in~\Cref{sec:varpro} are bounded below by zero, thereby satisfying the first assumption. The backtracking linesearch for $\starM$-optimization (\Cref{alg:starMopt}) enforces the second assumption if the objective function is sufficiently regular. A sufficient criterion for regularity is that the gradient of the objective function is Lipschitz continuous\footnote{To define  Lipschitz continuity precisely on Riemannian manifolds, we need tools such as parallel transport because $\myGrad \overline{\Phi}(\bfM)$ and $\myGrad \overline{\Phi}(\bfM')$ lie in different tangent spaces.  Our proof will not rely on this subtlety, hence we omit these details; see~\cite{BoumalAbsilCoralia2018:globalManifold} for further reading.}; that is, there exists some constant $L> 0$ such that for all $\bfM, \bfM' \in \Mcal$, we have
	\begin{align}\label{eq:LipschitzContinuity}
	\|\myGrad \overline{\Phi}(\bfM) - \myGrad \overline{\Phi}(\bfM')\|_F \le L \|\bfM - \bfM'\|_F.
	\end{align}
	
We analyze $\starM$-optimization for the $t$-linear regression (\Cref{proto:regression}).  We will show that the reduced objective function, $\overline{\Phi}_{\rm reg}$, has a Lipschitz continuous gradient by proving that Riemannian Hessian, $\myHess \overline{\Phi}$, has a bounded operator norm for all $\bfM\in \Mcal$. We will then discuss of the convergence behavior implied by this Lipschitz continuity and illustrate the behavior through empirical examples.

We introduce two new definitions to assist proving the boundedness Riemannian Hessian, the  $\starM$-pseudoinverse and $\starM$-operator norm, both relying on the $t$-SVDM.
\begin{definition}[$\starM$-psuedoinverse]\label{def:starMpseudo} 
Given a $t$-rank-$r$ tensor $\TA \in \Rbb^{n_1 \times n_2 \times n_{3}}$ and $t$-SVDM $\TA = \TU \starM \TS \starM \TV^\top$, the $\starM$-pseudoinverse, $\TA^\dagger$, is
	\begin{align}
	\TA^{\dagger} = \TV \starM \TS^\dagger \starM \TU^\top 
	 \qquad \text{where} \qquad 
	\widehat{(\TS^\dagger)}_{:,:,i}= (\widehat{\TS}_{:,:,i})^\dagger \text{ for $i = 1,\dots, n_3$.}
	\end{align}
 \end{definition}

 \begin{definition}[$\starM$-operator norm]\label{def:operatorNorm}
For $\TA \in \Rbb^{n_1\times n_2\times n_3}$, the $\starM$-operator norm is the globally largest singular value in the transform domain\footnote{We note that this definition can be generalized to any operator norm.  We chose the $2$-norm because of the connection to the singular values.}; that is,
	\begin{align}
		\|\TA\| = \max_{i\in \{1,\dots,n_3\}} \sigma_{1}(\widehat{\TA}_{:,:,i})
	\end{align}
where $\sigma_j(\bfZ)$ returns the $j$-th largest singular value of the matrix $\bfZ$. 
\end{definition}

Like the matrix equivalent, the Frobenius norm is an upper bound for the operator norm; i.e., $\|\TA\| \le \|\TA\|_F$. Further note that the $\starM$-operator norm of the pseudoinverse returns the inverse of the smallest, nonzero singular value; that is, 
	\begin{align}
	\|\TA^{\dagger}\| = \max_{i\in \{1,\dots n_3\}} \{\sigma_{r_i}(\widehat{\TA}_{:,:,i})^{-1} \mid r_i = \rank(\widehat{\TA}_{:,:,i})\}.
	\end{align}
We include examples of $\starM$-operator norms of tensors in~\Cref{sec:boundedOperatorExamples}.   We now have the tools to prove the Lipschitz continuity of the Riemannian gradient of the reduced $t$-linear regression function. 

\begin{mytheo}{Lipschitz continuity of $\myGrad \overline{\Phi}_{\rm reg}$}{lipschitzContinuity}
Let  $\TA$ be a data tensor and assume its $\starM$-pseudoinverse  has a bounded $\starM$-operator norm.  Then, $\myGrad \overline{\Phi}_{\rm reg}$ is Lipschitz continuous. 
\end{mytheo}

\begin{proof} 
By~\cite[Corollary 5.47]{boumal2023intromanifolds}, if the $\nabla^2 \overline{\Phi}_{\rm reg}$ is bounded in the Frobenius norm, then $\myHess \overline{\Phi}_{\rm reg}$ will be bounded. This is because the canonical norm on $\Ocal_{n_3}$ is the Frobenius norm~\cite{EdelmanAriasSmith1998} and the Riemannian Hessian is a projection of the Euclidean Hessian.  By~\cite[Corollary 10.47]{boumal2023intromanifolds}, a bounded Riemannan Hessian implies the Riemannian gradient is Lipschitz continuous.  

Proving the boundedness of the Euclidean Hessian relies on vectorizing the $\starM$-product with respect to $\bfM$, differentiating using the product rule, and eliminating the dependence on $\bfM$ in the Frobenius norm using the orthogonality of the transformation. Because the $\starM$-product uses $\bfM$ three times (\Cref{defn:starMProduct}), we have to vectorize in three different ways. For the approximation, $\TA \starM \TX_{\rm reg}(\bfM)$, the solution also depends on $\bfM$, thus requiring another vectorization step. In addition, we must differentiate through the solution, $\TX_{\rm reg}(\bfM)$, which we do using implicit differentiation. For clarity here, we leave the technical details to~\Cref{sec:hessianBound}. 

\end{proof}

\begin{mycorollary}{Convergence of $\starM$-optimization for $\overline{\Phi}_{\rm reg}$}{convergence}
Under the same assumptions as~\Cref{thm:lipschitzContinuity}, $\starM$-optimization will return an iterate $\bfM$ with $\|\myGrad \overline{\Phi}_{\rm reg}(\bfM)\|_F \le \varepsilon$ in at most $\Ocal(1/\varepsilon^2)$ iterations. 
\end{mycorollary}

\begin{proof}
Convergence of $\starM$-optimization is guaranteed by the Lipschitz continuity of the Riemannian gradient (\Cref{thm:lipschitzContinuity}). The convergence rate follows directly from~\cite[Corollary 2.9 and Theorem 2.11]{BoumalAbsilCoralia2018:globalManifold}. 
\end{proof}

We note that~\cite[Sections 4.3 and 4.4]{boumal2023intromanifolds} require less rigorous assumptions for convergence. In our case, proving the stronger assumption of Lipschitz continuity was more straightforward.  We show the convergence properties in the following examples. 

\begin{example}[Convergence of $\starM$-optimization]\label{exam:nonconvergence}
Consider tensors $\TA\in \Rbb^{3\times 1\times 2}$ and $\TB\in \Rbb^{3\times 1\times 2}$ given by
	\begin{align}
	\TA_{:,:,1} = \begin{bmatrix}
	1 & 0 \\ 0 & 1 \\ 0 & 0
	\end{bmatrix}, \quad  
	\TA_{:,:,2} = \begin{bmatrix}
	0 & 0 \\ 1 & 0 \\ 0 & 1
	\end{bmatrix}, \quad
	\TB_{:,1,1} = \begin{bmatrix}
	1 \\ 1 \\ 1
	\end{bmatrix},  \text{and} \quad
	\TB_{:,1,2} = \begin{bmatrix}
	1 \\ 1 \\ 1 
	\end{bmatrix}.
	\end{align}
We can parameterize any $2\times 2$ rotation matrix $\bfQ(\bftheta)$ by angle $\theta\in [0,2\pi)$. The univariate reduced objective function in terms of $\theta$ is given by
	\begin{align}\label{eq:angleObjFctn}
	\bar{\phi}(\theta) \equiv \tfrac{1}{2}\|\TA \star_{\bfQ(\theta)} \TX(\theta) - \TB\|_F^2 = 3 - \frac{16}{7+\cos(4\theta)}.
	\end{align}
To avoid notational confusion, we use $\bar{\phi}$ to denote the  angle-based objective function.As expected, $\bar{\phi}$ is twice continuously differentiable with bounded derivatives.  

We show the convergence behavior in~\Cref{fig:nonconvergence} and observe multiple phenomena. First, the function yields four local minima at $\theta = \frac{n\pi}{4}$ for $n=1,3,5$, and $7$ and four local maxima at $\theta = \frac{n \pi}{4}$ for $n=0,2,4,6$. The four different options correspond to all possible invariants of $\bar{\phi}$ for $2\times 2$ orthogonal matrices (\Cref{sec:uniqueInvariancePrototype}). Second, we observe that we always converge to a local optimum. 
This empirically supports the proof that the Riemannian gradient is Lipschitz continuous (\Cref{thm:lipschitzContinuity}).

\begin{figure}
\centering

	\includegraphics[width=\linewidth]{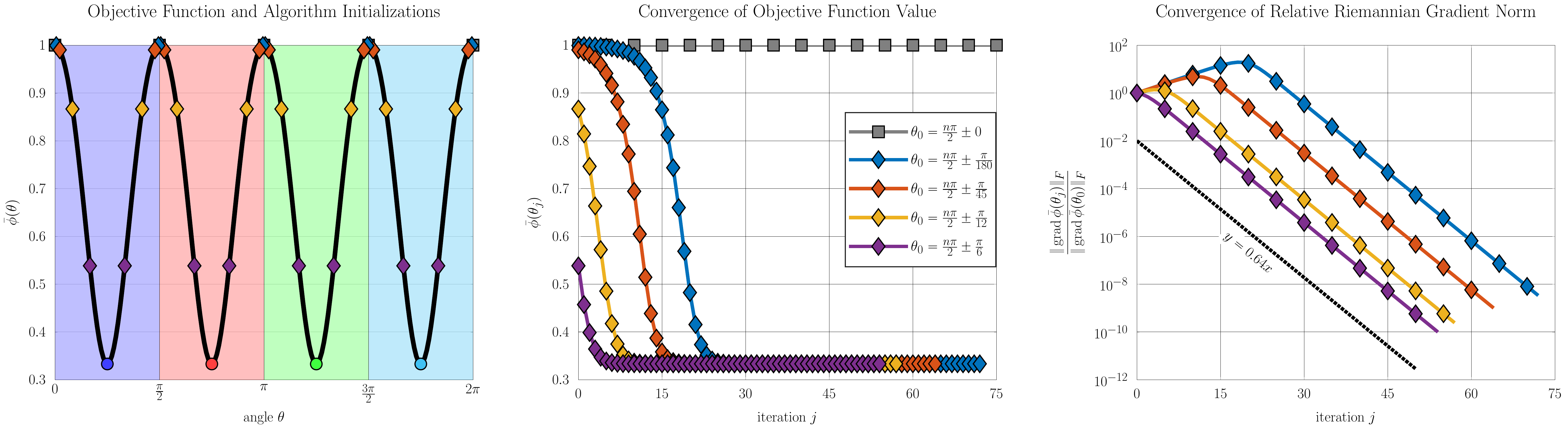}

	\caption{Convergence of $\starM$-optimization for~\eqref{eq:angleObjFctn}.  
(Left): Visualization of $\overline{\phi}$ with initial angles indicated (diamonds). There are four minima (colorful circles) and four maxima ({\bf \color{gray} gray} squares). (Middle): Convergence of $\overline{\phi}$ using a fixed step size of $\alpha=10^{-1}$. Here, $\alpha = \frac{1}{L^*}$ where $L^* = 10 \ge \max_{\theta} |\bar{\phi}''(\theta)|$. For all initial angles (save the maxima), $\starM$-optimization converges to the closest minimum, satisfying~\Cref{cor:convergence}. The periodicity of $\overline{\phi}$ leads to the same convergence values for initializations equidistant from the optima (i.e., diamonds of the same color follow the same convergence behavior). (Right): Convergence of the norm of $\myGrad \overline{\phi}$.  The convergence rate is asymptotically linear with $\|\myGrad(\theta_j)\|_F \approx 0.64 \|\myGrad(\theta_{j-1})\|_F$. Thus, $\|\myGrad(\theta_{j})\|_F \le \varepsilon$ in $\Ocal(\log(1 / \varepsilon))$ iterations, which is within the guarantees of~\Cref{cor:convergence}. 
}
\label{fig:nonconvergence}
\end{figure}

\end{example}

We now compare convergence $\starM$-optimization to an alternating descent algorithm for various problem sizes.  

\begin{example}[$\starM$-optimization vs. alternating descent]\label{exam:varproAD} 
We construct a synthetic $t$-linear regression problem where we control the true optimal solution.  We choose $n_3$ one-dimensional linear models with slope-intercept pairs $\alpha_i = \beta_i = -1 + 2\frac{i}{n_3}$ for $i=1,\dots,n_3$. We sample $100$ points uniformly from the lines and create two tensors in the transform domain, $\widehat{\TA} \in \Rbb^{n_1\times 2\times n_3}$ and $\widehat{\TB}\in \Rbb^{n_1\times 1\times n_3}$, such that 
	\begin{align}
	\widehat{\TA}_{:,:,i} = \begin{bmatrix} 1 & z_1^{(i)} \\ 1 & z_2^{(i)} \\ \vdots & \vdots \\ 1 & z_{n_1}^{(i)} \end{bmatrix}
	\quad \text{and} \quad
	\widehat{\TB}_{:,:,i} = \begin{bmatrix} \beta_i + \alpha_i z_1^{(i)} \\ \beta_i + \alpha_i z_2^{(i)}  \\ \vdots \\ \beta_i + \alpha_i z_{n_1}^{(i)} \end{bmatrix} \quad \text{for $i=1,\dots,n_3$.}
	\end{align}
We use $n_3=2^d$ for $d=1,\dots,4$, choose the underlying transformation $\bfM_{\rm true} = \bfC$ as the $2^d\times 2^d$ discrete Cosine transform\footnote{In {\sc Matlab}, \texttt{C = dct(eye(n3));}.}, and return to the spatial domain via
	\begin{align}
	\TA = \widehat{\TA} \times_3 \bfM_{\rm true}^\top \qquad \text{and} \qquad 
	\TB = \widehat{\TB} \times_3 \bfM_{\rm true}^\top
	\end{align}
This setup ensures that the optimal value is $\overline{\Phi}_{\rm reg}(\bfM_{\rm true}) = 0$.   

We compare the variable projection formulation of $\starM$-optimization to an alternating descent approach on the full optimization problem, which switches between updating the representation and the transformation; i.e., 
	\begin{align}
	\TX  \gets \TX - \beta \nabla_{\TX} \Phi_{\rm reg}(\bfM, \TX) \quad \text{and} \quad 
	\bfM \gets \retraction_{\bfM}(-\alpha \myGrad_{\bfM} \Phi_{\rm reg}(\bfM, \TX)).
	\end{align}
We use a backtracking line search to determine the step sizes $\alpha$ and $\beta$. We train until stopping criteria are met or a line search breaks and present the results in~\Cref{fig:varproAD}. 

We see that $\starM$-optimization converges for all choices of $\bfM$, empirically supporting~\Cref{cor:convergence}. The convergence is slower for larger $n_3$, which follows from the Lipschitz constant's dependence on the dimensions of the problem (\Cref{sec:hessianBound}). In comparison, alternating descent converges more slowly than $\starM$-optimization or fails to converge in terms of both iterations and wall clock time. This behavior reflects the strong coupling between the representation and the transformation. We explore this example for noisy data and provide additional geometric intuition in~\Cref{sec:leastSquares}. 

	\begin{figure}
	\centering
	
	\includegraphics[width=0.8\linewidth]{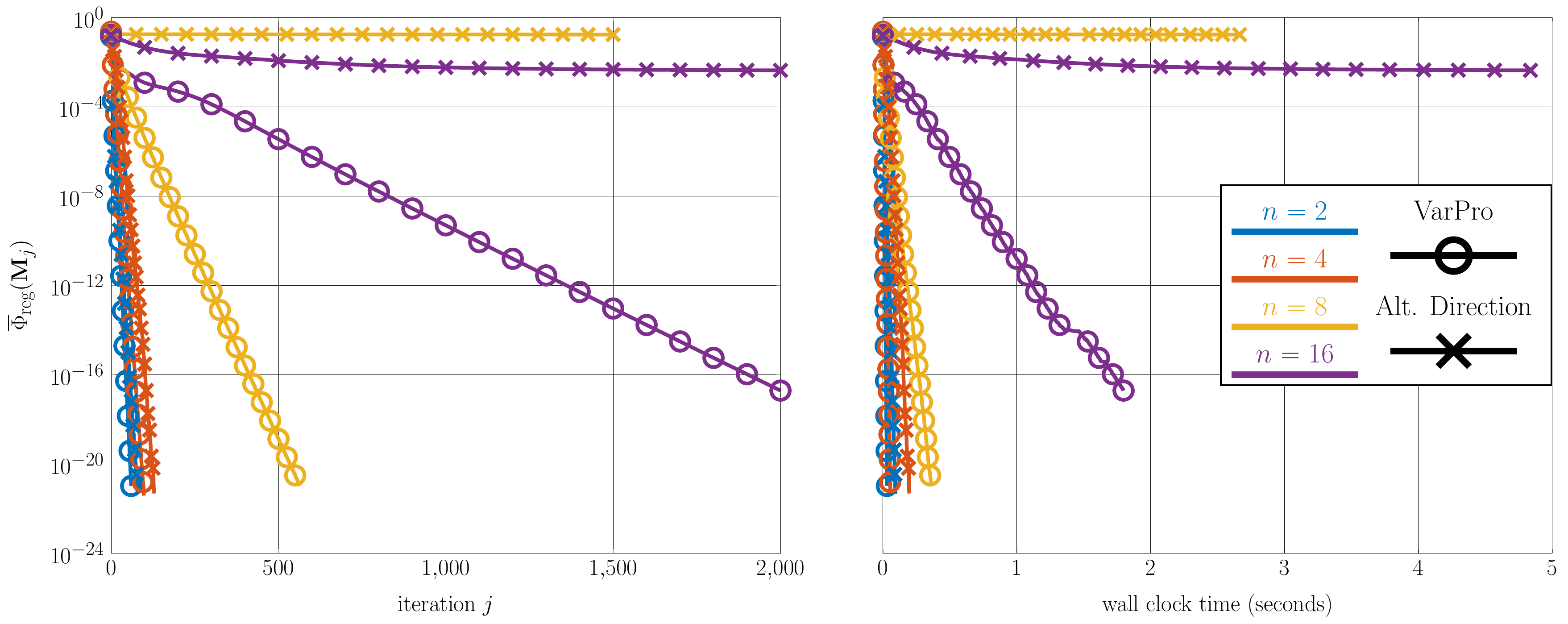}
	\caption{$\starM$-optimization (VarPro) (circles) vs. alternating descent ($\times$'s) for $t$-linear regression problems. We show the convergence of $\overline{\Phi}_{\rm reg}$ to learn optimal $n_3\times n_3$ transformations for various choices of size $n_3$ (colors)  vs. iteration (left) and wall clock time (right). Experiments were performed on a 2022 MacBook Pro laptop with an Apple M2 chip, $8$ cores, and $16$ gigabytes of memory using {\sc Matlab} R2022b.}
	\label{fig:varproAD}
	\end{figure}

\end{example}

\section{Numerical Experiments}
\label{sec:numerical}

We present several numerical examples to show the intuition behind and breadth of applications of $\starM$-optimization.  In~\Cref{sec:portfolio}, we extend to a constrained $t$-linear regression for financial index tracking. In~\Cref{sec:digits}, we compress image data using a $t$-SVDM and demonstrate the transferability of the learned transformation. In~\Cref{sec:rom}, we present a tensor reduced order modeling paradigm for parametric dynamical systems and explore the spectral implications of $\starM$-optimization. All code used to generate these results is publicly available at \url{https://github.com/elizabethnewman/star-M-opt}. 

\subsection{Experiment Setup and Parameters} 
The following experiments focus on comparing different choices of transformation for given objectives.  The heuristic transformations to which we compare are the identity matrix $\bfI$, the discrete cosine transform matrix $\bfC$, and a data-dependent matrix $\bfZ^\top$ that is the transposed left singular matrix of the mode-$3$ unfolded tensor; that is, $\TA_{(3)} = \bfZ \bfD \bfW^\top$.  In~\Cref{sec:rom}, we additionally compare to the original $t$-product, which uses the complex-valued discrete Fourier transform matrix\footnote{In practice, we use fast Fourier transforms (\texttt{AHat = fft(A,[],3);}   \texttt{A = ifft(AHat,[],3);}) instead of forming $\bfF$ explicitly.} $\bfF$ and a random orthogonal matrix\footnote{In {\sc Matlab}, \texttt{Q = orth(randn(n3));}} $\bfQ$.

\subsection{Fitting Financial Data with $t$-Linear Regression}
\label{sec:portfolio}

Indexing tracking is widely-studied in financial mathematics~\cite{Boyd_Vandenberghe_2018, Benidis2018:indexTracking}, with the goal to predict market trends using a small portfolio of securities. Mathematically, index tracking can be formulated as a linear regression problem of the form
	\begin{align}\label{eq:indexTrackingObjective}
	\min_{\bfx \in \Delta^m} \tfrac{1}{2}\|\bfA \bfx - \bfb\|_2^2
	\end{align}
where $\bfA\in \Rbb^{n_1 \times m}$ contains portfolio data over $n_1$ historic time points for $m$ securities and $\bfb \in \Rbb^{n_1}$ contains financial index data over the same time period. The entries of $\bfA$ and $\bfb$ are the percent change of return from a start date. The goal of~\eqref{eq:indexTrackingObjective} is to learn the weights of the portfolio $\bfx$ where $\bfx_i$ represents the proportion of investment in security $i$. This proportion is indicated by the constraints on $\bfx$, where $\bfx\in \Delta^m$ means $\bfx$ belongs to the unit simplex ($\sum_{i=1}^m\bfx_i = 1$ and $\bfx_i \ge 0$ for $i=1,\dots, m$).  

The matrix index tracking setup does not account for potential multidimensional relationships in the securities. To explore the benefits of capturing multilinearity, we build a tensor index tracking experiment by adding a third dimension given by sector. Instead of tracking the global S\&P 500 index, we track the equivalent index per sector.  In this setting, we can capture trends across sectors with a single set of weights. We construct an intentionally-diversified portfolio containing of ten stocks per sector ($110$ total), and minimize the index tracking objective function 
	\begin{align}\label{eq:tensorIndexTracking}
	\min_{\bfM\in \Mcal, \TX\in \Xcal_{\Delta}} \Phi_{\rm ind}(\bfM, \TX) \equiv \tfrac{1}{2}\|\TA \starM \TX- \TB\|_F^2 + \tfrac{\lambda}{2}\|\TX\|_F^2
	\end{align}
where $\TA_{i,j,k}$ contains the percent change of return on day $i$ of stock $j$ belonging to sector $k$ and  $\TB_{i,1,k}$ contains the percent change of return on day $i$ for sector index $k$.  The constraint set for $\TX$ is the tensor analog of the unit simplex given by
	\begin{align}
	\Xcal_{\Delta} = \left\{\TX\in \Rbb^{10\times 1\times 11} \; \left\vert \; \sum_{k=1}^{11} \sum_{i=1}^{10} \TX_{i,1,k}  = 1 \text{ and }\TX_{i,1,k} \ge 0\right.\right\}.
	\end{align}
The data are collected from from Yahoo!~Finance~\cite{yahoo} using code provided by~\cite{Lensky2023:yahooFinance} to download into {\sc Matlab}.  We query five months of training data (January to May 2023), and backtest the results for two and four months into the future (June to September 2023).  The stocks were selected among the highest market capped companies per sector; see~\Cref{app:stocks} for specifics. We run $\starM$-optimization for $100$ iterations and use a regularization parameter of $\lambda = 10^{-2}$ to prevent overfitting.

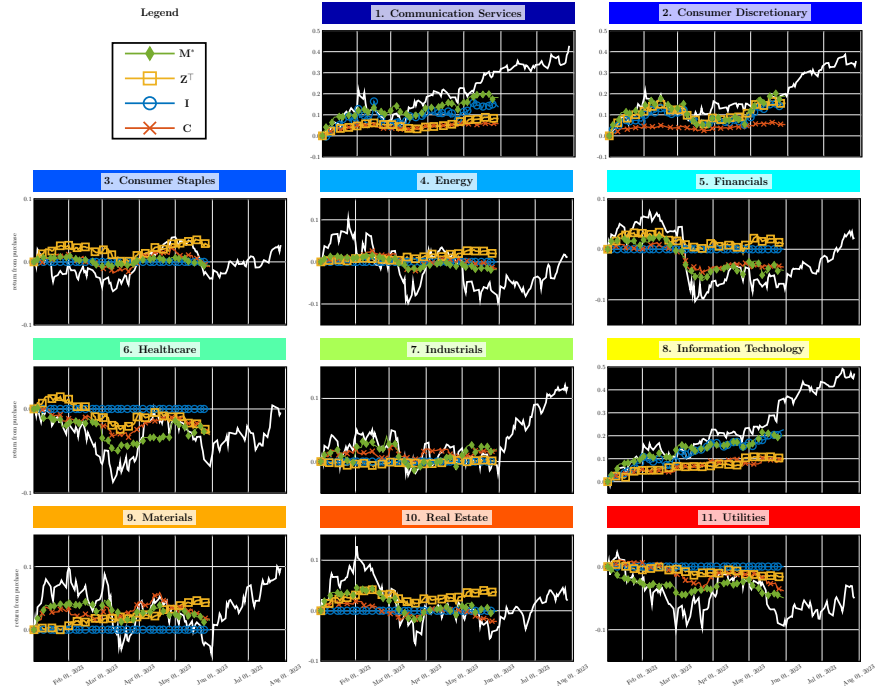
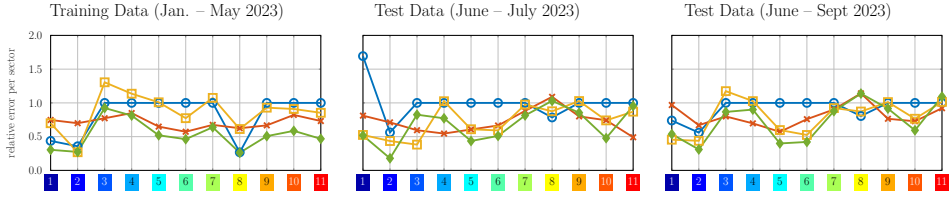
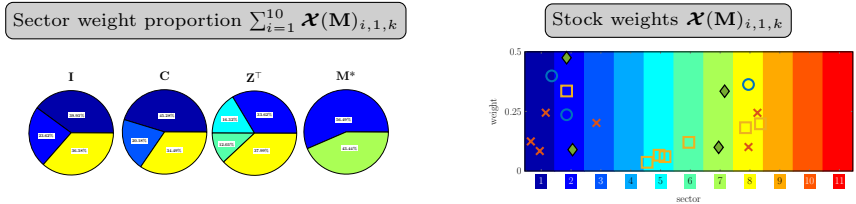
\begin{figure}
\centering

\subfloat[Tensor index tracking per sector, ordered alphabetically by sector.  The white line is the sector index per sector to track for the training time period of January 1, 2023 to May 31, 2023.  The index tracker from the learned $\bfM^*$ follows the market trends most faithfully during this time period. \label{fig:trainIndexTracking}]{
\begin{tikzpicture}
 \draw[white] (-0.485\linewidth,0) -- (0.485\linewidth,0);
\node at (0,0) {\includegraphics[width=0.925\linewidth]{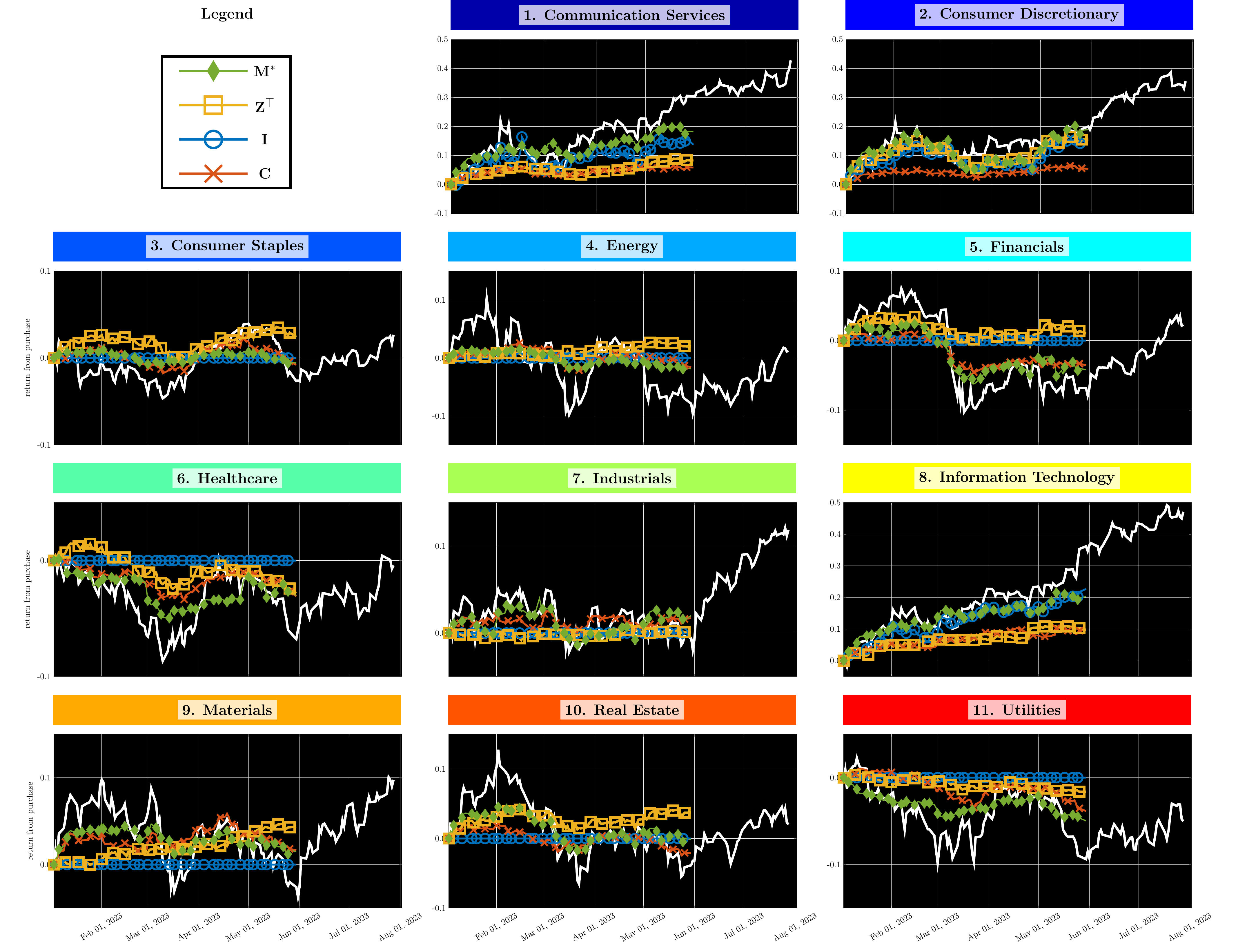}};
\end{tikzpicture}
}

\subfloat[Relative error per sector $\|(\TA \starM \TX(\bfM))_{:,:,i} - \TB_{:,:,i}\|_F / \|\TB_{:,:,i}\|_F$ for (left-to-right) historic training data and backtesting 2, 3, and 4 months into the future.  (Left): error on training period corresponding to plots in~\Cref{fig:trainIndexTracking}.  (Middle and Right): Backtesting two and four months into the future. 
\label{fig:relErrIndexTracking}]{
\includegraphics[width=0.98\linewidth]{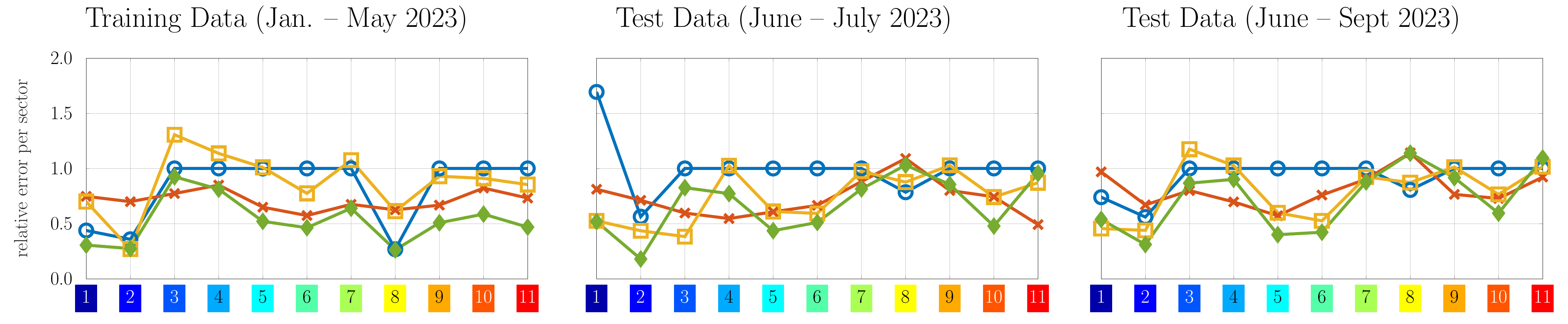}
}

\subfloat[Learned sector weight proportions $(\ge 5\%)$ and stock weights $(\ge 1\%)$. \label{fig:indexTrackingWeights}]
{
\begin{tikzpicture}
\scriptsize 
\draw[white] (-0.24\linewidth,0) -- (0.74\linewidth,0);
\node (n0) at (0,0) {\includegraphics[width=0.38\linewidth]{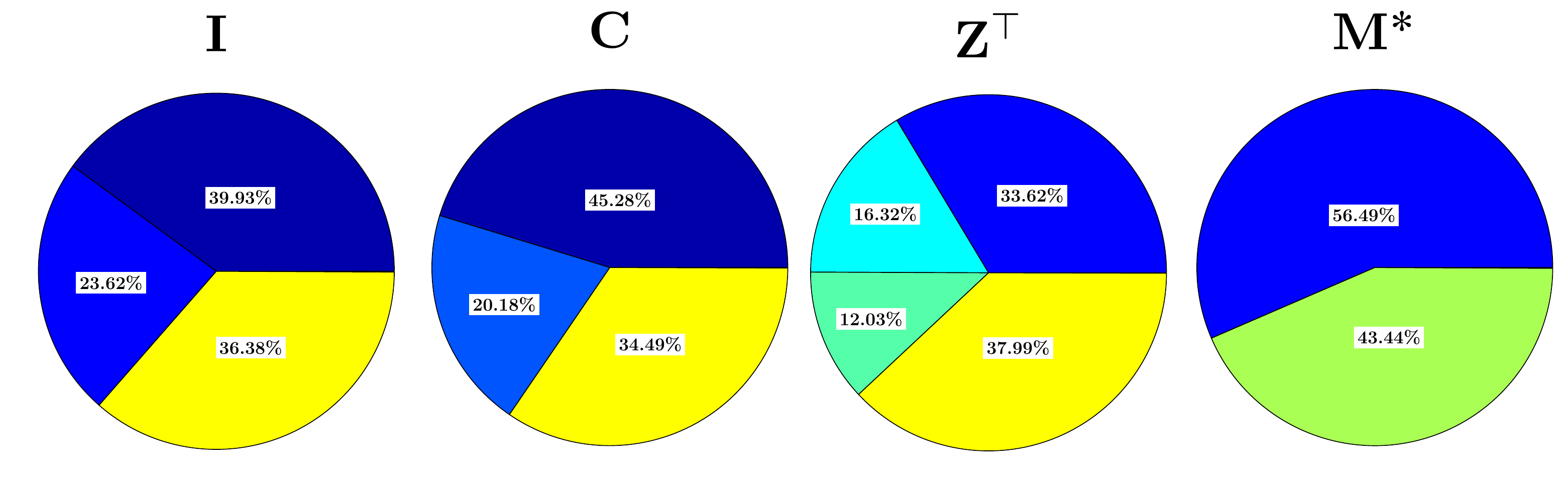}};
\node[right=1.0cm of n0.east, anchor=west] (n1) {\includegraphics[width=0.38\linewidth]{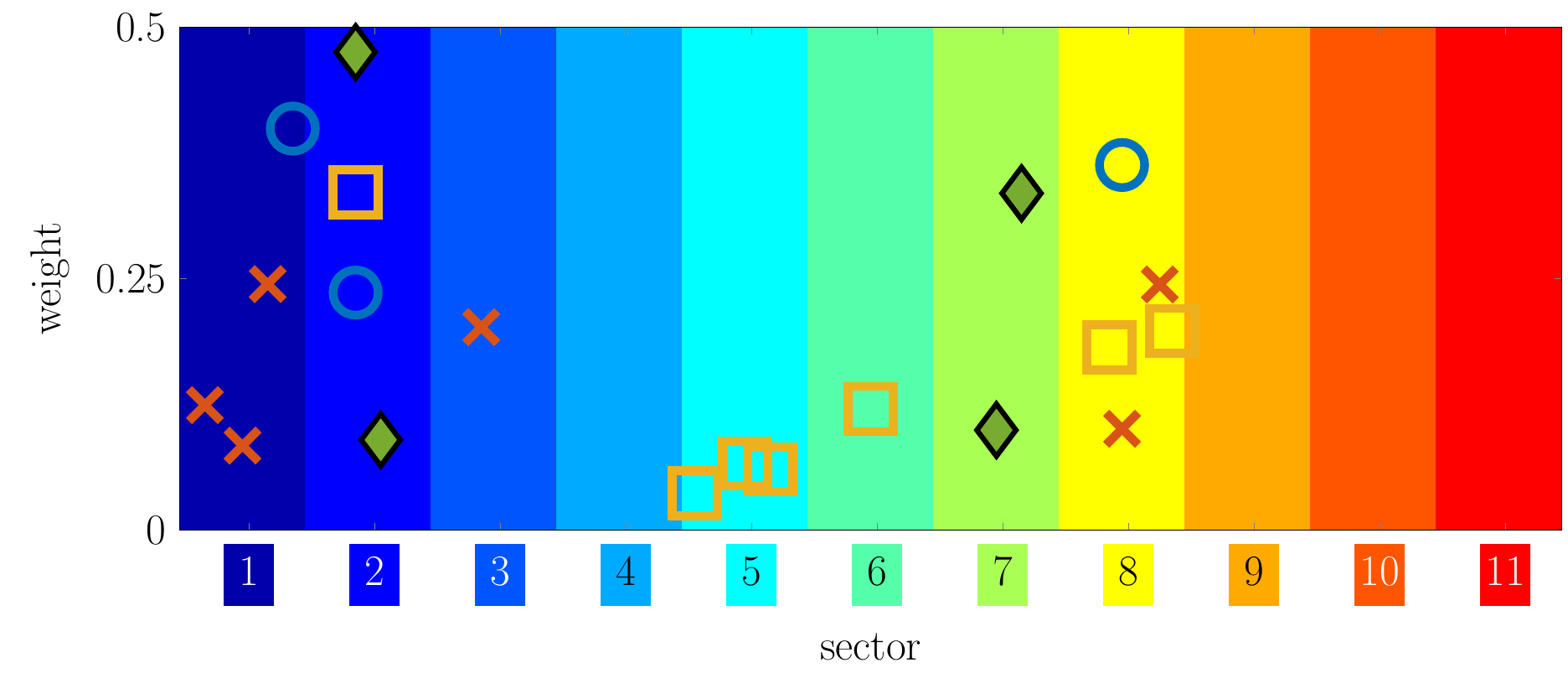}};

\node[above=0.0cm of n1.north, anchor=south, draw, rounded corners=0.1cm, fill=lightgray!75] (tt) {Stock weights $\TX(\bfM)_{i,1,k}$};
\node[draw, rounded corners=0.1cm, fill=lightgray!75] at (tt -| n0) {Sector weight proportion $\sum_{i=1}^{10}\TX(\bfM)_{i,1,k}$};

\end{tikzpicture}
}

\caption{Tensorized index tracking per sector for various choices of $\bfM$. The learned $\bfM^*$ tracks the sector indices best for the historic data and does comparatively well across sectors for future data. }
\label{fig:indexTracking}
\end{figure}

The results presented in~\Cref{fig:indexTracking} provide insight into the efficacy learning the $\starM$-product for tensor index tracking. Throughout our analysis, we use $\TX_{\rm ind}(\bfM)$ to denote learned portfolio weights for a given transformation $\bfM$ and $\widetilde{\TB}(\bfM) = \TA \starM \TX_{\rm ind}(\bfM)$ to denote the approximation of sector indices. 

\paragraph{Tracking Sector Indices}
In~\Cref{fig:trainIndexTracking}, $\widetilde{\TB}(\bfM^*)$ tracks trends in the training data more accurately than approximations associated with heuristic transformations. In comparison, $\widetilde{\TB}(\bfI)$ and  $\widetilde{\TB}(\bfZ^\top)$ are unable to capture trends per sector; e.g., for Healthcare, $\widetilde{\TB}(\bfI)$ remains nearly constant whereas $\widetilde{\TB}(\bfM^*)$ mirrors the downward trajectory. Similarly, $\widetilde{\TB}(\bfC)$ does not capture the magnitude of the trends as well as $\widetilde{\TB}(\bfM^*)$ (e.g., Consumer Discretionary). The left-most figure in~\Cref{fig:relErrIndexTracking} quantitatively confirms this result, showing the relative error per sector is lowest for $\widetilde{\TB}(\bfM^*)$.

\paragraph{Backtesting}
When forecasting two and four months ahead (middle and right figures, respectively, in~\Cref{fig:relErrIndexTracking}),
$\widetilde{\TB}(\bfM^*)$ achieves the lowest overall prediction error (Consumer Discretionary), performs the best or second best on at least eight out of $11$ sectors in both time periods, and produces the worst approximation only once (Utilities, four months ahead).  While the closest competitor, $\widetilde{\TB}(\bfC)$, has similar success, $\widetilde{\TB}(\bfM^*)$ achieves smaller overall relative errors, particularly for the longer prediction time. The data-dependent approximation, $\widetilde{\TB}(\bfZ^\top)$, achieved worse prediction results than $\widetilde{\TB}(\bfM^*)$ in all but one case for each of the backtesting time periods. We further note that $\widetilde{\TB}(\bfI)$ does not predict trends well, demonstrating the benefit of exploiting the correlation among sectors through tensor-tensor products for tensor index tracking.

\paragraph{Allocation of Portfolio Weights}
In~\Cref{fig:indexTrackingWeights} (left), we observe that $\TX_{\rm ind}(\bfM^*)$ allocates the majority of the weight in only two sectors, Healthcare and Consumer Discretionary. We conjecture that because Healthcare has a downward and flat trajectories and Consumer Discretionary has a mostly upward trajectory, tracking these sectors well enables subsequent combinations to capture behavior in other sectors.   In comparison, $\TX_{\rm ind}(\bfI)$ distributes weights to the sectors with the largest changes and market caps, resulting in suboptimal approximations of small market cap sectors. 

In~\Cref{fig:indexTrackingWeights} (right), we observe $\TX_{\rm ind}(\bfM^*)$ contains the largest weights overall. This indicates that the $\starM$-optimized weights are more amenable to sparsification than for the other cases. Comparatively, $\TX_{\rm ind}(\bfC)$ yields decent tracking performance and forecasting, but requires more smaller weights to achieve these results. Predicting market behavior with fewer stocks can be preferable, resulting in easier portfolio management. We mention that this is merely an observation; our setup does not explicitly encourage sparsity of the weights and differs from sparse index tracking~\cite{Benidis2018:indexTracking}. 

\paragraph{Variable Projection for Constrained Optimization} 
We learn the transformation $\bfM^*$ with a constrained inner optimization problem in~\eqref{eq:tensorIndexTracking},  demonstrating the generality of the $\starM$-optimization framework~\eqref{eq:bilevelOpt}. The necessary first-order optimality conditions are discussed in~\Cref{app:varproConstrained}. 

\paragraph{Comparison to Matrix Case}
We do not compare to the matrix index tracking problem~\eqref{eq:indexTrackingObjective}, which is a fundamentally different task. Specifically, the matrix case tracks a global index whereas the tensorized case tracks $11$ individual sector indices.

\subsection{Low-Rank Approximations of Digits}
\label{sec:digits}

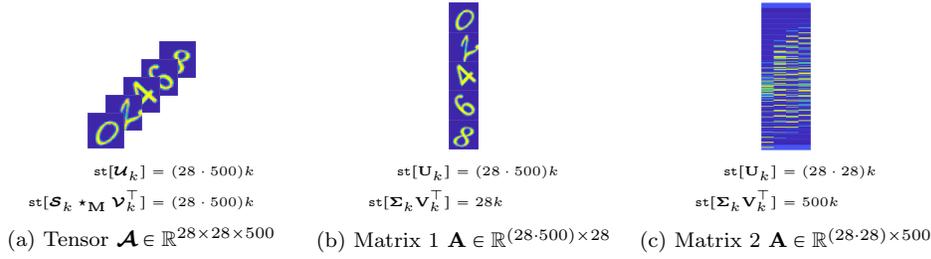
\begin{figure}
\centering
\subfloat[Tensor $\TA\in \Rbb^{28\times 28\times 500}$]{
\begin{tikzpicture} 
\tiny 
\draw[white] (-0.155\linewidth,0) -- (0.155\linewidth,0);
\node (n0) at (0,0) {\includegraphics[height=1.5cm]{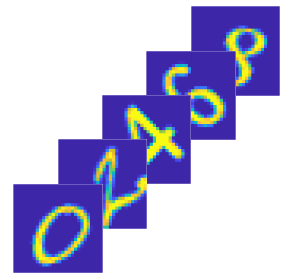}};
\node[below=0.0cm of n0.south, anchor=north] {$\begin{aligned} \texttt{st}[\TU_k] &= (28 \cdot 500)k \\ \texttt{st}[\TS_k \starM \TV_k^\top] &= (28 \cdot 500) k \end{aligned}$};
\end{tikzpicture}
}
\subfloat[Matrix 1 $\bfA\in \Rbb^{(28 \cdot 500)\times 28}$]{
\begin{tikzpicture} 
\tiny
\draw[white] (-0.155\linewidth,0) -- (0.155\linewidth,0);
\node (n0) at (0,0) {\includegraphics[height=2cm, trim={-3cm 0cm -3cm 0cm}, clip]{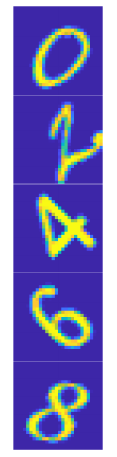}};
\node[below=0.0cm of n0.south, anchor=north] {$\begin{aligned} \texttt{st}[\bfU_k] &= (28 \cdot 500)k \\ \texttt{st}[\bfSigma_k \bfV_k^\top] &= 28 k \end{aligned}$};
\end{tikzpicture}
}
\subfloat[Matrix 2 $\bfA\in \Rbb^{(28 \cdot 28)\times 500}$]{
\begin{tikzpicture} 
\tiny 
\draw[white] (-0.155\linewidth,0) -- (0.155\linewidth,0);
\node (n0) at (0,0) {\includegraphics[height=2cm, trim={-3cm 0cm -3cm 0cm}, clip]{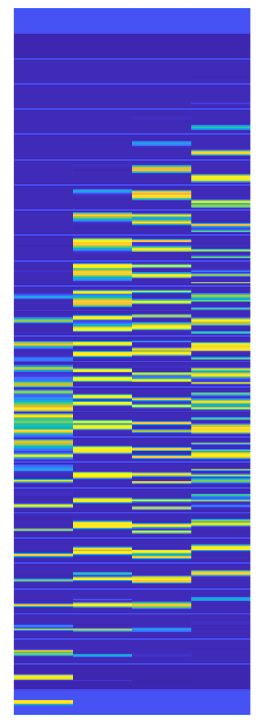}};
\node[below=0.0cm of n0.south, anchor=north] {$\begin{aligned} \texttt{st}[\bfU_k] &= (28 \cdot 28)k \\ \texttt{st}[\bfSigma_k \bfV_k^\top] &= 500 k \end{aligned}$};
\end{tikzpicture}
}
\caption{Illustration (not to scale) of Digits data for the three methods: tensor (images as frontal slices), matrix 1 (images concatenated vertically), and matrix 2 (vectorized images as columns).   The tensor and matrix 1 methods are directly comparable theoretically~\cite{Kilmer2021:pnas}. The matrix 2 method is typical for data compression. The function $\texttt{st}[\cdot]$ returns the number of floating point numbers needed to store the input.  The tensor method has the largest storage cost overall, but the cost to store the basis $\TU_k$ is the same as for  the matrix 1 method. }
\label{fig:digitsData}
\end{figure}

\begin{figure}
	\centering
	\subfloat[Visualizations of select images from approximation $\TA_k(\bfM)$ for $k=1$.  For each digit, we display the best and worst approximation based on the image relative error $\|\TA_{:,:,i} - \TX_{:,:,i}\|_F^2 / \|\TA_{:,:,i}\|_F$ for $i=1,\dots,500$ where $\TX = \TA_1(\bfM^*)$. The matrix approximations capture the ``average'' images and miss intricate details. 	Using the learned transformation qualitatively yields the closest approximations to the true digits, even in the worst-case scenarios.  \label{fig:digitsVisualization}]{
	\setlength\tabcolsep{0pt}
	\scriptsize
	\begin{tabular}{cc}
	best image for $\bfM^*$ & worst image for $\bfM^*$\\
	\includegraphics[width=0.49\linewidth]{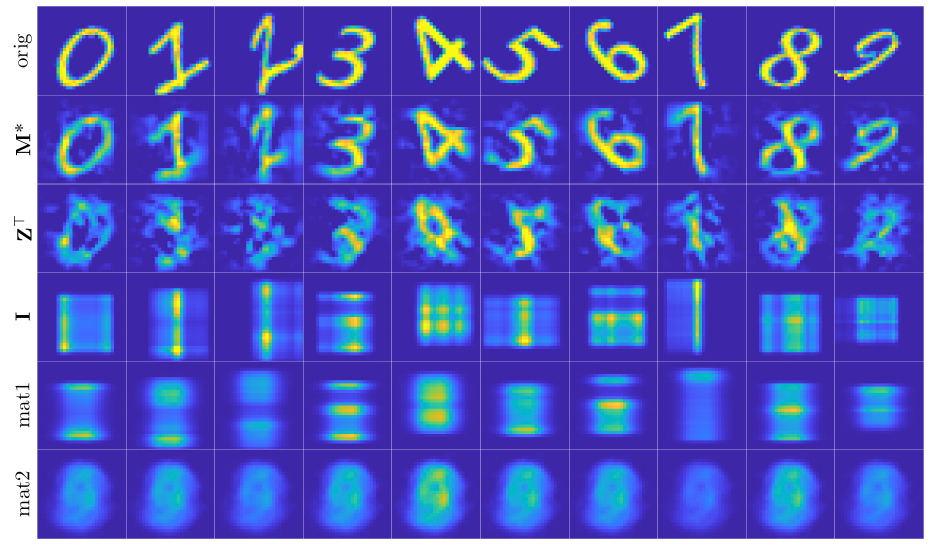}  &
	\includegraphics[width=0.49\linewidth]{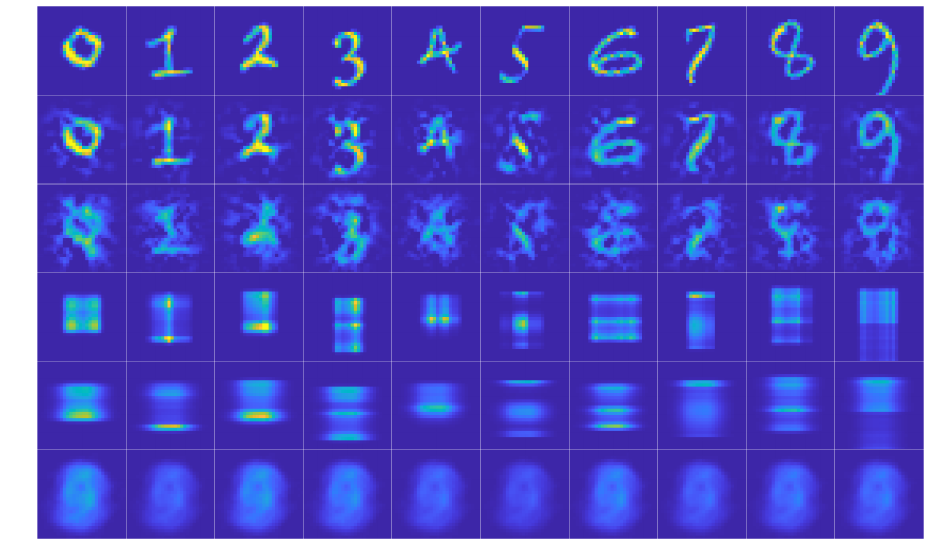}
	\end{tabular}}

	\subfloat[Structure of six frontal slices of $\widehat{\TA} = \TA \times_3 \bfM$ for various transformations.  To highlight the structure of the features in the transform domain, images are not shown on the same color scale.]{
	\begin{tikzpicture}
	
		\draw[white] (-0.18\linewidth,0) -- (0.78\linewidth,0);
		
		\node at (0,0) (n0) {\includegraphics[width=0.225\linewidth]{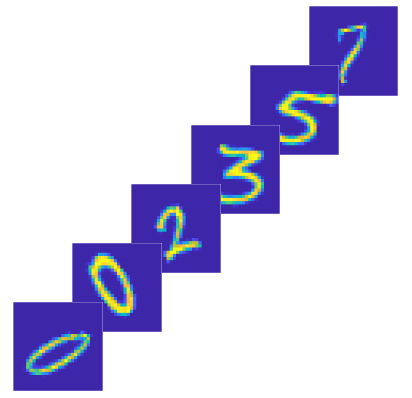}};
		\node[right=1cm of n0.east, anchor=west] (n1) {\includegraphics[width=0.225\linewidth]{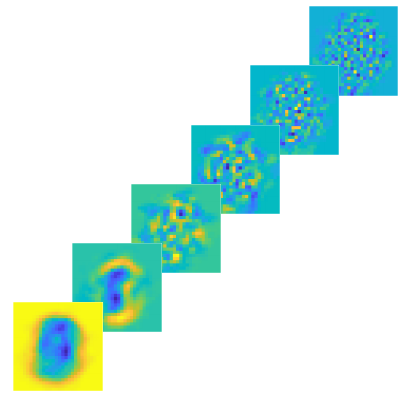}};
		\node[right=1cm of n1.east, anchor=west] (n2) {\includegraphics[width=0.225\linewidth]{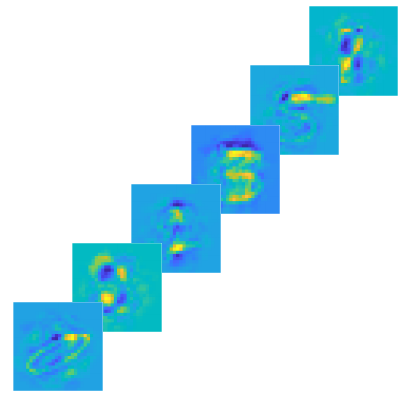}};
		
		\node[above left=-1.5cm of n0.north west, anchor=south east, rotate=45] {$\TA \times_3 \bfI$};
		\node[above left=-1.5cm of n1.north west, anchor=south east, rotate=45] {$\TA \times_3 \bfZ^\top$};
		\node[above left=-1.5cm of n2.north west, anchor=south east, rotate=45] {$\TA \times_3 \bfM^*$};
		
	\end{tikzpicture}
\label{fig:transformDomain}
	}

	 \caption{Approximation of Digits data and structure of learned features for various transformations $\bfM$ and truncation parameters $k$.}
	\label{fig:digitsApprox}
	\end{figure}
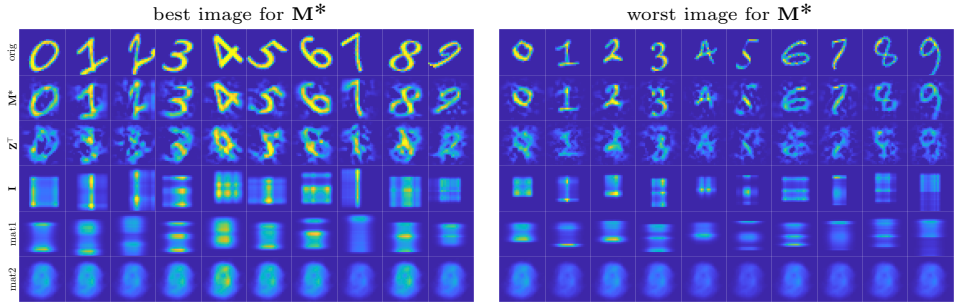
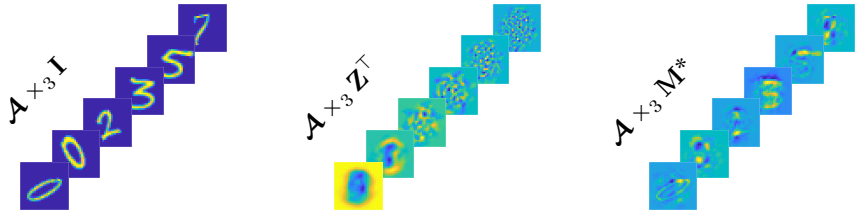

\begin{figure}
\centering

	\subfloat[Convergence and error of $\starM$-optimization. (Left): Convergence of the reduced objective function for various choices of truncation parameter $k$.  For every value of $k$, learning the transformation improves the approximation. (Right): global relative error for various transformations and truncations.   For all values of $k$, the learned transformation provides the best approximation. \label{fig:digitsConvergence}]{\begin{tikzpicture}
	\scriptsize
	\setlength\tabcolsep{3pt}
		\node[anchor=north] at (0,0) (n0) {\includegraphics[width=0.55\linewidth]{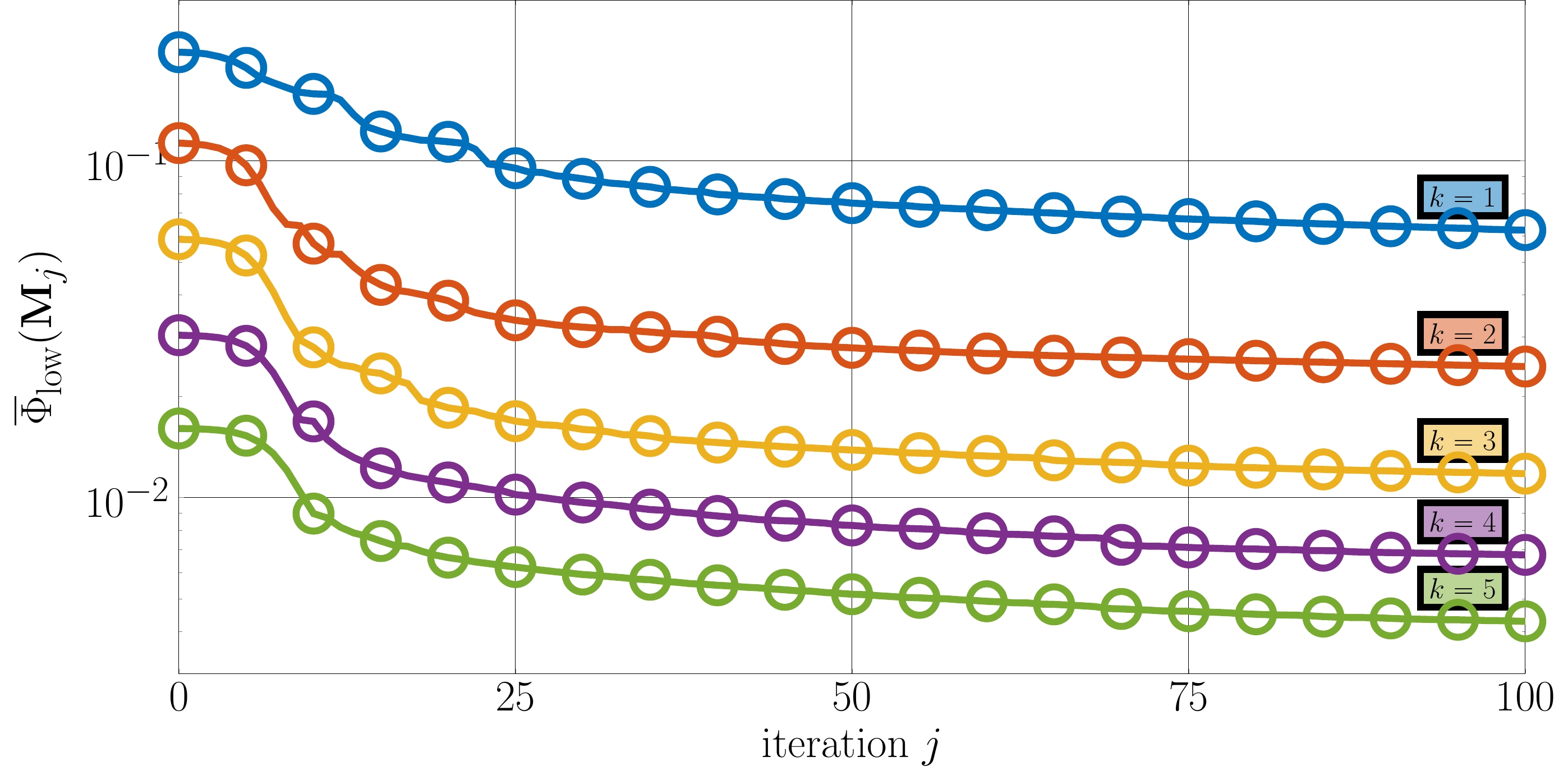}};
		\node[anchor=north, right=0.0cm of n0.east] (n1) { \scriptsize\begin{tabular}{|c||c|c|c|c|c|c|}
    \hline
    $k$ & mat1   &  mat2  &  $\bfI$ &   $\bfC$    &  $\bfZ^\top$ & $\bfM^*$\\
    \hline\hline
     \cellcolor{mycolor1!50}1	 & 0.73 & 0.80  & 0.65 & 0.63 & 0.61 & \cellcolor[gray]{0.7}  \bf 0.35\\
     \cellcolor{mycolor2!50}2 	& 0.63 & 0.77  & 0.48 & 0.48 & 0.46  &\cellcolor[gray]{0.7}  \bf 0.22\\
     \cellcolor{mycolor3!50}3 	&  0.54 & 0.74 & 0.34 & 0.38 & 0.35  &\cellcolor[gray]{0.7}  \bf 0.15\\
     \cellcolor{mycolor4!50}4 	& 0.45 & 0.72 & 0.23 & 0.30 & 0.28 &\cellcolor[gray]{0.7} \bf 0.12\\
     \cellcolor{mycolor5!50}5 	& 0.39 & 0.70 & 0.18 & 0.23 & 0.22 & \cellcolor[gray]{0.7} \bf 0.09\\
     \hline
    \end{tabular}};
    	\node[anchor=south, above=0.1cm of n1.north] {$\|\TA - \TA_k(\bfM)\|_F / \|\TA\|_F$};
	\end{tikzpicture}}

\subfloat[Heatmap of the transferability of $\bfM^*$ for various choices of $\bfM$ and truncation values $k$.   We create $50$ random batches of the Digits data and report the average relative error the approximation. Standard deviations were of the same order for all approximations.   For the retrained examples, we recompute the left-singular matrix or use an additional ten $\starM$-optimization iterations with two different initializations.  Darker hues indicate lower average relative error, and we bold the lowest relative error per $k$  with and without retraining in {bold}.\label{fig:transferability}]{\includegraphics[width=1\linewidth]{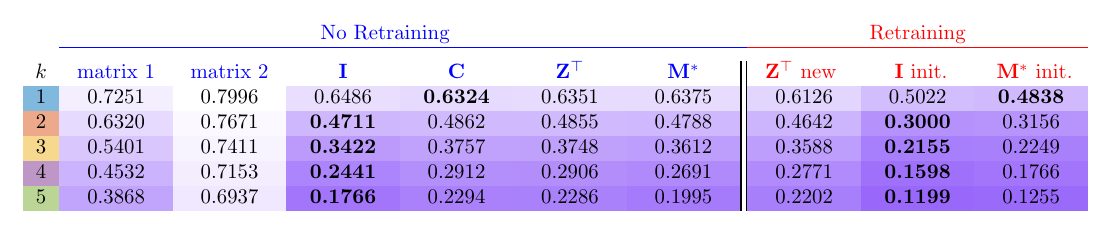}}

\caption{Convergence and transferability analysis for learning Digits data approximations for various truncation parameters $k$.}
\label{fig:transferability_convergence}
\end{figure}

Many data-driven algorithms that learn from large-scale, high-dimensional data hit the bottleneck of memory capacity. In this experiment, we consider compressing image data using the $t$-SVDM as a preprocessing step for machine learning.   We use the built-in {\sc Matlab} Digits dataset\footnote{Details on obtaining the Digits dataset can be found at \url{https://www.mathworks.com/help/deeplearning/ug/data-sets-for-deep-learning.html}.}, which consists of synthetically-generated $28\times 28$ grayscale images of rotated handwritten digits from $0$ to $9$.  We learn the optimal transformation for a batch of $50$ images per digit.  We store in both tensor and matrix formats; see~\Cref{fig:digitsData} for details. We run $\starM$-optimization with a line search for $100$ iterations, initialized using the identity matrix. We compare the approximation the truncated $t$-SVDM  and matrix SVD for various $\bfM$ and truncation parameters $k$ in~\Cref{fig:digitsApprox} and~\Cref{fig:transferability_convergence}. Throughout our analysis, we denote the truncated $t$-SVDM or SVD with $\TA_k(\bfM)$ or $\bfA_k$, respectively. 

\paragraph{Powerful Representations} 
In~\Cref{fig:digitsVisualization}, we observe qualitatively  that $\TA_1(\bfM^*)$ captures key digits features better than all other approximations. In contrast, $\TA_1(\bfI)$ computes rank-$1$ approximations of each image, resulting in rectangles that capture the locality and span of the digit, but none of the details. Conversely, $\TA_1(\bfZ^\top)$ captures the rotation and common curvature of the digits more than edge features. The $\starM$-optimized $\TA_1(\bfM^*)$, initialized using the identity matrix, captures both edges and curves, effectively combining the benefits of the two heuristic transformations.   We further demonstrate the learned features by examining the frontal slices in the transform domain in~\Cref{fig:transformDomain}. We see that $\TA\times_3 \bfZ^\top$ orders the frontal slices from low to high frequency (rotation) while $\TA \times_3 \bfM^*$ retains features of the original images.  Importantly, the frontal slices of  $\TA \times_3 \bfM^*$ appear closer to rank-$1$ than the original images, making the $t$-SVDM approximation more accurate.

\paragraph{Convergence and Sensitivity to Truncation Parameter}  
In~\Cref{fig:digitsConvergence} (left), we see that $\starM$-optimization convergences at roughly the same rate for all choices of $k$ and, as expected, reaches to lower values for less truncation (larger $k$).  
The table (right) shows that for every value of $k$, $\TA_1(\bfM^*)$ cuts the global approximation error roughly in half compared to the best-performing heuristic approach.  

\paragraph{Transferability}  In~\Cref{fig:transferability}, we show how well the learned transformation $\bfM^*$ generalizes to other similar batches of Digits data on average.  Without retraining, $\TA_k(\bfM^*)$ consistently achieves the second lowest errors. We note that $\TA_k(\bfI)$ yields the best approximations, but the other heuristic approaches are noticeably worse.  This shows that a good heuristic can outperform the learned transformation, but without a priori knowledge, the learned transformation can lead to reliable approximations.  Retraining with ten iterations $\starM$-optimization consistently produces the smallest overall relative error, including outperforming recomputing the left-singular matrix.

\subsection{Compressing Snapshot Data of Dynamical Systems}
\label{sec:rom}

Dynamical systems are central to simulating variety of real-world processes~\cite{HAY_BORGGAARD_PELLETIER_2009, ANTIL20121986, LieuFarhat2007}; however, solving parametric dynamical systems is computationally expensive. In cases where one must solve repeatedly, such as for parameter uncertainty quantification or in optimization constrained by partial differential equations (PDEs), the expense can be prohibitive. To decrease the cost of each PDE solve, projection-based reduced order modeling techniques have gained widespread use~\cite{BrennerGugercinWillcox2015:MOR}.  At a high level, the goal is to approximate high-dimensional state measurements $\bfx(t) \in \Rbb^{n_1}$ in a low-dimensional subspace $\bfx_k(t)\in \Rbb^k$ with $k \ll n_1$. Mathematically, we seek a ``good'' basis matrix $\bfU_k \in \Rbb^{n_1\times k}$ such that $\bfx(t) \approx \bfU_k \bfx_k(t)$ for all relevant times; that is, the solution in the reduced space well-approximates the dynamics of the original model.

The proper orthogonal decomposition (POD) is an SVD-based approach to form a global basis across the parameters~\cite{Sirovich1987}.  The heart of POD is the \emph{method of snapshots}, which builds a basis that captures dominant trends from simulations. Procedurally, we construct a snapshot matrix for a fixed parameter at discrete time points
	\begin{align}\label{eq:snapshotMatrix}
	\bfX(c_i) = \begin{bmatrix}
		\vert & \vert &  &\vert \\
		\bfx(t_0; c_i) & \bfx(t_1; c_i) & \cdots & \bfx(t_{n_2-1}; c_i) \\
		\vert & \vert &  &\vert 
		\end{bmatrix}
	\end{align}
where $\bfx(t_j; c_i)$ is the discretized high-dimensional state for parameter $c_i$ at time $t_j$.   To capture dynamical behavior across the parameters, traditional matricized POD concatenates local snapshot matrices in a global snapshot matrix
	\begin{align}
	\bfX = \begin{bmatrix}
		\bfX(c_1) & \bfX(c_2) & \cdots & \bfX(c_{n_3})
		\end{bmatrix},
	\end{align}
then forms a reduced basis via truncated matrix SVD $\bfX \approx \bfU_k \bfSigma_k \bfV_k^\top.$

In practice, we may expect some relationship across the parameter dimension, and may lose these correlations when we construct the global snapshot matrix. In this experiment, we consider storing each $\bfX(c_\ell)$ as a frontal slice; that is, 
	\begin{align}
	\TX_{:,:,i} = \bfX(c_i) \qquad \text{for $i=1,\dots, n_3$}. 
	\end{align}
The resulting snapshot tensor $\TX$ has dimensions corresponding to $\text{space} \times \text{time} \times \text{parameters}$.  Our goal is to generate a tensor POD basis $\TU_k$ by computing the $t$-SVDM of the snapshot tensor.   The quality of our basis using the relative error of the projection to the true states  
	\begin{align}
	\frac{\|\TX - \TU_k \starM \TU_k^\top \starM \TX\|_F }{\|\TX\|_F}.
	\end{align}

We setup a homogeneous two-dimensional wave equation from the {\sc Matlab} PDE Toolbox\footnote{\url{https://www.mathworks.com/help/pde/ug/wave-equation.html}}, parameterized by wave speed $c$. The specifics can be found in~\Cref{sec:wave_equation}. 

Snapshots of the solution to this hyperbolic PDE will be highly compressible because of the periodic behavior.  Different wave speeds will yield the same wave amplitudes and patterns at different time scales.  Thus, this is an ideal test problem to explore $\starM$-based POD because of the commonalities across parameters. 

\begin{figure}
\centering
\scriptsize 
        \begin{tikzpicture}
            \node (n0) {\includegraphics[width=0.35\linewidth]{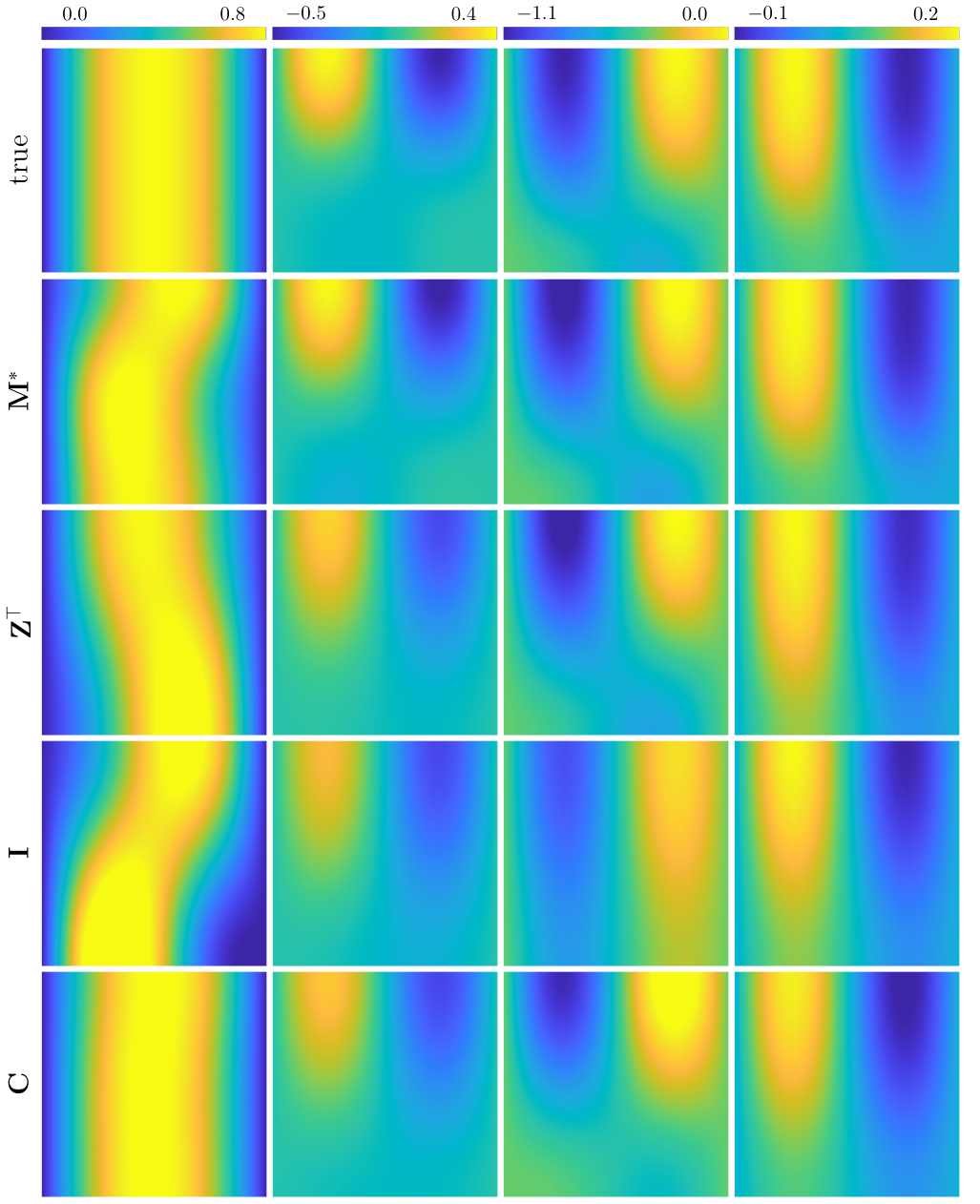}};
            \node[right=0.0cm of n0.south east, anchor=south west] (n1) {\includegraphics[width=0.35\linewidth]{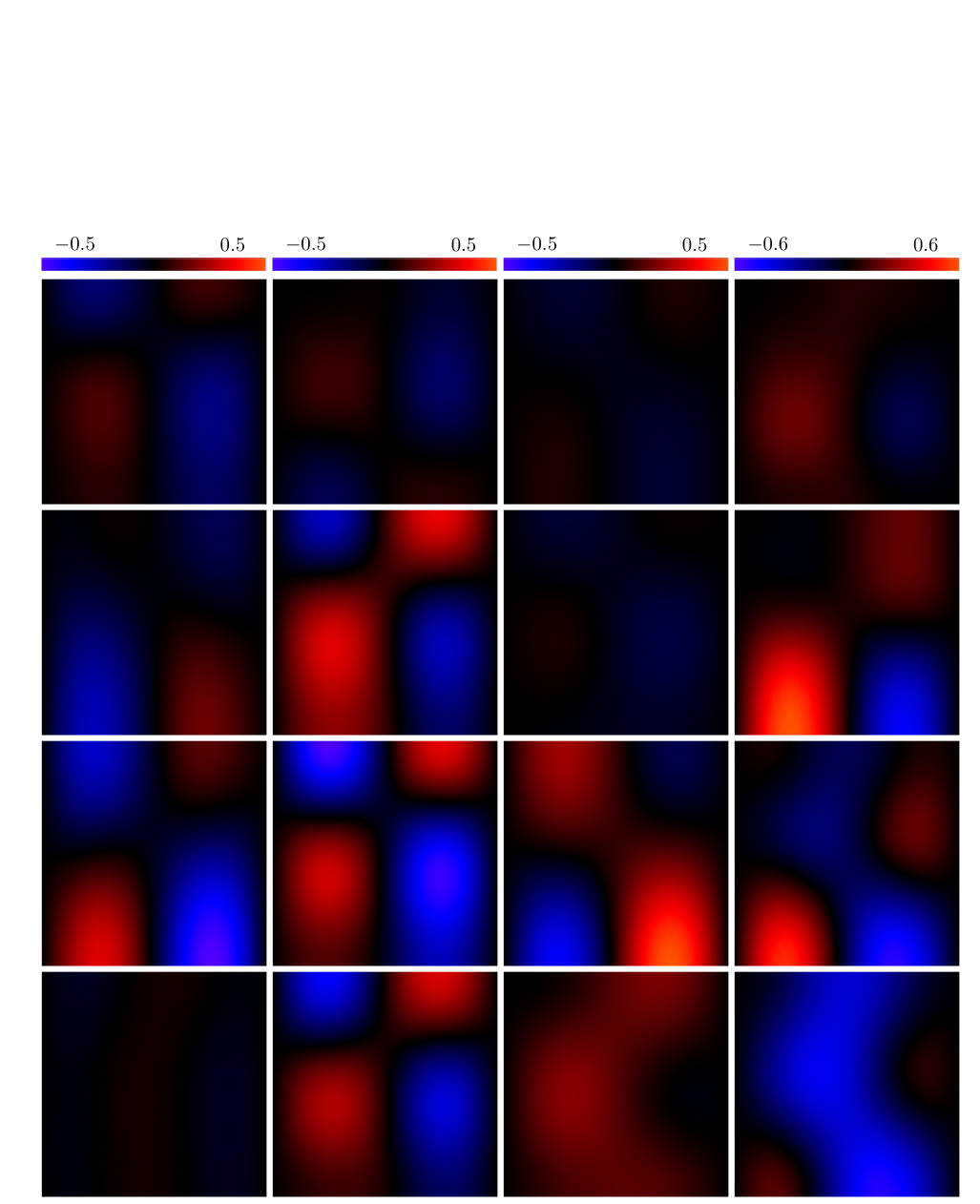}};
            \draw[line width=2pt, ->, shorten <=1cm, shorten >=0.75cm]  (n0.south west) --node[below, midway] {time}  (n0.south east);
            \node[above=0.0cm of n0.north, anchor=south, draw, rounded corners=0.1cm, fill=lightgray!75] {Approximations $\TX_k(\bfM)_{:,:,i}$};
            \node[above=0.0cm of n1.north, anchor=south, draw, rounded corners=0.1cm, fill=lightgray!75] {Error $\TX_k(\bfM)_{:,:,i} - \TX_{:,:,i}$};

            \draw[line width=2pt, ->, shorten <=1cm, shorten >=0.75cm]  (n1.south west) --node[below, midway] {time}  (n1.south east);

        \end{tikzpicture}

\medskip

\begin{tikzpicture}
\scriptsize
\node at (0,0) (n1) {\includegraphics[width=0.725\linewidth]{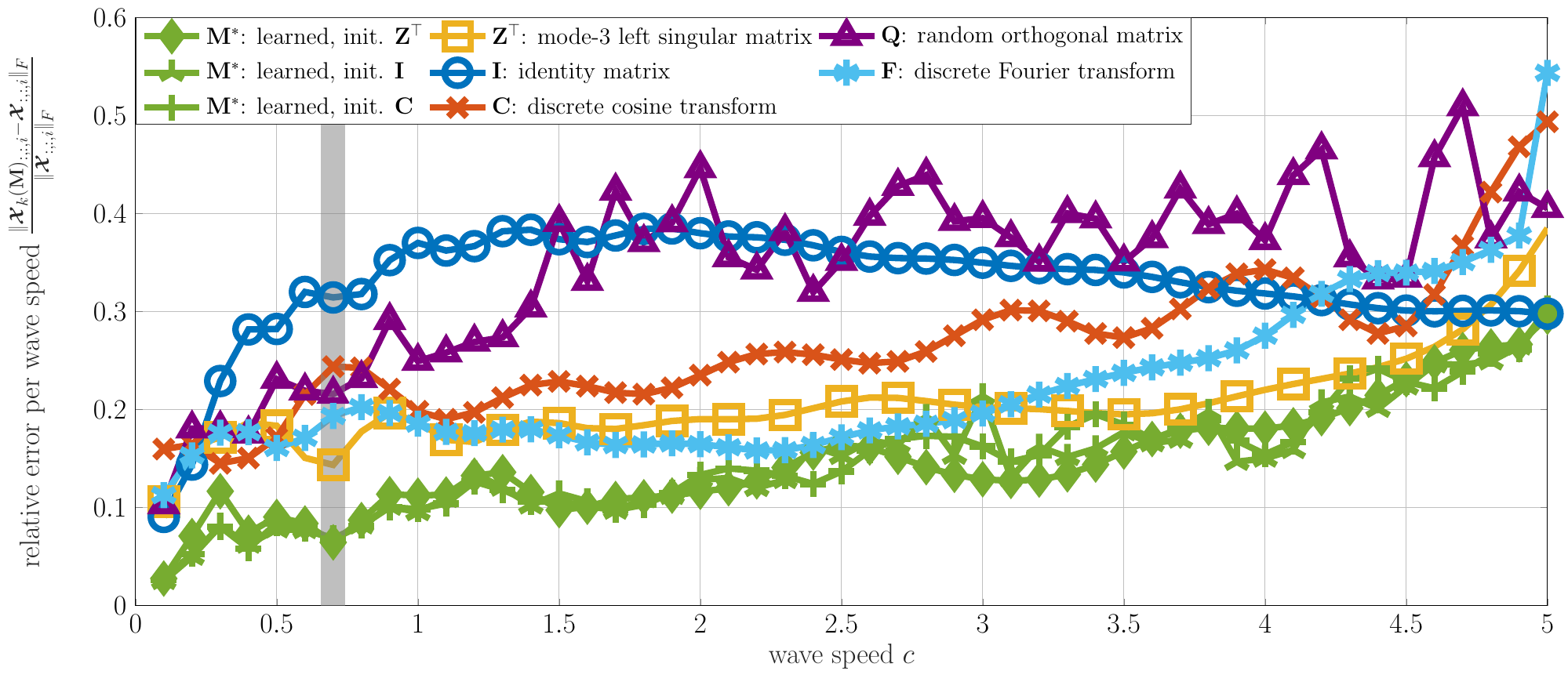}};
\node[above=0.0cm of n1.north, anchor=south, draw, rounded corners=0.1cm, fill=lightgray!75] (tt) {Relative Error per Wave Speed};
          \node[right=0.0cm of n1.east, anchor=west] (t) {
          \scriptsize
          \renewcommand{\arraystretch}{1.05}
          \setlength\tabcolsep{4pt}
          \begin{tabular}{|c||r|r|}
          \hline
	  & \multicolumn{1}{c}{$\bfM_0$} &\multicolumn{1}{|c|}{ \cellcolor{mycolor5!75}{$\bfM^*$}} \\
	 \hline\hline
	 \cellcolor{mycolor3!75}{ $\bfZ^\top$ }& 0.184 &  0.124\\
	\cellcolor{mycolor1!50}{$\bfI$} & 0.293 & \bf 0.122\\
	\cellcolor{mycolor2!75}{$\bfC$} & 0.232 &\bf 0.122\\
	\cellcolor{mycolor6!50}{$\bfF$} &  0.199 & \----\\
	\cellcolor{mycolor4!50}{$\bfQ$} &  0.291 & \----\\
	\hline
\end{tabular}};

	\node[above=0.0cm of t.north, anchor=south] {$\|\TX - \TX_k(\bfM)\|_F/\|\TX\|_F$};
	
	\node[draw, rounded corners=0.1cm, fill=lightgray!75] at (tt -| t) {Global Relative Error};

\end{tikzpicture}

\caption{Comparisons of tensor POD bases of a homogeneous 2D wave equation parameterized by wave speed $c$.  (Top): for one wave speed  ($c=0.7$) and truncation $k=2$, we show the approximation of PDE solution (left) and the error (right) where \red{\bf red} indicates overestimation and \blue{\bf blue} indicates underestimation. The {\color{gray} \bf gray} vertical line in the plot below indicates the wave speed used for the top visualization.  We display the case with learn $\bfM^*$ initialized with $\bfM_0 = \bfI$.  (Bottom): Relative error per wave speed (left) and the global relative error (right).   The learned $\bfM^{*}$ consistently produces the lowest error for multiple initial guesses.}
\label{fig:pdeApprox}

\end{figure}
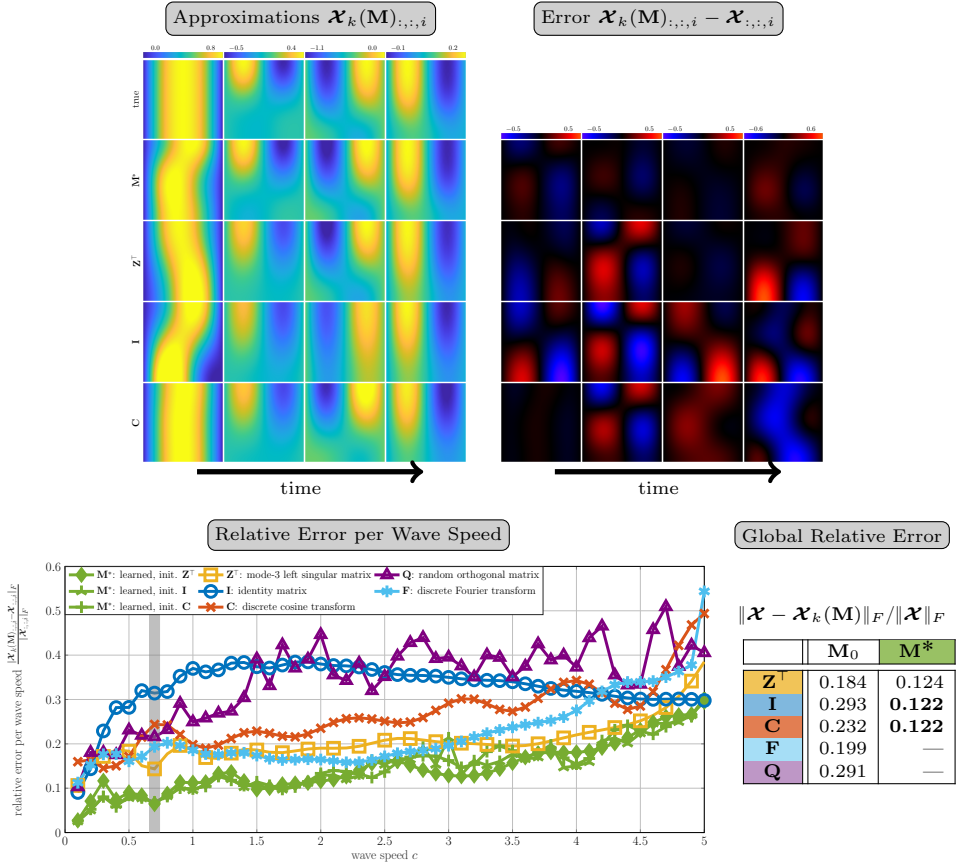

\begin{figure}
\centering

	\begin{tikzpicture}

\scriptsize
\node (n0) at (0,0) {};

\foreach[count=\a] \M in {Z, I, C}{
	\pgfmathsetmacro{\b}{int(\a-1)}
	\node[below=0.0cm of n\b.south, anchor=north] (n\a) {\includegraphics[width=0.18\linewidth, trim={1cm 1cm 0cm 0cm}, clip]{figures/rom2/romSpectrum_\M_init_\M.pdf}};
	\node[right=0.0cm of n\a.east, anchor=west] (m\a) {\includegraphics[width=0.18\linewidth, trim={1cm 1cm 0cm 0cm}, clip]{figures/rom2/romSpectrum_MOpt_init_\M.pdf}};
	
}

\node[left=-0.15cm of n2.west, anchor=south, rotate=90] {percent of energy};

\node (mm) at ($(n3.south)!0.5!(m3.south)$) {};
\node[below=-0.25cm of mm.south, anchor=north] (xx) {singular value index};

\node[above=0.0cm of n1.north, anchor=south] (x0) {initial $\widehat{\TX}_k(\bfM_0)$};
\node[above=0.0cm of m1.north, anchor=south] (xopt) {learned $\widehat{\TX}_k(\bfM^*)$};

\node[draw, rounded corners=0.1cm, fill=lightgray!75] at ($($(x0)!0.5!(xopt)$) + (0,0.5)$) {Percent of Energy Retained};

\node[left=0.5cm of n1.west, anchor=south, rotate=90] {$\bfM_0=\bfZ^\top$};
\node[left=0.5cm of n2.west, anchor=south, rotate=90] {$\bfM_0=\bfI$};
\node[left=0.5cm of n3.west, anchor=south, rotate=90] {$\bfM_0=\bfC$};

\node[right=0.25cm of m2.east, anchor=north, rotate=90] (cb) {\includegraphics[width=0.4\linewidth]{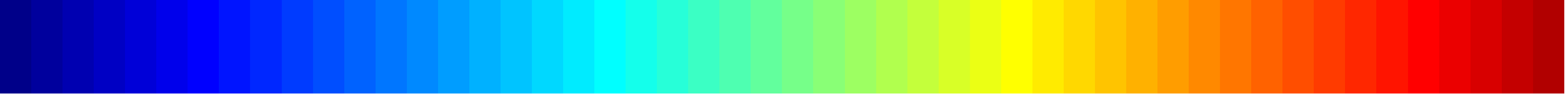}};
\node[above=0.0cm of cb.west, anchor=north] {$i=1$};
\node[above=0.0cm of cb.east, anchor=south] {$i=50$};
\node[left=0.0cm of cb.south, anchor=south, rotate=270] {frontal slice index $i$};

\foreach[count=\a] \M/\errA/\errB in {
	Z/\colorbox{white}{$0.1842$}/\colorbox{white}{$0.1239$},
	I/\colorbox{white}{$0.2927$}/\colorbox{yellow!50}{$\bf 0.1218$}, 
	C/\colorbox{white}{$0.2315$}/\colorbox{white}{$0.1224$}
	}{
	\pgfmathsetmacro{\b}{int(\a-1)}
	\node[right=1.5cm of m\a.north east, anchor=north west] (nn\a) {\includegraphics[width=0.15\linewidth]{figures/rom2/matrix_\M.png}};

	\node[right=0.0cm of nn\a.east, anchor=west] (mm\a) {\includegraphics[width=0.15\linewidth]{figures/rom2/MOpt--init-\M.png}};
	
	\ifthenelse{\a = 3}{
		\draw[<->, thick, shorten <= 0.15cm, shorten >= 0.15cm] ($(nn\a.south west)+(0,-0.16)$) -- node[midway, fill=white] {$n_3$} ($(nn\a.south east)+(0,-0.16)$);
		
		\draw[<->, thick, shorten <= 0.15cm, shorten >= 0.15cm] ($(nn\a.south west)+(-0.16,0)$) -- node[midway, fill=white] {$n_3$} ($(nn\a.north west)+(-0.16,0)$);
	}{}

}

 \node[right=0.0cm of mm2.east, anchor=west] (c) {\includegraphics[height=0.4\linewidth, trim={2cm 0cm 0cm 0cm}, clip]{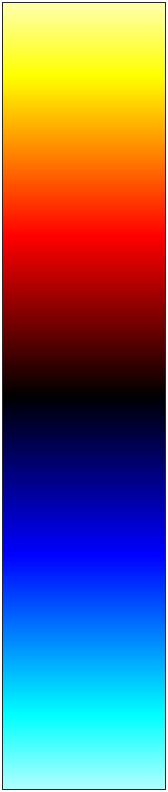}};
 \node[right=-0.15cm of c.south east, anchor=west] {$-1$};
  \node[right=-0.15cm of c.north east, anchor=west] {$+1$};
  \node[right=-0.15cm of c.east, anchor=west] {$0$};

\node (t1) at (xopt -| nn1) {initial $\bfM_0$};
\node (t2) at (xopt -| mm1) {learned $\bfM^*$};

\node[draw, rounded corners=0.1cm, fill=lightgray!75] at ($($(t1)!0.5!(t2)$) + (0,0.5)$) {Structure of Transformations};

\end{tikzpicture}

\caption{Comparison of spectra and transformation matrices for various algorithm initializations. 
(Left): Percentage of energy captured by $k$ singular values per frontal slice in transform domain.  The cumulative energy of the $i$-th frontal slice is given by $\|\widehat{\TX}_k(\bfM)_{:,:,i}\|_F^2. $ Approximations under $\bfM^*$ (right column) generally capture more energy per frontal slice at $k=2$ than the heuristic cases. For $\bfZ^\top$ (top row), the energy for the back frontal slices (redder colors) is negligible ($\approx \Ocal(10^{9})$ less energy than the first frontal slices), hence we see little change in those spectra. (Right): Visualizations of the learned $n_3\times n_3$ transformations $\bfM$ where $n_3 = 50$ for various initial guesses.  
}

\label{fig:romSpectra}
\end{figure}

Using {\sc Matlab}'s PDE toolbox, we solve~\eqref{eq:wave_pde} using finite elements for $50$ different wave speeds, equispaced between $c=0.1$ and $c=5$.  We construct snapshots from $31$ equispaced timepoints from $t=0$ to $t=5$.  
The resulting snapshot tensor $\TX$ is of size $493\times 31\times 50$, where $493$ is the number of finite element nodes.  We run $\starM$-optimization for $1000$ iterations using truncation parameter $k=2$ and present the results in~\Cref{fig:pdeApprox} and~\Cref{fig:romSpectra}. We will use $\TX_k(\bfM)$ to denote the low-$t$-rank approximation to the snapshot tensor for a given $\bfM$. 

\paragraph{Approximations Per Parameter} In~\Cref{fig:pdeApprox}, we observe that $\TX_2(\bfM^*)$ achieves a lower approximation error per parameter compared to the competing representations \emph{and} harnesses multilinear relationships to obtain local consistency.  Similar performance is achieved for three different initializations, demonstrating robustness of $\starM$-optimization to the initial guess. The learned approximation reduces the global relative error by more than $30\%$ compared to the strongest competitors, $\TX_2(\bfZ^\top)$ and $\TX_2(\bfF)$.  The least accurate approximations, $\TX_2(\bfI)$ and $\TX_2(\bfQ)$, do not exploit commonalities among wave solutions across parameters, further demonstrating the need to leverage multilinear correlations for quality low-$t$-rank representations.

\paragraph{Analysis of Spectra} In~\Cref{fig:romSpectra} (left), we observe that the learned transformation adjusts the spectra such that more energy is given to the first $k=2$ singular values. For various initial guesses, approximation $\widehat{\TX}_k(\bfM^*)$ can capture about $95\%$ of the energy content per frontal slice using the two largest singular values.  In comparison, the other methods consist of more frontal slices that capture only about $80\%$ of the energy using $t$-rank-$2$ approximations.

\paragraph{Structure of Learned $\bfM^*$} 
In~\Cref{fig:romSpectra} (right),  we observe that $\bfM^*$ retains properties of the initial guess and modifies the necessary rows to improve the global accuracy of the approximation. Each row of the transformations maps to a frontal slice in the transform domain; that is, $\widehat{\TY}_{:,:,i} \equiv \TY\times_3 \bfM_{i,:}$ for $i=1,\dots,n_3$.  The modified rows reflect changes in the spectra in~\Cref{fig:romSpectra} (left). For example, for $\bfM_0 = \bfC$, the first $20$ to $30$ darker blue energy curves in $\widehat{\TX}(\bfM^*)$ increase the most, which corresponds to the significant changes of the first $20$ to $30$ rows of $\bfM^*$.

\paragraph{Comparison to Matrix POD}   In order to compare fairly to the matrix POD case, we need to develop methods to remove the parameter dimension (mode-$3$) of the tensor basis, $\TU_k$. However, this would require new theoretical developments beyond the scope of the presented work. We leave this to future work.

\section{Conclusions}
\label{sec:conclusions}

We presented a new $\starM$-optimization paradigm that leverages the matrix mimeticity of the $\starM$-framework to simultaneously learn the underlying tensor algebra and obtain optimal representations of multilinear data.  In~\Cref{sec:invariantProduct}, we proved invariant properties of the $\starM$-product under row permutations of the transformation.In~\Cref{sec:varpro}, we introduced our variable-projection framework and highlighted the generality of our formulation through two prototype problems: $t$-linear regression and low-$t$-rank approximations. We derived the gradients needed to learn orthogonal transformations through Riemannian optimization, including differentiating through the $t$-SVDM (\Cref{sec:tsvdmDerivatives}) and proved that $\starM$-optimization will converge for $t$-linear regression (\Cref{cor:convergence}). In~\Cref{sec:numerical}, we applied $\starM$-optimization to a wide range of applications, including index tracking, image compression, and reduced order modeling.

This paper serves as the seminal work of the $\starM$-optimization framework and opens the door to many new research directions. Currently, we are restricted to inner optimization problems that have analytic solutions. Following approaches in~\cite{Newman2021:GNvpro}, we will generalize to convex inner problems and apply $\starM$-optimization to, e.g., tensor completion. In future work, we will extend $\starM$-optimization to complex-valued and higher-order tensors, following similar extensions seen in~\cite{Kilmer2021:pnas, Keegan2022}. To scale up to larger problems, we will reduce the computational and storage costs by exploring more compressible variants of the $t$-SVDM~\cite{Kilmer2021:pnas}, randomized versions of least squares and the $t$-SVDM~\cite{murray2023randomized, zhang2016randomized}, and matrix-free implementations of the transformation. We will further explore second-order variants of Riemannian optimization to accelerate convergence and reduce training time~\cite{EdelmanAriasSmith1998}.  Finally, we will apply $\starM$-optimization to new applications, such as X-ray spectroscopy~\cite{Townsend:22} and hyperspectral image analysis~\cite{Cui:2017indianPines}.

\bibliographystyle{siamplain}
\bibliography{references}


\appendix

\input{sections/S_00_introduction}

 \input{sections/S_01_underlying_algebra}

\input{sections/S_09_riemannian_optimization}

\input{sections/S_08_fundamental_derivatives}

\input{sections/S_02_tsvdm_derivatives}

\input{sections/S_02_02_tsvdm_derivatives}

\input{sections/S_03_prototype_invariance}

\input{sections/S_04_bounded_operator_norm}

\input{sections/S_07_bounded_hessian}

	\input{sections/S_07_01_outline}

\input{sections/S_07_02_bounded_gradient}

\input{sections/S_07_03_bounded_hessian_implication}
	\input{sections/S_07_04_technical_details}

\input{sections/04_01_geometric_intuition}

\input{sections/S_05_varpro_unit_simplex}

\input{sections/S_06_stocks}

\input{sections/S_10_wave_equation}

\end{document}

%% file: ex_shared.tex
\ifpdf
  \DeclareGraphicsExtensions{.eps,.pdf,.png,.jpg}
\else
  \DeclareGraphicsExtensions{.eps}
\fi

\headers{Optimal Tensor Algebras}{E.Newman and K. Keegan}

\title{Optimal Matrix-Mimetic Tensor Algebras via Variable Projection\thanks{Submitted to the editors~\today.
\funding{The work by E. Newman was partially supported by the National Science Foundation (NSF) under
grants [DMS-2309751] and the work by K. Keegan was supported by the Department of Energy Computational Science Graduate Fellowship [DE-SC0023112]. }}}

\author{Elizabeth Newman\thanks{Emory University, Atlanta, GA
  (\email{elizabeth.newman@emory.edu}, \href{https://math.emory.edu/~enewma5/}{https://math.emory.edu/\textasciitilde{enewma5}/}).}
\and Katherine Keegan\thanks{Emory University, Atlanta, GA
  (\email{katherine.emiri.keegan@emory.edu}, \url{https://katiekeegan.org/}).}}

\input{macros}


%% file: macros.tex
\usepackage{amsmath, amsfonts, amssymb} 

\usepackage{graphicx}
\usepackage{tikz, pgfplots, pgflibraryplotmarks, pgfkeys, ifthen}
\usepackage{pgf-pie} 
\usetikzlibrary{calc}
\usepackage{amssymb, amsfonts, amsmath, latexsym, dsfont}
\usepackage{comment}
\usepackage{array}
\pgfplotsset{compat = 1.3}
\usepgfplotslibrary{polar}
\usepackage{xfrac}
\usetikzlibrary{decorations.text}

\usetikzlibrary{shapes.arrows, arrows, decorations.markings, positioning}
\usetikzlibrary{calc}
\usetikzlibrary{3d}
\usepackage{graphicx} 
\usepackage{xcolor} 
\usepackage{colortbl}
\usepackage{multicol} 
\usepackage{hyperref} 
\usepackage{mathabx}
\usepackage{cleveref}
\usepackage{subfig}
\usepackage{array}
\usepackage{wrapfig}
\usepackage{bbm}
\usepackage{lipsum}
\usepackage{arydshln}
\usepackage{enumitem}
\usepackage{algorithm, algorithmicx, algpseudocode}

\numberwithin{equation}{section}

\newtheorem{example}[theorem]{Example}

\theoremstyle{remark}
\newtheorem*{remark}{Remark}

\definecolor{EmoryBlue}{RGB}{1, 33, 105} 
\definecolor{EmoryDarkBlue}{RGB}{12, 35, 64} 
\definecolor{EmoryMediumBlue}{RGB}{0, 51, 160} 
\definecolor{EmoryLightBlue}{RGB}{0, 125, 186} 
\definecolor{EmoryYellow}{RGB}{242, 169, 0} 
\definecolor{EmoryGold}{RGB}{181, 133, 0} 
\definecolor{EmoryMetallicGold}{RGB}{132, 117, 78}

\definecolor{mycolor0}{rgb}{1, 1, 1}%
\definecolor{mycolor1}{rgb}{0.00000,0.44700,0.74100}%
\definecolor{mycolor2}{rgb}{0.85000,0.32500,0.09800}%
\definecolor{mycolor3}{rgb}{0.92900,0.69400,0.12500}%
\definecolor{mycolor4}{rgb}{0.49400,0.18400,0.55600}%
\definecolor{mycolor5}{rgb}{0.46600,0.67400,0.18800}%
\definecolor{mycolor6}{rgb}{0.30100,0.74500,0.93300}%

\definecolor{jet1}{rgb}{0.0000,0.0000,0.6667}
\definecolor{jet2}{rgb}{0.0000,0.0000,1.0000}
\definecolor{jet3}{rgb}{0.0000,0.3333,1.0000}
\definecolor{jet4}{rgb}{0.0000,0.6667,1.0000}
\definecolor{jet5}{rgb}{0.0000,1.0000,1.0000}
\definecolor{jet6}{rgb}{0.3333,1.0000,0.6667}
\definecolor{jet7}{rgb}{0.6667,1.0000,0.3333}
\definecolor{jet8}{rgb}{1.0000,1.0000,0.0000}
\definecolor{jet9}{rgb}{1.0000,0.6667,0.0000}
\definecolor{jet10}{rgb}{1.0000,0.3333,0.0000}
\definecolor{jet11}{rgb}{1.0000,0.0000,0.0000}

\definecolor{mygreen}{RGB}{111, 113, 88}

\usepackage[most]{tcolorbox}
\usepackage{cleveref}

\makeatletter
\crefformat{tcb@cnt@mytheo}{theorem~#2#1#3}
\Crefformat{tcb@cnt@mytheo}{Theorem~#2#1#3}
\crefformat{tcb@cnt@mylemma}{lemma~#2#1#3}
\Crefformat{tcb@cnt@mylemma}{Lemma~#2#1#3}
\crefformat{tcb@cnt@mycorollary}{corollary~#2#1#3}
\Crefformat{tcb@cnt@mycorollary}{Corollary~#2#1#3}
\crefformat{tcb@cnt@myproto}{prototype problem~#2#1#3}
\Crefformat{tcb@cnt@myproto}{Prototype Problem~#2#1#3}
\crefformat{tcb@cnt@myalg}{algorithm~#2#1#3}
\Crefformat{tcb@cnt@myalg}{Algorithm~#2#1#3}
\makeatother

\newtcbtheorem[auto counter, number within=section]{mytheo}{Theorem}%
{colback=EmoryBlue!5,colframe=EmoryBlue, fonttitle=\bfseries}{thm}

\newtcbtheorem[auto counter, number within=section]{mylemma}{Lemma}%
{colback=gray!5,colframe=gray, fonttitle=\bfseries}{lem}

\newtcbtheorem[auto counter, number within=section]{mycorollary}{Corollary}%
{colback=gray!5,colframe=gray, fonttitle=\bfseries}{cor}

\newtcbtheorem[auto counter, number within=section]{myproto}{Prototype Problem}%
{colback=EmoryGold!5,colframe=EmoryGold, fonttitle=\bfseries}{proto}

\newtcbtheorem[auto counter, number within=section]{myalg}{Algorithm}%
{colback=white!5,colframe=black, fonttitle=\bfseries}{alg}

\DeclareMathOperator*{\rank}{rank}
\DeclareMathOperator*{\Trank}{t-rank}
\DeclareMathOperator{\myImplicitRank}{implicit-rank}
\DeclareMathOperator*{\diag}{diag}
\DeclareMathOperator{\myGrad}{grad}
\DeclareMathOperator*{\myHess}{Hess}

\DeclareMathOperator*{\myTriu}{triu}
\DeclareMathOperator*{\myBidiag}{bidiag}
\DeclareMathOperator*{\myTril}{tril}
\DeclareMathOperator*{\myCirc}{circ}
\DeclareMathOperator*{\argmin}{argmin}

\DeclareMathOperator*{\myVec}{vec}

\DeclareMathOperator{\myFold}{fold}

\DeclareMathOperator{\myBdiag}{bdiag}
\DeclareMathOperator*{\subjectto}{s.t.}

\DeclareMathOperator{\retraction}{Retr}

\newcommand{\starM}{\star_{\bfM}}

\newcommand{\red}[1]{{\color{red} #1}}
\newcommand{\blue}[1]{{\color{blue} #1}}

\newcommand{\cyan}[1]{{\color{cyan} #1}}

\newcommand{\orange}[1]{{\color{orange} #1}}
\newcommand{\violet}[1]{{\color{violet} #1}}
\newcommand{\magenta}[1]{{\color{magenta} #1}}

\def\mydefb#1{\expandafter\def\csname bf#1\endcsname{\mathbf{#1}}}
\def\mydefallb#1{\ifx#1\mydefallb\else\mydefb#1\expandafter\mydefallb\fi}
\mydefallb aAbBcCdDeEfFgGhHiIjJkKlLmMnNoOpPqQrRsStTuUvVwWxXyYzZ\mydefallb

\def\mydefb#1{\expandafter\def\csname #1bb\endcsname{\mathbb{#1}}}
\def\mydefallb#1{\ifx#1\mydefallb\else\mydefb#1\expandafter\mydefallb\fi}
\mydefallb aAbBcCdDeEfFgGhHiIjJkKlLmMnNoOpPqQrRsStTuUvVwWxXyYzZ\mydefallb

\def\mydefb#1{\expandafter\def\csname #1cal\endcsname{\mathcal{#1}}}
\def\mydefallb#1{\ifx#1\mydefallb\else\mydefb#1\expandafter\mydefallb\fi}
\mydefallb aAbBcCdDeEfFgGhHiIjJkKlLmMnNoOpPqQrRsStTuUvVwWxXyYzZ\mydefallb

\def\mydefb#1{\expandafter\def\csname T#1\endcsname{\boldsymbol{\mathcal{#1}}}}
\def\mydefallb#1{\ifx#1\mydefallb\else\mydefb#1\expandafter\mydefallb\fi}
\mydefallb aAbBcCdDeEfFgGhHiIjJkKlLmMnNoOpPqQrRsStTuUvVwWxXyYzZ\mydefallb

\def\mydefgreek#1{\expandafter\def\csname bf#1\endcsname{\text{\boldmath$\mathbf{\csname #1\endcsname}$}}}
\def\mydefallgreek#1{\ifx\mydefallgreek#1\else\mydefgreek{#1}%
   \lowercase{\mydefgreek{#1}}\expandafter\mydefallgreek\fi}
\mydefallgreek {alpha}{Alpha}{beta}{Beta}{gamma}{Gamma}{delta}{Delta}{epsilon}{Epsilon}{zeta}{Zeta}{eta}{Eta}{theta}{Theta}{iota}{Iota}{kappa}{Kappa}{lambda}{Lambda}{mu}{Mu}{nu}{Nu}{omicron}{Omicron}{pi}{Pi}{rho}{Rho}{sigma}{Sigma}{tau}{Tau}{upsilon}{Upsilon}{phi}{Phi}{xi}{Xi}{chi}{Chi}{psi}{Psi}{omega}{Omega}\mydefallgreek


%% file: sections/S_00_introduction.tex
\section{Introduction}

The supplementary material provides thorough details that provide intuition and support theory presented in the manuscript. 
In~\Cref{tab:supplementary_material_organization}, we describe the organization of the supplementary material sections and the correspondence with sections in the manuscript. 

\begin{table}
\centering
\tiny
\caption{Organization of supplementary material sections and corresponding to manuscript sections.}
\label{tab:supplementary_material_organization}
\renewcommand{\arraystretch}{1.5}
\begin{tabular}{|l|l|}
\hline
\bf Supplementary Material & \bf Manuscript\\
\hline\hline
\Cref{sec:underlyingStructure}:  & \Cref{sec:underlyingAlgebra}:\\
 \nameref{sec:underlyingStructure} &  \nameref{sec:underlyingAlgebra}\\
\hline
\Cref{sec:riemannian_optimization_background}: & \Cref{sec:optimizationOnManifold}:\\
 \nameref{sec:riemannian_optimization_background} &   \nameref{sec:optimizationOnManifold}\\
\hline
\Cref{sec:fundamental_derivatives}: & \Cref{sec:starMDerivatives}: \\
 \nameref{sec:fundamental_derivatives}  &  \nameref{sec:starMDerivatives}\\
\hline
\Cref{app:svdDerivation}: & \Cref{sec:tsvdmDerivatives}: \\
 \nameref{app:svdDerivation}  &  \nameref{sec:tsvdmDerivatives}\\
\hline
\Cref{sec:objectiveInvariance}: & \Cref{sec:uniqueInvariancePrototypeShort}:\\
 \nameref{sec:objectiveInvariance} &  \nameref{sec:uniqueInvariancePrototypeShort}\\
\hline
\Cref{sec:boundedOperatorExamples}: & \Cref{sec:analysis}:\\
 \nameref{sec:boundedOperatorExamples} &  \nameref{sec:analysis}\\
\hline
\Cref{sec:hessianBound}: & \Cref{thm:lipschitzContinuity}:\\
 \nameref{sec:hessianBound} &  \nameref{thm:lipschitzContinuity}\\
\hline
\Cref{sec:leastSquares}: & \Cref{exam:varproAD}:\\
 \nameref{sec:leastSquares} &  \nameref{exam:varproAD}\\
\hline
\Cref{app:varproConstrained}: & \Cref{sec:portfolio}:\\
 \nameref{app:varproConstrained} &  \nameref{sec:portfolio}\\
\hline
\Cref{app:stocks}: & \Cref{sec:portfolio}:\\
 \nameref{app:stocks} &  \nameref{sec:portfolio}\\
\hline
\Cref{sec:wave_equation}: & \Cref{sec:rom}:\\
 \nameref{sec:wave_equation} & \nameref{sec:rom}\\
\hline
\end{tabular}
\end{table}

\subsection{Outline of Supplementary Material}

The supplementary material is organized as follows. 
In~\Cref{sec:underlyingStructure}, we present concrete examples of tensor algebras for various choices of transformation $\bfM$. 
In~\Cref{sec:riemannian_optimization_background}, we provide some intuition for Riemannian optimization. 
In~\Cref{sec:fundamental_derivatives}, we provide formulas and derivations for fundamental tensor operations. 
In~\Cref{app:svdDerivation}, we differentiate through the $t$-SVDM.

%% file: sections/S_01_underlying_algebra.tex
\section{Underlying Tensor Algebra Structure}
\label{sec:underlyingStructure}

The underlying algebra of the $\starM$-framework is defined by $\bfR_{\bfM}[\cdot]$ in~\eqref{eqn:starMtubalVectorized}.
We identify some simple cases where we can write the algebra explicitly in~\Cref{tab:underlyingStructure}. 
While this list is incomplete, it demonstrates that the choice of algebra can significantly alter tubal multiplication. 

\begin{table}
\centering
\caption{Examples of tensor algebras for various choices of transformations. 
The matrix subalgebra changes significantly based on the transformation. 
DFT stands for discrete Fourier transform. 
Here, $\mathbf{1}$ is the constant vector of all ones, $\myTril(\cdot, k)$ and $\myTriu(\cdot, k)$ are the lower/upper triangular parts of a matrix up to the sub/super-diagonal $k$ where $k = 0$ is the main diagonal. }
\label{tab:underlyingStructure}

\scriptsize 
\begin{tabular}{lccc}
Name & $\bfM$ & $\bfR_{\bfM}[\bfa]$ & $3\times 3$ case\\
\hline &&&\\[0.25em]
Identity & $\bfI$ & $\diag(\bfa)$ & $\tiny \begin{bmatrix} a_1 \\ & a_2 \\ && a_3 \end{bmatrix}$\\[2em]

DFT & $\bfF$ & $\myCirc(\bfa)$ & $\tiny \begin{bmatrix} a_1 & a_3 & a_2 \\ a_2 & a_1 & a_3 \\ a_3 & a_2 & a_1 \end{bmatrix}$\\[2em]

Summation & $\bfS = \myTril(\mathbf{1}\mathbf{1}^\top)$ & $\myTril(\bfa \mathbf{1}^\top,-1) + \diag(\bfS \bfa)$ & $\tiny \begin{bmatrix} a_1 & 0 & 0 \\ a_2 & a_1 + a_2 & 0 \\ a_3 & a_3 & a_1 + a_2 + a_3 \end{bmatrix}$\\[2em]

Finite Difference & $\bfD = \myBidiag(1,-1)$ & $2\myTriu(\bfb \mathbf{1}\bfa^\top) - \diag(\bfa)$ & $\tiny \begin{bmatrix} a_1 & 2a_2 & 2 a_3 \\ 0 & a_2 & 2 a_3 \\ 0 & 0 & a_3 \end{bmatrix}$

\end{tabular}
\end{table}

%% file: sections/S_09_riemannian_optimization.tex
\section{Riemannian Optimization Intuition}
\label{sec:riemannian_optimization_background}

The following derivation of Riemannian gradient descent on the orthogonal group $\Ocal_{n_3}$ is based on~\cite{Absil2008,  boumal2023intromanifolds, EdelmanAriasSmith1998}. 
Given a Euclidean search direction $\bfS \in \Rbb^{n_3\times n_3}$ (e.g., Euclidean gradient), we decompse the search direction into two components in Euclidean space $\bfS = \bfS_T + \bfS_N$ where $\bfS_T$ lies tangent to the manifold at the current iterate $\bfM$ and $\bfS_N$ lies in the normal direction.   
For $\Ocal_{n_3}$, the tangent space at the current iterate is
	\begin{align}\label{eq:tangentSpace}
		T_{\bfM} \Ocal_{n_3} := \left\{\bfM\bfOmega \mid \bfOmega\in \Rbb^{n_3\times n_3} \text{ and } \bfOmega = -\bfOmega^\top\right\}.
	\end{align}
One intuition comes from differentiating the constraint $\bfM^\top \bfM = \bfI$ with respect to $\bfM$.    
This yields the expression $(\delta \bfM)^\top \bfM + \bfM^\top (\delta \bfM) = \bf0$ where $\delta \bfM$ lies in the tangent space at the point $\bfM$. 
Because $\bfM$ is invertible, we can safely set $\delta \bfM = \bfM\bfOmega$ for some matrix $\bfOmega$.  
Plugging this form for $\delta \bfM$ into the differentiated constraint, we find that $\bfOmega$ is antisymmetric, as desired.     

\begin{figure}
\centering

\includegraphics[width=0.5\linewidth]{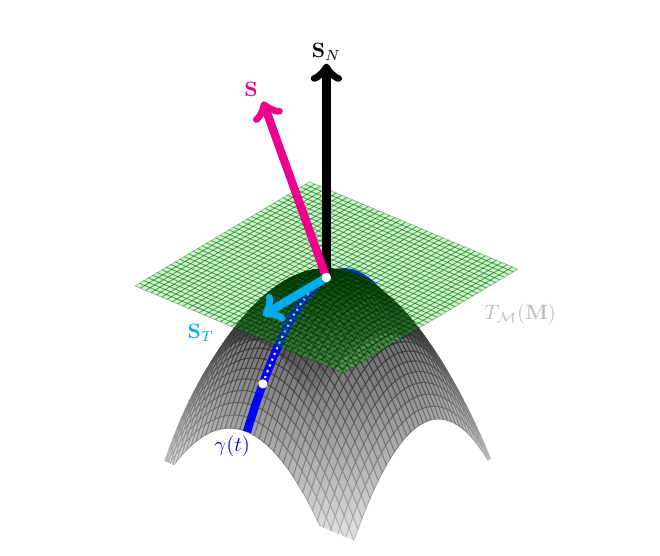}

\caption{Illustration of optimization on a manifold.  
We update our current iterate (white dot at top) along the manifold ({\color{gray}{gray}}) following the geodesic $\gamma(t)$ ({\bf \color{blue}{blue}}). 
From an implementation perspective, we first compute the Euclidean gradient ({\bf \magenta{magenta}} arrow $\bfS$) and decompose it into a tangent direction ({\bf \cyan{cyan}} arrow {$\bfS_T$} on the {\bf \color{mygreen}{green}} tangent bundle {$T_\Mcal(\bfM)$}) and a normal direction ({\bf black} arrow $\bfS_N$).  
We step along the tangent direction and then retract onto the manifold (large, lower white dot). 
}
\label{fig:manifoldOptimization}

\end{figure}

To find the tangent  component of $\bfS$, we first represent the search direction $\bfS$ as 
	\begin{align}
	\bfS = \bfM (\bfOmega + \bfW) \qquad \Longrightarrow \qquad  \bfM^\top \bfS = \bfOmega + \bfW
	\end{align}
where $\bfOmega$ is antisymmetric and $\bfW$ is symmetric.  
We rewrite $\bfM^\top \bfS$ as the sum of an antisymmetric and symmetric matrix via
	\begin{align}
	\bfM^\top \bfS = \underbrace{\frac{\bfM^\top \bfS - \bfS^\top \bfM}{2}}_{\bfOmega} +  \underbrace{\frac{\bfM^\top \bfS + \bfS^\top \bfM}{2}}_{\bfW}.
	\end{align}
We define the Riemannian gradient by projecting $\bfS$ onto the tangent space via 
	\begin{align}\label{eqn:projection2}
	\myGrad \overline{\Phi}(\bfM) &:= \bfM\left(\frac{\bfM^\top \bfS - \bfS^\top \bfM}{2}\right)
	\end{align}
We return to the manifold using a retraction, specifically the exponential mapping
	\begin{align}\label{eq:exponentialMapping3}
	\bfM \gets \retraction_{\bfM}(-\alpha \myGrad \overline{\Phi}(\bfM))  \qquad \text{where} \qquad \retraction_{\bfM}(\bfM\bfOmega) = \bfM \exp(\bfOmega)
	\end{align}
where $\bfOmega \in T_{\bfM} \Ocal_{n_3}$ and $\exp: \Rbb^{n_3\times n_3} \to \Rbb^{n_3\times n_3}$ is a matrix function~\cite{Higham2008:matrixFunctions}. 
We note that many other retractions are possible, such as Cayley transformation and QR factorization~\cite{Absil2008}.  
We focus on the exponential mapping in this work and leave other mappings for future research.  
We illustrate manifold optimization in~\Cref{fig:manifoldOptimization}.

%% file: sections/S_08_fundamental_derivatives.tex
\section{Fundamental $\starM$-Product Derivatives} 
\label{sec:fundamental_derivatives}

We assume all tensors and matrices are of compatible sizes for the given operations. 
We denote the direction in which we apply a derivative with respect to a variable $z$ by $\delta z$, which has the same size as the original variable $z$.  
For a scalar- or vector-valued function $f: \Ucal \to \Vcal$ mapping from input space $\Ucal$ to output space $\Vcal$, we will denote the Euclidean gradient as $\nabla f: \Vcal \to \Ucal$.  

The key derivatives we require to understand subsequent formulations are the derivative of the mode-$3$ product and the $\starM$-product.   

\begin{mylemma}{mode-$3$ product derivatives}{mode3Grad}
Given $\TB := f(\TA,\bfM) = \TA\times_3 \bfM$, the directional derivatives of the mode-$3$ product are given by 
\begin{subequations}
	\begin{align}
		\nabla_{\TA} f(\delta\TB)  &= \delta \TB \times_3 \bfM^\top   \label{eq:mode3gradA}\\
		\nabla_{\bfM} f(\delta \TB) &=(\delta \TB)_{(3)} (\TA_{(3)})^\top \label{eq:mode3gradM}
	\end{align}
\end{subequations}
\end{mylemma}

\begin{proof}
The proof follows from expressing $\TA \times_3 \bfM = \myFold_3(\bfM \TA_{(3)})$. 
See~\cite{NewmanTNN2018} for details.
\end{proof}

 \begin{mylemma}{$\starM$-product derivatives}{starMGrad}
Given $\TC := f(\TA,\TB,\bfM) = \TA \starM \TB$, the directional derivatives of the $\starM$-product are
\begin{subequations}
	\begin{align}
		\nabla_{\TA} f(\delta \TC) 
			&=   \delta \TC \starM \TB^\top \label{eq:starMGradA} \\
		\nabla_{\TB} f(\delta \TC) 
			&=  \TA^\top \starM \delta\TC \label{eq:starMGradB} \\
			\begin{split}
		\nabla_{\bfM} f(\delta \TC) 
			&=\bfM \left[\TC_{(3)} (\delta\TC_{(3)})^\top   \right.\\
			&\qquad +\left.\left(\nabla_{\TA} f(\delta \TC) \right)_{(3)} (\TA_{(3)})^\top+ \left(\nabla_{\TB} f(\delta \TC) \right)_{(3)} (\TB_{(3)})^\top\right] \label{eq:starMGradM}
		\end{split}
	\end{align}
\end{subequations}
\end{mylemma}

\begin{proof}
The proof follows from expressing the $\starM$-product as
	\begin{align}
	\TC = \TA \starM \TB = \left[\widehat{\TA} \smalltriangleup \widehat{\TB}\right] \times_3 \bfM^\top,
	\end{align}
where $\widehat{\TY} = \TY \times_3 \bfM$.  See~\cite{NewmanTNN2018} for details of computing~\eqref{eq:starMGradA} and~\eqref{eq:starMGradB}.  
We derive~\eqref{eq:starMGradM} as follows. 
From~\Cref{lem:mode3Grad}, the first term is the derivative with respect to $\bfM^\top$, given by the transpose of~\eqref{eq:mode3gradM}.  
The second and third terms respectively combine~\eqref{eq:starMGradA} and~\eqref{eq:starMGradB} with~\eqref{eq:mode3gradM}. 

\end{proof}

%% file: sections/S_02_tsvdm_derivatives.tex
\section{Differentiating through the $t$-SVDM}
\label{app:svdDerivation}

\def\dU{\delta \bfU}
\def\dS{\delta \bfSigma}
\def\dV{\delta \bfV}
\def\dA{\delta \bfA}


The goal of this section is to compute the derivatives through the truncated $t$-SVDM. 
We start by presenting the reverse mode derivatives of a matrix SVD based on~\cite{Townsend2016, Giles1948, Wan2019:complexSVD, minka1997old}.  
Let $r = \min(n_1,n_2)$ and let $\bfA = \bfU \bfSigma \bfV^\top$ be the economic SVD where $\bfU \in \Rbb^{n_1\times r}$ and $\bfV\in \Rbb^{n_2\times r}$ have orthonormal columns and $\bfSigma\in \Rbb^{r\times r}$ is a diagonal matrix with real, nonnegative diagonal entries in decreasing order.  
First, we perturb the factor matrices by a small amount, e.g., $\bfU + \varepsilon\dU$ where $\varepsilon>0$ is small, the matrix $\bfA$ and the orthogonality constraints are perturbed accordingly via
	\begin{align}
	\begin{split}
	\bfA + \varepsilon\dA+ \Ocal(\varepsilon^2) 
		&= (\bfU + \varepsilon\dU)(\bfSigma + \varepsilon\dS)(\bfV + \varepsilon\dV^\top)\\
		&=\bfU\bfSigma\bfV^\top +\varepsilon\left( \dU \bfSigma \bfV^\top + \bfU \dS \bfV^\top + \bfU \bfSigma \dV^\top\right) + \Ocal(\varepsilon^2).
	\end{split}
	\end{align}
Then, the gradients are
	\begin{subequations}\label{eq:svdFactorDerivatives}
	\begin{alignat}{2}
	&\nabla_{\bfA} \bfU(\dU) &&=[\bfU (\bfF \odot (\bfU^\top \dU - \dU^\top \bfU)) \bfSigma + (\bfI_{n_1} - \bfU \bfU^\top) \dU \bfSigma^{-1}] \bfV^\top\\
	&\nabla_{\bfA} \bfSigma_k(\dS) &&= \bfU(\bfI_r \odot \dS) \bfV^\top\\
	&\nabla_{\bfA} \bfV(\dV) &&= \bfU [\bfSigma (\bfF \odot (\bfV^\top\dV- \dV^\top \bfV) \bfV^\top + \bfSigma^{-1}\dV (\bfI_{n_2} - \bfV\bfV^\top)]
	\end{alignat}
	\end{subequations}

%% file: sections/S_02_02_tsvdm_derivatives.tex
\subsection{$t$-SVDM Derivative Formula}
\label{app:tsvdmDerivation}
We leverage the matrix SVD derivative in~\Cref{app:svdDerivation} to compute the $t$-SVDM version. 
Suppose $\TA \in \Rbb^{n_1 \times n_2 \times n_3}$ with $t$-SVDM $\TA = \TU \starM \TS \starM \TV^\top$ for $\bfM\in \Ocal_{n_3}$.   
Let $\TX = \TU_k \starM \TS_k \starM \TV_k^\top$ be the truncation $t$-SVDM with $t$-rank $k$.  
Our goal is to differentiate through the $t$-SVDM with respect to $\bfM$.  
Recall, computing the (truncated) $t$-SVDM requires the following steps: 
	\begin{enumerate}[label=(\arabic* $\downarrow$)]
	 \item Move $\TA$ into the transform domain: $\widehat{\TA} = \TA\times_3 \bfM$
	 \item Compute matrix SVDs: $\widehat{\TA}_{:,:,i} = \widehat{\TU}_{:,:,i}  \widehat{\TS}_{:,:,i} \widehat{\TV}_{:,:,i}^\top$ for $i=1,\dots,n_3$.
	 \item Truncate the factors: $\widehat{\TA}_k = \widehat{\TU}_{:,1:k,:} \smalltriangleup \widehat{\TS}_{1:k,1:k,:}\smalltriangleup\widehat{\TV}_{:,1:k,:}^\top$
	 \item Return to the spatial domain: $\TA_k = \widehat{\TA}_k \times_3 \bfM^\top$
	\end{enumerate}

Computing the $t$-SVDM derivative with respect to $\bfM$ in reverse mode requires differentiation in the opposite order, starting with step (4):
\begin{enumerate}[label=(\arabic*)]
\item Differentiate through the $\TA_k$ using~\eqref{eq:mode3gradA}.
\item Pad with zeros to reverse the truncation. 
\item Differentiate through the facewise product of the factors using~\eqref{eq:starMGradA} and~\eqref{eq:starMGradB}. 
\item Differentiate through $\widehat{\TA}_k$ using~\eqref{eq:svdFactorDerivatives} for each frontal slice. 
This results in~\eqref{eq:tsvdmFactorDerivatives}. 
\item Differentiate through the mode-$3$ product using~\eqref{eq:mode3gradA}.
\end{enumerate}

In our current implementation, we compute the derivative of the full SVD and zero out the directions pertaining to the small singular values. 
This seemed to give more robust performance in our algorithm than the direct low-rank versions of the derivatives. 
We will focus on low-rank implementations in future versions of the code as part of the push to improve efficiency.

%% file: sections/S_03_prototype_invariance.tex
\section{Invariance of Prototype Objective Functions}
\label{sec:objectiveInvariance}

Building from the invariant properties of the $\starM$-product in~\Cref{sec:invariantProduct}, we identify modifications of $\bfM$ that will not change the prototype objective function value. 
By~\Cref{lem:starMInvariance} and~\Cref{thm:starMInvariance}, both prototype objective functions are invariant to negative left multiplication of a permutation matrix. 

\begin{mytheo}{Invariance of $\overline{\Phi}_{\rm reg}$}{tlinearRegressionInvariance2}
The reduced objective function $\overline{\Phi}_{\rm reg}$ is invariant to permutations and negations of the rows of $\bfM$. \end{mytheo}

\begin{proof} 
The proof for permutation of the rows follows from~\Cref{thm:starMInvariance}.    
The proof of negation of the rows proceeds as follows.

Let $\bfD$ be a diagonal matrix and every diagonal entry is either $+1$ or $-1$. 
First, we observe that $\TY \star_{\bfD} \TZ = (\TY \smalltriangleup \TZ) \times_3 \bfD$. 
This means the $\bfD$ affects only the sign of the frontal slices under the facewise product. 
Using this observation, we show that the $t$-linear regression solutions are equal \emph{in the transform domain} for either transformation; that is, $\TX_{\rm reg}(\bfD\bfM) \times_3 (\bfD \bfM) = \widehat{\TX}_{\rm reg}(\bfM)$ where $\widehat{\TY} = \TY \times_3 \bfM$, 
We start from the normal equations with the $\star_{\bfD\bfM}$-product and use similar logic as in~\eqref{eq:productOfMatricesMprod} to write the product in the $\starM$-transform domain:
	\begin{subequations}
	\begin{align}
	 (\TA^\top \star_{\bfD\bfM} \TA) \star_{\bfD\bfM} \TX &= \TA^\top \star_{\bfD\bfM} \TB\\
	[(\widehat{\TA}^\top \smalltriangleup \widehat{\TA}) \smalltriangleup\widehat{\TX}] \times_3 \bfD\bfM^\top
		&=[\widehat{\TA}^\top \smalltriangleup \widehat{\TB}] \times_3 \bfD\bfM^{\top},
	\end{align}
	\end{subequations}
Because $\bfD\bfM^\top$ is orthogonal, we can cancel the mode-$3$ multiplication on both sides.  
Thus, the optimal solution in the transform domain depends only on $\bfM$, not on $\bfD$, which in turn means the transform-domain solutions are under the $\star_{\bfD\bfM}$- and $\starM$-products.   
Let $\widehat{\TX}_{\rm reg}^* = \widehat{\TX}_{\rm reg}(\bfM) = \TX_{\rm reg}(\bfD\bfM) \times_3 (\bfD\bfM)^\top$. 

We use the equivalence of solutions in the transform domain to show $\overline{\Phi}_{\rm reg}(\bfD\bfM) = \overline{\Phi}_{\rm reg}(\bfM)$.  
Starting from the $\bfD\bfM$ case, we obtain
	\begin{align}
	\begin{split}
	\overline{\Phi}_{\rm reg}(\bfD\bfM)  
		&= \tfrac{1}{2}\|\TA \star_{\bfD\bfM} \TX_{\rm reg}(\bfD\bfM) - \TB\|_F^2 \\
		&= \tfrac{1}{2} \|(\widehat{\TA} \smalltriangleup  \widehat{\TX}_{\rm reg}^*)  \times_3 \bfD \bfM^\top - \widehat{\TB} \times_3 \bfD \bfM^\top \|_F^2.
	\end{split}
	\end{align}
Using the orthogonal invariance of the Frobenius norm, we have
	\begin{align}
	\begin{split}
	\overline{\Phi}_{\rm reg}(\bfD\bfM)  & =  \tfrac{1}{2} \|(\widehat{\TA} \smalltriangleup \widehat{\TX}_{\rm reg}^*) - \widehat{\TB} \|_F^2 \\
	&=   \tfrac{1}{2} \|({\TA} \starM {\TX}_{\rm reg}(\bfM)) - {\TB} \|_F^2\\
	&= \overline{\Phi}_{\rm reg}(\bfM).
	\end{split}
	\end{align}

\end{proof}

\begin{mytheo}{Invariance of $\overline{\Phi}_{\rm low}$}{lowRankInvariance}
The reduced low-$t$-rank objective function $\overline{\Phi}_{\rm low}$ is invariant to negation and row permutations of $\bfM$.  
\end{mytheo}

\begin{proof} 
The proof follows from~\Cref{lem:starMInvariance},~\Cref{thm:starMInvariance}, and~\Cref{thm:eckartYoung}.  Specifically, the solution, $\TX_{\rm low}(\bfM) = \TU_{:,1:k,:} \starM \TS_{1:k,1:k,:} \starM \TV_{:,1:k,:}^\top$ requires two $\starM$-products, which are invariant under negation and permutations. 
\end{proof}


\subsection{Uniqueness of the Prototype Problem Representations}
\label{sec:uniqueInvariancePrototype}

Each prototype problem offers uniqueness properties of the representation $\TX(\bfM)$ and invariance to modifications of $\bfM$. 
We require a few additional $\starM$-based definitions for these results. 

\begin{definition}[$\starM$-invertible]\label{def:starMinvertible} 
A tensor $\TA \in \Rbb^{n_1 \times n_1 \times n_{3}}$ if $\starM$-invertible if there exists a tensor $\TB\in  \Rbb^{n_1 \times n_1 \times n_{3}}$ such that $\TA \starM \TB = \TB \starM \TA = \TI$.  
In practice, we compute the inverse in the transform domain with
	\begin{align*}
	\widehat{\TB}_{:,:,i} = (\widehat{\TA}_{:,:,i})^{-1} \qquad \text{for $i=1,\dots,n_3$}.
	\end{align*}
We denote the inverse with $\TB = \TA^{-1}$.  
\end{definition}

Note that the $t$-rank captures the number of singular tubes, but each tube requires storage of up to $n_3$ nonzero entries. 
We define a different notion of the rank to quantify the total number of nonzero singular values needed to form the $t$-SVDM:
\begin{definition}[$\starM$-implicit rank]\label{def:implicitRank} 
The implicit rank of a tensor $\TA$ is the sum of the ranks of each frontal slice in the transform domain; that is, for $\widehat{\TA} = \TA \starM \bfM$, we have
	\begin{align}
	\myImplicitRank(\TA, \bfM) = \sum_{k=1}^{n_3} \rank(\widehat{\TA}_{:,:,k}).
	\end{align}
\end{definition}

We describe key properties below. 

\begin{mytheo}{Uniqueness of $\TX_{\rm reg}(\bfM)$}{tLinearRegUnique}
Given model tensor $\TA\in \Rbb^{n_1\times p\times n_3}$ with $n_1\le p$ and  data $\TB\in \Rbb^{n_1\times n_2\times n_3}$,  fix $\bfM\in \Mcal$.  
If $\TA$ has full implicit rank, then $\TX_{\rm reg}(\bfM)$ is unique.
\end{mytheo}

\begin{proof}
Solving the $t$-linear regression problem is equivalent solving the $\starM$-analog of the normal equations
	\begin{align}\label{eq:tensorNormalEquations}
	(\TA^\top \starM \TA) \starM \TX = \TA^\top \starM \TB.
	\end{align}
The tensor  $\TA^\top \starM \TA$ is invertible if each frontal slice in the transform domain is invertible; i.e., if $\widehat{\TA}_{:,:,i}^\top \widehat{\TA}_{:,:,i}$ is  invertible for $i=1,\dots, n_3$.  
Because $\TA$ has full implicit rank (\Cref{def:implicitRank}) by assumption, each $\widehat{\TA}_{:,:,i}$ has full matrix rank, and hence each matrix  $\widehat{\TA}_{:,:,i}^\top \widehat{\TA}_{:,:,i}$ is invertible. 
Thus, the solution is uniquely given by 
	\begin{align}
	\TX_{\rm reg}(\bfM) = (\TA^\top \starM \TA)^{-1} \starM \TA^\top \starM \TB
	\end{align}    
Note that if $\TA$ does not have full implicit rank for a given transformation $\bfM$, there are infinitely many solutions, $\TX_{\rm reg}(\bfM)$, to the $t$-linear regression problem. 
In this case, we use the {\sc Matlab} default of selecting the solution for $\TX_{\rm reg}(\bfM)$ that has the greatest number of zeros in the transform domain\footnote{See~\url{https://www.mathworks.com/help/dsp/ref/qrsolver.html} for details.}. 
\end{proof}

\begin{mytheo}{Uniqueness of $\TX_{\rm low}(\bfM)$}{lowRankUnique}
Consider a tensor $\TA\in \Rbb^{n_1\times n_2\times n_3}$.  
The $t$-SVDM $\TA = \TU \starM \TS \starM \TV^\top$ is formed by computing matrix SVDs in the transform domain; that is, 
	\begin{align}
	\widehat{\TA}_{:,:,i} = \widehat{\TU}_{:,:,i} \widehat{\TS}_{:,:,i} \widehat{\TV}_{:,:,i}^\top \qquad \text{for $i=1,\dots,n_3$.}
	\end{align}
Thus, the $t$-SVDM is unique up to properties of the matrix SVD for each frontal slices in the transform domain. 
\end{mytheo}

\begin{proof} 
The matrix SVD $\bfA = \bfU \bfSigma \bfV^\top$ is unique up to sign changes (i.e., $\bfU\bfD$ and $\bfV\bfD^{-1}$ where $\bfD$ is a diagonal matrix with $\pm 1$ on the diagonal) and to orthogonal transformations of subspaces spanned by singular vectors corresponding to the same singular value~\cite[Section 7.1.3]{GoluVanl96}.  
The $t$-SVDM has the same uniqueness properties for each frontal slice in the transform domain. 
\end{proof}

%% file: sections/S_04_bounded_operator_norm.tex
\section{Boundedness of the $\starM$-Operator Norms}
\label{sec:boundedOperatorExamples}

A tensor $\TA$ has a bounded $\starM$-operator norm if $\|\TA\| \le C$ for all $\bfM\in \Mcal$ where $C$ is a constant independent of $\bfM$.  
We present two simple tensors and the corresponding $\starM$-pseudoinverses (\Cref{def:starMpseudo}), one that has an unbounded operator norm (\Cref{def:operatorNorm}) and one that has a bounded operator norm.

\begin{example}[Unbounded $\starM$-pseudoinverse operator norm]Consider the $2\times 2\times 2$ tensor
	\begin{align}
	\TA_{:,:,1} &= \begin{bmatrix} 1 & 0 \\ 0 & 1 \\ \end{bmatrix} 
	& \TA_{:,:,2} &= \begin{bmatrix} 1 & 0 \\ 0 & 0 \\ \end{bmatrix} 
	\end{align}
For $\bfM = \bfI$, we have $\|\TA\| = 1$.
Suppose $\bfM$ is a $2\times 2$ rotation matrix
	\begin{align}
	\bfM(\theta) = \begin{bmatrix} \cos(\theta) & -\sin(\theta)\\ \sin(\theta) & \cos(\theta) \end{bmatrix}.
	\end{align}
Then, in the transform domain, we have
	\begin{align}
	\widehat{\TA}_{:,:,1} &= \begin{bmatrix} \cos(\theta) - \sin(\theta) & 0 \\ 0 & \cos(\theta) \\ \end{bmatrix} 
	& \widehat{\TA}_{:,:,2} &= \begin{bmatrix} \cos(\theta) - \sin(\theta) & 0 \\ 0 & \sin(\theta) \\ \end{bmatrix}.
	\end{align}
The $\starM$-pseudoinverse operator norm is unbounded; specifically, $\|\TA^\dagger\| = \frac{1}{\sin \theta} \to \infty$ as $\theta \downarrow 0$. 
\end{example}

\begin{example}[Bounded $\starM$-pseudoinverse operator norm] Consider the $2\times 2 \times 2$ tensor
	\begin{align}
	\TA_{:,:,1} &= \begin{bmatrix} 1 & 1 \\ -1 & 1 \\ \end{bmatrix} 
	& \TA_{:,:,2} &= \begin{bmatrix} 1 & -1 \\ 1 & 1\\ \end{bmatrix}.
	\end{align} 
For any orthogonal $\bfM \in \Mcal$, we have
	\begin{align} 
	\widehat{\TA}_{:,:,i} = (m_{i1} + m_{i2})\begin{bmatrix} 1 & 0\\ 0 & 1 \end{bmatrix} + (m_{i1} - m_{i2}) \begin{bmatrix} 0 & 1\\ -1 & 0 \end{bmatrix}
	\end{align}
where $m_{ij} = \bfM_{ij}$ for $i,j=1,2$. 
The singular values of the transformed frontal slices are constant with 
	\begin{align}
	\widehat{\sigma}_{j,j,i} = \sqrt{(m_{i1} + m_{i2})^2 + (m_{i1} - m_{i2})^2} = \sqrt{2}. 
	\end{align}
for $j =1,2$. 
Note that we use the orthogonality of $\bfM$ for the simplification. 
Thus, $\|\TA^\dagger\| = \frac{1}{\sqrt{2}}$, which is constant for any choice of $\bfM\in \Mcal$. 

\end{example}

%% file: sections/S_07_bounded_hessian.tex
\section{Bounded Hessian for Reduced $t$-Linear Regression}
\label{sec:hessianBound}

This section serves as a detailed proof of~\Cref{thm:lipschitzContinuity}.  
To prove the boundedness of the Euclidean Hessian relies,  we rely on vectorizing matrix-matrix and tensor-tensor products using the Kronecker product definition from~\cite{Petersen07thematrix}. 
Specifically, given $\bfA\in \Rbb^{n_1\times p}$ and $\bfB\in \Rbb^{p\times n_2}$, we have
	\begin{align}
	\myVec(\bfA \bfB) &= (\bfB^\top \otimes \bfI_{n_1})\myVec(\bfA)	 = (\bfI_{n_2} \otimes \bfA) \myVec(\bfB) \in \Rbb^{n_1n_2}.
	\end{align}
\begin{remark}
For notational simplicity in the subsequent derivation, we remove the subscript ``reg'' and use $\overline{\Phi}$ and $\TX(\bfM)$ to denote the $t$-linear regression objective function and optimal solution, respectively. 
\end{remark}

%% file: sections/S_07_01_outline.tex
\subsection{Outline of Proof}
Our goal is to show the Euclidean Hessian, $\nabla^2 \overline{\Phi}$, is bounded in the Frobenius norm; this will show that the Riemannian Hessian, $\myHess \overline{\Phi}$, on the orthogonal group is bounded. 
As described in~\cite[Corollary 5.47, p. 111]{boumal2023intromanifolds}, the Riemannian Hessian is a projection of the Euclidean Hessian with an additional correction from Euclidean gradient information. 
Following the presentation in~\cite[Section 2.4.5]{EdelmanAriasSmith1998}, the Riemannian Hessian on the orthogonal group is
	\begin{align}
	\myHess \overline{\Phi}(\bfM) &= \nabla^2 \overline{\Phi}(\bfM) + \tfrac{1}{2}((\bfM \otimes \nabla \overline{\Phi}(\bfM)^\top) + (\nabla \overline{\Phi}(\bfM) \otimes \bfM^\top))
	\end{align}
We abuse notation in our definition slightly. 
Typically, Hessians act on two tangent directions, denoted $\bfDelta_1, \bfDelta_2 \in \Rbb^{n_3\times n_3}$. 
Here, we have defined the Hessian to act on vectorized tangents; that is,
	\begin{align}
	\begin{split}
	&\myHess \overline{\Phi}(\bfM)[\bfDelta_1, \bfDelta_2] \\
	&\qquad = \myVec(\bfDelta_1)^\top [\nabla^2 \overline{\Phi}(\bfM) + \tfrac{1}{2}((\bfM \otimes \nabla \overline{\Phi}(\bfM)^\top) + (\nabla \overline{\Phi}(\bfM) \otimes \bfM^\top))]\myVec(\bfDelta_2).
	\end{split}
	\end{align} 
The action of the Hessian on the vectorized tangents is the same as the non-vectorized case.  
This notation will be convenient to prove the boundedness of the Hessian as an operator. 
	
The Riemannian Hessian is bounded above in the Frobenius norm (which serves as an upper bound of the operator norm) if the Euclidean gradient and the Euclidean Hessian are bounded independent of $\bfM$. 
We derived the bound on the Euclidean gradient in~\Cref{lem:boundedGradient}. 
The key to bounding  the Euclidean Hessian is to differentiate through the approximation $\widetilde{\TB}(\bfM) = \TA\starM \TX(\bfM)$, vectorized with respect to the transformation $\bfM$. 
In~\Cref{defn:starMProduct}, we see that $\bfM$ appears three times in the $\starM$-product, to compute $\widehat{\TA}$ and $\widehat{\TX}(\bfM)$ and to return to the spatial domain using $\bfM^\top$. 
The approximation has one additional dependence on $\bfM$ because the optimal solution, $\TX(\bfM)$, depends on the transformation.  
Thus, the Jacobian of the vectorized approximation\footnote{In the vectorized formulation, the Jacobian is the transpose of the gradient.  In our proof, we will use the most convenient representation of the derivative based on context.} will be composed of four terms:
	\begin{align}\label{eq:approxJacobian}
	\begin{split}
	\Jcal \myVec(\widetilde{\TB}(\bfM)) &= \underset{\color{gray}\widehat{\TA}}{\underline{J_1(\TA, \TX(\bfM), \bfM)}} + \underset{\color{gray}\widehat{\TX}(\bfM)}{\underline{J_2(\TA, \TX(\bfM), \bfM)}} \\
	&\qquad+ \underset{\color{gray}\bfM^\top}{\underline{J_3(\TA, \TX(\bfM), \bfM)}} + \underset{\color{gray}\TX(\bfM)}{\underline{J_*(\TA, \TX(\bfM), \bfM)}}.
	\end{split}
	\end{align} 
The notation $J_*$ indicates differentiation through the optimal solution, $\TX(\bfM)$. 
In~\Cref{sec:boundedHessianTechnicalDerivation}, we break down the terms of the Euclidean Hessian and argue that the bound on the approximation is sufficient to show the Hessian is bounded as well.  
In~\Cref{sec:boundedHessianHelperLemmas}, we show that the four terms in~\eqref{eq:approxJacobian} are bounded independent of $\bfM$.

%% file: sections/S_07_02_bounded_gradient.tex
\subsection{Bounded Euclidean Gradient}
\label{sec:boundedGradient}

We start by showing the Euclidean gradient is bounded independent of $\bfM$. 

\begin{mylemma}{Bounded Euclidean gradient of $\overline{\Phi}_{\rm reg}$}{boundedGradient}
Let  $\TA \in \Rbb^{n_1\times n_2\times n_3}$ be a data tensor and assume its $\starM$-pseudoinverse  has a bounded $\starM$-operator norm. 
Then, the Euclidean gradient of the reduced $t$-linear regression function, $\nabla \overline{\Phi}_{\rm reg}$, is bounded.
\end{mylemma}

\begin{proof}
Using the first-order optimality condition for the full problem, 
	\begin{align}
	\nabla_{\TX} \Phi_{\rm reg}(\bfM, \TX(\bfM)) = \TA^\top \starM \TR  = \bf0, 
	\end{align}
we can simplify the Euclidean gradient in~\eqref{eq:starMGradM} in~\Cref{lem:starMGrad} by omitting the final term and obtain
	\begin{align}\label{eq:gradBoundGoal}
	\nabla \overline{\Phi}_{\rm reg}(\bfM) = \bfM\left[(\TA\starM \TX(\bfM))_{(3)} (\TR_{(3)})^\top + (\TR \starM \TX(\bfM)^\top)_{(3)} (\TA_{(3)})^\top \right].
	\end{align}
where $\TR = \TA \starM \TX(\bfM) - \TB$ depends on $\bfM$.  
Using the orthogonal invariance and submultiplicativity of the Frobenius norm and the triangle inequality, we obtain the bound
	\begin{align}
	\|\nabla \overline{\Phi}_{\rm reg}(\bfM)\|_F \le \|\TA\starM \TX(\bfM)\|_F \|\TR\|_F + \|\TR\|_F \|\TX(\bfM)\|_F \|\TA\|_F.
	\end{align}
We will show the norm of every variable -- the residual, $\TR$, the approximation, $\TA \starM \TX(\bfM)$, and the solution, $\TX(\bfM)$ -- can be bounded independent of $\bfM \in \Mcal$. 
We derive these bounds term by term below.

\begin{itemize}
\item {\bf Residual Bound:} Because $\TX(\bfM)$ is optimal, we have
	\begin{align}\label{eq:resBound}
	\|\TA \starM \TX(\bfM) - \TB\|_F \le \|\TA \starM \TX - \TB\|_F \qquad \text{for all }\TX\in \Xcal.
	\end{align}
In particular, the above inequality holds for $\TX = \bf0$, thus 
	\begin{align}
	\|\TR\|_F = \|\TA \starM \TX(\bfM) - \TB\|_F \le \|\TB\|_F.
	\end{align}

\item {\bf Approximation Bound:} Notice that $\TA \starM \TX(\bfM) = \TR + \TB$. Using the triangle inequality, we have
	\begin{align}
	\|\TA \starM \TX(\bfM) \|_F  = \|\TR + \TB\|_F \le \|\TR\|_F + \|\TB\|_F \le 2 \|\TB\|_F
	\end{align}

\item {\bf Solution Bound:} 
For each frontal slices in the transform domain, we compute the least squares solution by solving the normal equations via
	\begin{align}
	\widehat{\TX}_{:,:,i}(\bfM) = (\widehat{\TA}_{:,:,i}^\top \widehat{\TA}_{:,:,i})^\dagger \widehat{\TA}_{:,:,i}^\top \widehat{\TB}_{:,:,i}.
	\end{align}
By the submultiplicativity of the Frobenius norm, we have
	\begin{align}
	\|\widehat{\TX}_{:,:,i}(\bfM)\|_F  &\le \|(\widehat{\TA}_{:,:,i}^\top \widehat{\TA}_{:,:,i})^\dagger\|_F \|\widehat{\TA}_{:,:,i}\|_F \|\widehat{\TB}_{:,:,i}\|_F
	\end{align}
The assumption of the bounded operator norm on the pseudoinverse requires that $\|\TA^\dagger\| \le C$ for all $\bfM\in \Mcal$. 
By definition, this means the inverse of the smallest nonzero singular value of $\TA$ is bounded above for all transformations. 
Applying this assumption and definition, we have 
	\begin{align}
	 \|(\widehat{\TA}_{:,:,i}^\top \widehat{\TA}_{:,:,i})^\dagger\|_F = \sqrt{\sum_{j=1}^{\min(n_1,n_2)} \widehat{\sigma}_{j,j,i}^{-2}} \le C \sqrt{\min(n_1,n_2)}
	\end{align}
where $\widehat{\sigma}_{j,j,i}$ is the $j$-th singular value of $\widehat{\TA}_{:,:,i}$. 
The same bound holds for each frontal slice in the transform domain.  Thus, we can bound the solution in the transform domain via
	\begin{align}
	\|\widehat{\TX}(\bfM)\|_F \le C \sqrt{n_3\min(n_1,n_2)} \|\widehat{\TA}\|_F \|\widehat{\TB}\|_F
	\end{align}
Because $\bfM$ is orthogonal and the Frobenius norm is orthogonally invariant, the same bound holds with all tensors in the spatial domain. 
\end{itemize}

Combining all of the bounds, we find an upper bound for~\eqref{eq:gradBoundGoal} to be
	\begin{align}
	\|\nabla \overline{\Phi}_{\rm reg}(\bfM)\|_F \le 
	2\|\TB\|_F^2 + C \sqrt{n_3\min(n_1,n_2)} \|\TB\|_F^2 \|\TA\|_F^2.
	\end{align}

Importantly, this bound is independent of $\bfM \in \Mcal$. 
We do not claim this bound is tight; we are concerned with existence and independence of the transformation. 

\end{proof}

A subtle nuance in~\Cref{lem:boundedGradient} is that $\TA$ is not required to have full implicit rank under $\bfM$. 
This was a requirement to obtain a unique solution in~\Cref{thm:tLinearRegUnique}. 
However, uniqueness of the solution is not needed to bound the gradient. 

To speak briefly about the assumption on the $\starM$-operator norm, we expect a sufficient number and variety of data points to form $\TA$ such that the $\starM$-operator norm is bounded for $t$-linear regression problems.  
In practice, we have found the $\starM$-operator norm for $\bfM = \bfI$ is a good estimate. 

%

%% file: sections/S_07_03_bounded_hessian_implication.tex
\subsection{Bounded Approximation Implies Bounded Hessian}
\label{sec:boundedHessianTechnicalDerivation} 

In our setup, we consider the data tensors $\TA \in \Rbb^{n_1\times p\times n_3}$ and $\TB\in \Rbb^{n_1\times n_2 \times n_3}$ with solution $\TX(\bfM)\in \Rbb^{p \times n_2\times n_3}$. 
From~\eqref{eq:gradBoundGoal}, the gradient is
	\begin{subequations}
	\begin{align}
	\bfG_1(\bfM) &= (\TA\starM \TX(\bfM))_{(3)} (\TR(\bfM)_{(3)})^\top\\
	\bfG_2(\bfM) &= (\TR(\bfM) \starM \TX(\bfM)^\top)_{(3)} (\TA_{(3)})^\top\\
	\nabla \overline{\Phi}_{\rm reg}(\bfM) &= \bfM \bfG(\bfM)
	\end{align}
	\end{subequations}
where $\bfG(\bfM) = \bfG_1(\bfM) + \bfG_2(\bfM)$.
We vectorize the Euclidean gradient in two ways
	\begin{subequations}
	\begin{align}
	\myVec(\nabla \overline{\Phi}_{\rm reg}(\bfM) ) 
		&=(\bfG(\bfM)^\top \otimes \bfI_{n_3})\myVec(\bfM)\\
		&=(\bfI_{n_3}\otimes \bfM) \myVec(\bfG(\bfM)).
	\end{align}
	\end{subequations}
We can further vectorize $\bfG(\bfM)$ term-by-term as follows:
	\begin{subequations}
	\begin{align}
	\myVec(\bfG_1(\bfM)) 
		&= (\TR(\bfM)_{(3)} \otimes \bfI_{n_3}) \bfP_{n_1n_2n_3}^\top\myVec(\TA\starM \TX(\bfM))\\
		&=(\bfI_{n_3} \otimes  (\TA\starM \TX(\bfM))_{(3)})  \bfP_{n_1n_2n_3}^\top \myVec(\TR(\bfM))
	\end{align}
	\end{subequations}
and
	\begin{align}
	\myVec(\bfG_2(\bfM))
		&= (\TA_{(3)} \otimes \bfI_{n_3}) \bfP_{n_1pn_3}^\top\myVec(\TR(\bfM) \starM \TX(\bfM)^\top)
	\end{align}
where $\bfP_{n_1pn_3}$ is defined in~\Cref{def:mode3_to_tensor}. 
Differentiating through the gradient, the Hessian of the vectorized format is
	\begin{align}
	\nabla^2 \overline{\Phi}_{\rm reg}(\bfM) = (\bfG(\bfM) \otimes \bfI_{n_3}) + \nabla \myVec(\bfG(\bfM)) (\bfI_{n_3}\otimes \bfM^\top)
	\end{align}
Here, $\nabla \myVec(\bfG(\bfM)) \in \Rbb^{n_3^2 \times n_3^2}$. 
The gradient of $\bfG(\bfM)$ is further broken down into the following:
	\begin{align}
	\begin{split}
	 \nabla \myVec(\bfG_1(\bfM))  
	 	&= \nabla  \myVec(\TA\starM \TX(\bfM)) \bfP_{n_1n_2n_3}  (\TR(\bfM)_{(3)}^\top \otimes \bfI_{n_3})\\
	 	&\qquad +\nabla  \myVec(\TR(\bfM)) \bfP_{n_1n_2n_3} (\bfI_{n_3} \otimes  (\TA\starM \TX(\bfM))_{(3)}^\top)
	\end{split}\\
	 \nabla \myVec(\bfG_2(\bfM))
	 	&=\nabla \myVec(\TR(\bfM) \starM \TX(\bfM)^\top)  \bfP_{n_1pn_3} (\TA_{(3)}^\top \otimes \bfI_{n_3})
	\end{align}
We note that $\nabla  \myVec(\TR(\bfM))  =  \myVec(\TA\starM \TX(\bfM)) $.  
Through the submultiplicativity of the Frobenius norm, the triangle inequality, and the proven boundedness of the gradient in~\Cref{lem:boundedGradient}, we conclude that if the gradient of the approximation, $\nabla \myVec(\TA \starM \TX(\bfM))$, and, relatedly, the gradient of the solution, $\nabla \TX(\bfM)$, are bounded independent of $\bfM$, then the Euclidean Hessian is bounded as well; that is, 
	\begin{align}
	\|\nabla^2 \overline{\Phi}_{\rm reg}(\bfM) \|_F  \le f(n_1,n_2,p,n_3, \|\TA\|_F, \|\TB\|_F).
	\end{align}
where $f$ is a function that depends on the dimensions and norms of the model and observation tensors for the given problem.   
Importantly, this function does not depend on the transformation, $\bfM$. 
We do not claim this bound is tight, we only claim that such a bound exists.  
We prove that $\nabla \widetilde{\TB}(\bfM) =\nabla(\TA \starM \TX(\bfM))$ is bounded in the next section. 

%% file: sections/S_07_04_technical_details.tex
\subsection{Technical Details}
\label{sec:boundedHessianHelperLemmas}

We begin with two definitions to connect vectorized versions of various unfoldings and permutations. 

\begin{definition}[Vectorized Mode-$3$-to-Tensor Permutation Matrix]\label{def:mode3_to_tensor}
Given $\TC\in \Rbb^{m_1\times m_2\times m_3}$, we define the permutation matrix $\bfP_{m_1m_2m_3} \in \Rbb^{m_1m_2m_3 \times m_1m_2m_3}$ that maps a vectorized mode-$3$ unfolding to a vectorized tensor; that is, 
	\begin{align}
	\bfP_{m_1m_2m_3}\myVec(\TC_{(3)}) = \myVec(\TC) 
	= \begin{bmatrix}
	\myVec(\TC_{:,:,1})\\
	\vdots\\
	\myVec(\TC_{:,:,n_3})
	\end{bmatrix}.
	\end{align}	
Because $\bfP_{m_1m_2m_3}$ is orthogonal, $\bfP_{m_1m_2m_3}^\top$ performs the reverse mapping. 
\end{definition}

\begin{definition}[Vectorized Transposition]\label{def:transpose_to_matrix}
Given $\bfZ\in \Rbb^{m_1\times m_2}$, we define the permutation matrix $\bfQ_{m_1m_2} \in \Rbb^{m_1m_2 \times m_1m_2}$ that maps a vectorized transposed matrix to the vectorized non-transposed version; that is, 
	\begin{align}
	\bfQ_{m_1m_2}\myVec(\bfZ^\top) = \myVec(\bfZ).
	\end{align}
Because $\bfQ_{m_1m_2}$ is orthogonal, $\bfQ_{m_1m_2}^\top$ performs the reverse mapping. 
\end{definition}

We present the vectorization of the facewise product~\eqref{eq:facewise} as the following lemma:

\begin{mylemma}{Vectorization of the facewise product}{vecfacewise}
Given $\TA \in \Rbb^{n_1\times p\times n_3}$ and $\TB\in \Rbb^{p\times n_2 \times n_3}$, we have
	\begin{align}
	\begin{split}
		\myVec(\TA \smalltriangleup \TB) 
			&=\myBdiag(\TB_{:,:,i}^\top \otimes \bfI_{n_1}) \myVec(\TA)\\
			&=\myBdiag(\bfI_{n_2} \otimes \TA_{:,:,i}) \myVec(\TB)
	\end{split}
	\end{align}
where for $\TC\in \Rbb^{m_1\times m_2\times m_3}$, we have
	\begin{align}
	\myBdiag(\TC_{:,:,i})  = 
	\begin{bmatrix}
	\TC_{:,:,1}\\
	& \TC_{:,:,2}\\
	&& \ddots \\
	&&& \TC_{:,:,n_3}
	\end{bmatrix}  \in \Rbb^{m_1m_2m_3\times m_1m_2m_3}.
	\end{align} 
\end{mylemma}

\begin{proof}
The proof follows from~\cite{KoldaBader2009} and~\cite{Petersen07thematrix}.  
\end{proof}

We additionally prove an upper bound of the Frobenius norm of the $\starM$-product:
\begin{mylemma}{$\starM$-product upper bound}{starMBounded}
Given $\TA \in \Rbb^{n_1\times p\times n_3}$ and $\TB\in \Rbb^{p\times n_2 \times n_3}$, for any $\bfM\in \Mcal$, we have 
	\begin{align}
	\|\TA \starM \TB\|_F \le \|\TA\|_F \|\TB\|_F.
	\end{align}
\end{mylemma}

\begin{proof}
By the orthogonal invariance of the Frobenius norm, we have $\|\TA \starM \TB\|_F = \|\widehat{\TA} \smalltriangleup \widehat{\TB}\|_F$.  
Expressing the square Frobenius norm in terms of frontal slices, we have
	\begin{subequations}
	\begin{align}
	\|\widehat{\TA} \smalltriangleup \widehat{\TB}\|_F^2 
		&= \sum_{i=1}^{n_3} \|\widehat{\TA}_{:,:,i} \widehat{\TB}_{:,:,i}\|_F^2\\
		&\le \sum_{i=1}^{n_3} \|\widehat{\TA}_{:,:,i}\|_F^2 \|\widehat{\TB}_{:,:,i}\|_F^2\\
		&\le  \left(\sum_{i=1}^{n_3} \|\widehat{\TA}_{:,:,i}\|_F^2 \right) \left(\sum_{i=1}^{n_3} \|\widehat{\TB}_{:,:,i}\|_F^2\right)\\
		&=\|\widehat{\TA}\|_F^2 \|\widehat{\TB}\|_F^2.
	\end{align}
	\end{subequations}
Using orthogonal invariance again completes the proof. 
\end{proof}

For standard operations, we compile a list of upper bounds in~\Cref{tab:standardBounds}.

\begin{table}
\centering
\caption{Standard bounds where $\Mcal = \Ocal_{n_3}$, the set of orthogonal matrices. 
The bounds come from the submultiplicativity and orthogonal invariance of the Frobenius norm}
\label{tab:standardBounds}

\scriptsize
\begin{tabular}{|c|c|c|c|}
\hline
Name & Notation & Bound\\
\hline\hline

\Cref{def:mode3_to_tensor} & $\bfP_{m_1m_2m_3} \in \Rbb^{m_1m_2m_3\times m_1m_2m_3}$ & $\|\bfP_{m_1m_2m_3} \|_F \le \sqrt{m_1m_2m_3}$\\
\hline
\Cref{def:transpose_to_matrix} & $\bfQ_{m_1m_2} \in \Rbb^{m_1m_2\times m_1m_2}$ & $\|\bfQ_{m_1m_2} \|_F \le \sqrt{m_1m_2}$\\
\hline
\Cref{lem:vecfacewise} & $\myBdiag(\TC_{:,:,i} \otimes \bfI_{q}) \in \Rbb^{m_1q m_3\times q m_2 m_3}$ & $\|\myBdiag(\TC_{:,:,i} \otimes \bfI_{q}) \|_F \le \sqrt{q} \|\TC\|_F$\\
	&  $\myBdiag( \bfI_{q} \otimes \TC_{:,:,i}) \in \Rbb^{m_1q m_3\times q m_2 m_3}$ & $\|\myBdiag(\bfI_{q} \otimes \TC_{:,:,i}) \|_F \le \sqrt{q} \|\TC\|_F$\\
\hline
\Cref{lem:starMBounded} & $\TA \starM \TB\in \Rbb^{n_1\times n_2\times n_3}$ & $\|\TA\starM \TB\|_F \le \|\TA\|_F \|\TB\|_F$\\
\hline
transformation & $\bfM\in \Mcal \subset \Rbb^{n_3\times n_3}$ & $\|\bfM\|_F = \sqrt{n_3}$\\
				& $(\bfI_q \otimes \bfM) \in \Rbb^{qn_3 \times qn_3}$ &$ \|(\bfI_q \otimes \bfM) \|_F = \sqrt{qn_3}$\\
				& $(\bfM \otimes \bfI_q) \in \Rbb^{qn_3 \times qn_3}$ & $\|(\bfM \otimes \bfI_q) \|_F = \sqrt{qn_3}$\\
\hline
transform domain & $\widehat{\TC}\in \Rbb^{m_1\times m_1\times m_3}$ & $\|\widehat{\TC}\|_F = \|\TC\|_F$\\
\hline				
\end{tabular}

\end{table}

We now derive the Jacobians and corresponding upper bounds for vectorized versions of key operations used to form the approximation, $\TA\starM \TX(\bfM)$. 

\begin{mylemma}{Jacobians of $\starM$-product}{vecMprod}
Given $\TA \in \Rbb^{n_1\times p\times n_3}$ and $\TB\in \Rbb^{p\times n_2 \times n_3}$, we can vectorize with respect to $\bfM$ as
	\begin{subequations}
	\begin{align}
		\myVec(\TA \starM \TB) 
			&=J_1(\TA,\TB,\bfM) \myVec(\bfM)\\
			&=J_2(\TA,\TB,\bfM) \myVec(\bfM) \\
			&=J_3(\TA,\TB,\bfM)\myVec(\bfM).
	\end{align}
	\end{subequations}
Each $J_i$ is the Jacobian of the vectorized $\starM$-product with respect to the $i$-th argument, defined explicitly as
	\begin{subequations}
	\begin{align}
	\begin{split}
			&J_1(\TA, \TB, \bfM) = \\
			 &\qquad 
			 	\underset{\color{black}(1)}{\underline{\bfP_{n_1n_2n_3}}} 
				\underset{\color{red}(2)}{\red{\underline{(\bfI_{n_1n_2} \otimes \bfM^\top)}}}
				\underset{\color{blue}(3)}{\blue{\underline{\myBdiag(\widehat{\TB}_{:,:,i}^\top \otimes \bfI_{n_1})}}}
				\underset{\color{violet}(4)}{\violet{\underline{\bfP_{n_1pn_3} (\TA_{(3)}^\top \otimes \bfI_{n_3})}}}
		\end{split}\\ \nonumber\\
		\begin{split}
			&J_2(\TA, \TB, \bfM) =\\
			&\qquad 
				\underset{\color{black}(1)}{\underline{\bfP_{n_1n_2n_3}}}
				\underset{\color{red}(2)}{\red{\underline{(\bfI_{n_1n_2} \otimes \bfM^\top)} }}
				\underset{\color{blue}(3)}{\blue{\underline{\myBdiag(\bfI_{n_2} \otimes \widehat{\TA}_{:,:,i})}}}
				\underset{\color{violet}(4)}{\violet{\underline{\bfP_{pn_2n_3} (\TB_{(3)}^\top \otimes \bfI_{n_3})}}}
		\end{split}\\ \nonumber\\
		\begin{split}
			&J_3(\TA, \TB, \bfM) =\\
			&\qquad 
			\underset{\color{black}(1)}{\underline{\bfP_{n_1n_2n_3}}}
			\underset{\color{red}(2)}{\red{\underline{((\widehat{\TA} \smalltriangleup \widehat{\TB})_{(3)}^\top \otimes \bfI_{n_3}) \bfQ_{n_3n_3}}}}
		\end{split}
	\end{align}
	\end{subequations}
where $\bfQ_{m_1m_2} \in \Rbb^{m_1m_2 \times m_1m_2}$ transposes a matrix in vectorized form.  

\end{mylemma}

\begin{proof}
The proof follows from~\cite{KoldaBader2009} and~\cite{Petersen07thematrix}. 
The {\bf black} permutation matrix, labeled as (1), connects the vectorized $\starM$-product to its vectorized mode-$3$ unfolding; that is, 
	\begin{align}
	 \myVec(\TA \starM \TB) = \bfP_{n_1n_2n_3}\myVec(\bfM^\top (\widehat{\TA} \smalltriangleup \widehat{\TB})_{(3)}).
	\end{align}
The \red{\bf red} terms, labeled as (2), come from vectorizing $\bfM^\top (\widehat{\TA} \smalltriangleup \widehat{\TB})_{(3)}$ with respect to each variable. 
The \blue{\bf blue} terms, labeled as (3), come from vectorizing the facewise product $\widehat{\TA} \smalltriangleup \widehat{\TB}$ and using~\Cref{lem:vecfacewise}. 
The \violet{\bf violet} terms, labeled as (4), come from vectorizing the mode-$3$ unfolding of the tensors in the transform domain with respect to $\bfM$ via
	\begin{align}
	\myVec(\widehat{\TC}) = \bfP_{m_1m_2m_3} \myVec(\bfM \TC_{(3)}).
	\end{align}
Following the bounds described in~\Cref{tab:standardBounds} and using the submultiplicativity and orthogonal invariance of the Frobenius norm, we have
	\begin{subequations}
	\begin{align}
	\|J_1(\TA,\TB,\bfM)\|_F &\le n_1^{2}n_2n_3^{2}p^{1/2} \|\TA\|_F\|\TB\|_F \\
	\|J_2(\TA,\TB,\bfM)\|_F  &\le n_1n_2^2n_3^2 p^{1/2} \|\TA\|_F\|\TB\|_F\\
	\|J_3(\TA,\TB,\bfM)\|_F &\le n_1^{1/2}n_2^{1/2}n_3^2\|\TA\|_F \|\TB\|_F,
	\end{align}
	\end{subequations}
which are all bounded above independent of $\bfM$. 

\end{proof}

\begin{mylemma}{Vectorization of $\TA(\bfM)\starM \TB(\bfM)$}{vecSpecial}
Given $\TA(\bfM) \in \Rbb^{n_1\times p\times n_3}$ and $\TB(\bfM)\in \Rbb^{p\times n_2 \times n_3}$, we can vectorize with respect to $\bfM$ as
	\begin{subequations}
	\begin{align}
	\myVec(\TA(\bfM) \starM \TB(\bfM)) &=J_4(\TB(\bfM),\bfM) \myVec(\TA(\bfM))\\
							      &=J_5(\TA(\bfM),\bfM) \myVec(\TB(\bfM))
	\end{align}
	\end{subequations}
Each $J_i$ comes from vectorizing with respect to each tensor and is defined as 
	\begin{subequations}
	\begin{align}
	\begin{split}
	&J_4(\TB(\bfM),\bfM) = \\
	&\qquad \bfP_{n_1n_2n_3} (\bfI_{n_1n_2} \otimes \bfM^\top) \myBdiag(\widehat{\TB}_{:,:,i}(\bfM)^\top \otimes \bfI_{n_1}) \bfP_{n_1pn_3}
	 \underset{\color{magenta}(5)}{\magenta{\underline{(\bfI_{n_1p}\otimes \bfM) \bfP_{n_1pn_3}}}}
	\end{split}\\\nonumber\\
	\begin{split}
	&J_5(\TA(\bfM),\bfM) =\\
	&\qquad \bfP_{n_1n_2n_3} (\bfI_{n_1n_2} \otimes \bfM^\top) \myBdiag(\bfI_{n_2} \otimes \widehat{\TA}_{:,:,i}(\bfM)) \bfP_{pn_2n_3}
	\underset{\color{magenta}(5)}{\magenta{\underline{(\bfI_{pn_2}\otimes \bfM) \bfP_{pn_2n_3}}}}
	\end{split}
	\end{align}
	\end{subequations}
\end{mylemma}

\begin{proof}
The {\bf black} parts of the terms come from~\Cref{lem:vecMprod}. 
The \magenta{\bf magenta} terms, labeled as (5), come from vectorizing $\bfM \TC(\bfM)_{(3)}$ with respect to the tensor, $\TC(\bfM)$.  
The bounds follow from~\Cref{lem:vecMprod} and we get
	\begin{subequations}
	\begin{align}
	 \|J_4(\TB(\bfM),\bfM)\|_F &\le n_1^{3}  n_2n_3^{5/2} p^{3/2}\|\TB(\bfM)\|_F\\
	 \|J_5(\TA(\bfM),\bfM)\|_F &\le n_1  n_2^3 n_3^{5/2} p^{3/2}\ \|\TA(\bfM)\|_F.
	\end{align}
	\end{subequations}
Note that these bounds \underline{do depend on $\bfM$}. 
\end{proof}

\begin{mylemma}{Gradient of $\TX(\bfM)$}{gradX}
Given $\TX(\bfM) \in \Rbb^{p\times n_2\times n_3}$ is the solution to $t$-linear regression, the gradient is 
{
	\begin{align}
	\begin{split}
	&\nabla \myVec(\TX(\bfM)) =\\
	&\qquad J_5(\TA^\top \starM \TA,\bfM)^\dagger[\nabla (\TA^\top \starM \TB)\\
	&\qquad \qquad - (J_1(\TA^\top \starM \TA, \TX(\bfM),\bfM) + J_4(\TX(\bfM),\bfM))]
	\end{split}
	\end{align}
}
\end{mylemma}

\begin{proof}
We derive this through implicit differentiation. 
We start with the (vectorized) normal equations
	\begin{align}
	\nabla \myVec((\TA^\top \starM \TA) \starM \TX(\bfM)) = \nabla \myVec(\TA^\top \starM \TB)
	\end{align}
The left-hand side (LHS) and right-hand side (RHS) can be broken down as
	\begin{align}
	\text{LHS:} &\quad J_1(\TA^\top \starM \TA, \TX(\bfM),\bfM) + J_4(\TX(\bfM),\bfM) + J_5(\TA^\top \starM \TA,\bfM) \nabla \myVec(\TX(\bfM))\\
	 \text{RHS:} &\quad J_1(\TA^\top,\TB,\bfM) + J_2(\TA^\top,\TB,\bfM) + J_3(\TA^\top,\TB,\bfM)
	\end{align}
Solving for $\nabla \myVec{\TX(\bfM)}$, we get
	\begin{align}
	\begin{split}
	&\nabla \myVec{\TX(\bfM)} = \\
	&\qquad J_5(\TA^\top \starM \TA,\bfM)^\dagger[\text{RHS}  - (J_1(\TA^\top \starM \TA, \TX(\bfM),\bfM) + J_4(\TX(\bfM),\bfM) )]
	\end{split}
	\end{align}
where 
	\begin{align}
	\begin{split}
	&J_5(\TA^\top \starM \TA,\bfM)^\dagger  =\\
		&\qquad \bfP_{ppn_3}^\top (\bfI_{pp} \otimes \bfM) \bfP_{ppn_3}\myBdiag(\bfI_{p} \otimes 
		\underset{\orange{(6)}}{\orange{\underline{(\widehat{\TA}_{:,:,i}^\top \widehat{\TA}_{:,:,i})^\dagger})}} (\bfI_{pp}\otimes \bfM^\top) \bfP_{ppn_3}^\top
	\end{split}
	\end{align} 
In \Cref{lem:boundedGradient}, we assume that $\TA^\dagger$ has a bounded operator norm; that is, $\|\TA^\dagger\| \le C$ for all $\bfM\in \Mcal$. 
As a result,  the {\bf \orange{orange}} term, labeled as (6), is also bounded above:
	\begin{align}\label{eq:jacAdjointBound}
	\|J_5(\TA^\top \starM \TA,\bfM)^\dagger\|_F &\le p^{11/2}n_3^{5/2} \|{(\TA^\top  \starM \TA)^\dagger}\|_F \le p^{11/2}n_3^{5/2} {\color{teal} \underline{C \sqrt{n_3 \min(n_1,p)}}}
	\end{align}
The underlined {\color{teal} \bf teal} term is derived in the proof of~\Cref{lem:boundedGradient}. 
Using the submultiplicativity and triangle inequality of the Frobenius norm, the gradient of the solution is bounded above by
	{\scriptsize 
	\begin{subequations}
	\begin{align}
	\begin{split}
	&\|\myVec(\TX(\bfM))\|_F \\
	&\qquad \le \|J_5(\TA^\top \starM \TA,\bfM)^\dagger\|_F\|\text{RHS} -  (J_1(\TA^\top \starM \TA, \TX(\bfM),\bfM) + J_4(\TX(\bfM),\bfM) ))\|_F
	\end{split}\\
		&\qquad \le  \|J_5(\TA^\top \starM \TA,\bfM)^\dagger\|_F\left(\|\text{RHS}\|_F + \|J_1(\TA^\top \starM \TA, \TX(\bfM),\bfM)\|_F + \|J_4(\TX(\bfM),\bfM)\|_F\right)
	\end{align}
	\end{subequations}
	}
Following the bounds in~\Cref{lem:vecMprod},~\Cref{lem:vecSpecial}, and~\eqref{eq:jacAdjointBound} as well as the bounds in~\Cref{lem:starMBounded} and~\Cref{tab:standardBounds}, we have shown that the gradient of the solution is bounded independent of $\bfM$. 
\end{proof}

%% file: sections/04_01_geometric_intuition.tex
\section{Geometric Intuition of $\starM$-Optimization}
\label{sec:leastSquares}

We use the setup from~\Cref{exam:varproAD} for $n_3 = 2$ and add noise in the spatial domain via 
	\begin{subequations}
	\begin{align}
	\TA = \widehat{\TA} \times_3 \bfM_{\rm true}^\top + \eta \TN_{\TA} \qquad \text{and} \qquad 
	\TB = \widehat{\TB} \times_3 \bfM_{\rm true}^\top + \eta \TN_{\TB}
	\end{align}
	\end{subequations}
where $\eta \ge 0$ is the noise level and $\TN_{\bfC}$ is a random tensor the same size of $\TC$ where each entry is drawn from the standard normal distribution.     
An equivalent matrix-vector regression problem minimizes $\|\bfA \bfx - \bfb\|_2$ with model matrix and observations
	\begin{align}
	\bfA = \begin{bmatrix}
	\TA_{:,:,1}\\
	\TA_{:,:,2}
	\end{bmatrix}\in \Rbb^{2n_1\times 2} \qquad \text{and} \qquad 
	\bfb = \begin{bmatrix}
	\TB_{:,1,1}\\
	\TB_{:,1,2}
	\end{bmatrix} \in \Rbb^{2n_1 \times 1}.
	\end{align}
	
\begin{figure}
\centering

\subfloat[Convergence of $\starM$-optimization for various noise levels.  
For each noise level, we initialize $\bfM_0$ with the same random orthogonal matrix. 
(Left): convergence of the reduced $t$-linear objective function value.  
For each noise level $\eta$, $\starM$-optimization converges to the noise level.
(Right): convergence of the solution $\bfM_j$. 
We compare to all equivalent permutations and negations of $\bfM_{\rm true}$, denoted by the set $\Ccal$, and report the smallest error. \label{fig:toyRegressionConvergenceIllustration}]{\begin{tikzpicture}
	
	\node at (0,0) (n0) {\includegraphics[width=0.45\linewidth]{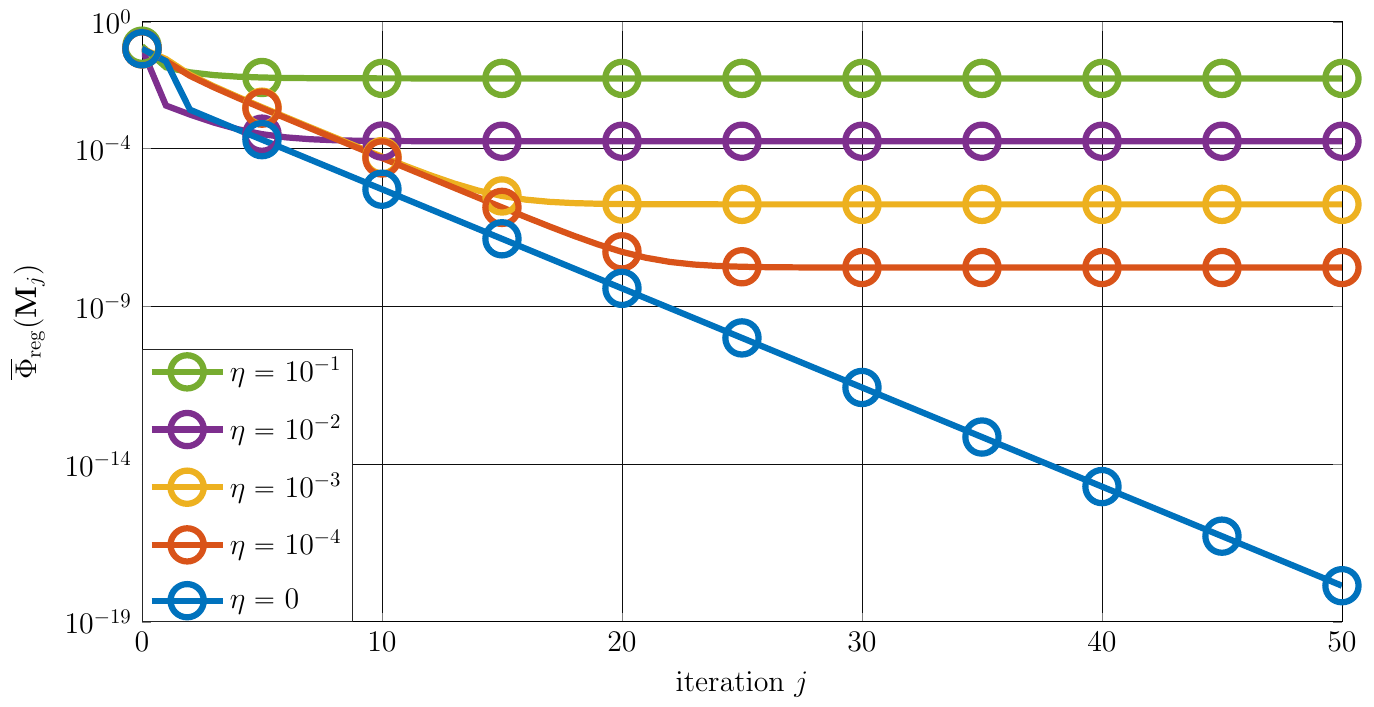}};

	\node[right=0.0cm of n0.east, anchor=west] (n1) {\includegraphics[width=0.45\linewidth]{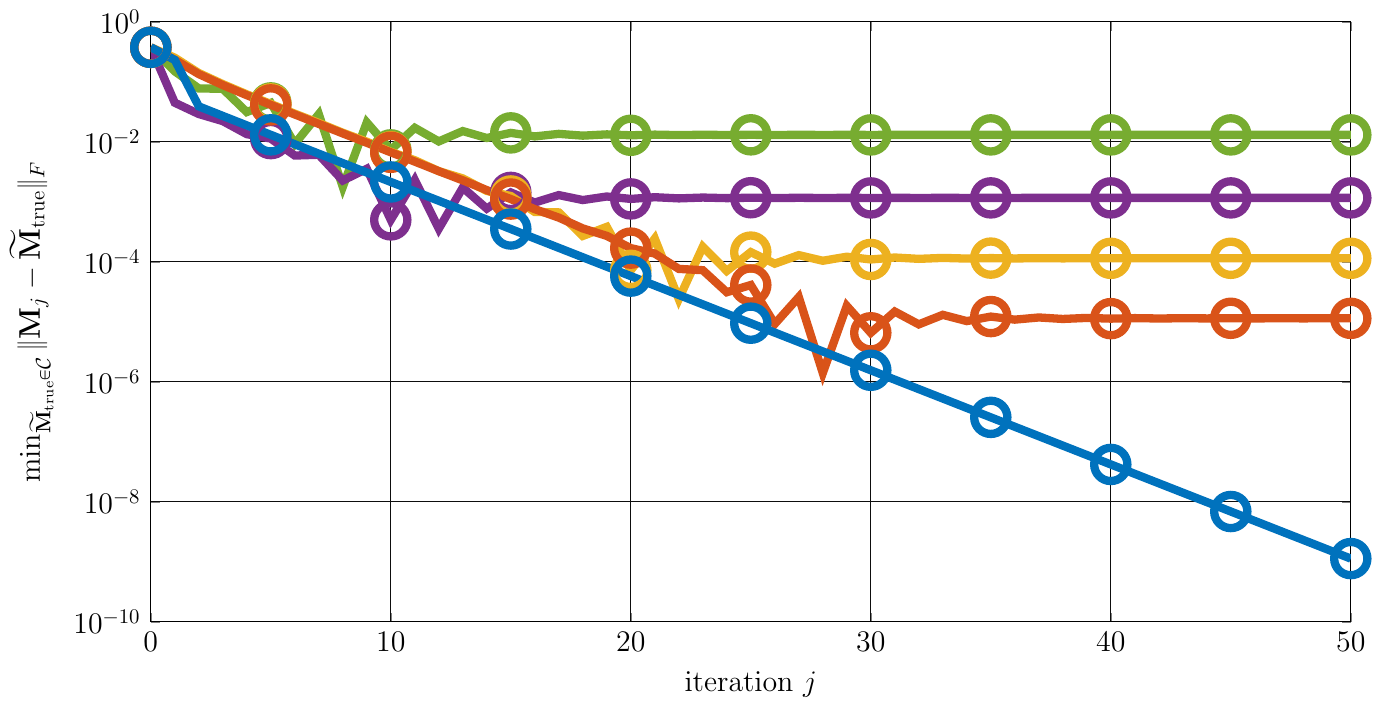}};
	
	\node at ($(n0)!0.5!(n1)$) (c) {};
\end{tikzpicture}
}

\subfloat[Illustration of regression data in the transform domain. 
The scatter points correspond to $(\widehat{\TA}_{:,1,i}, \widehat{\TA}_{:,2,i}, \widehat{\TB}_{:,:,i})$ where \red{\bf red} points come from the first frontal slice ($i = 1$) and \blue{\bf blue} points from the second frontal slice ($i=2$). 
The hyperplanes are the learned models with colors corresponding to the two frontal slices. 
(Left): Solutions for various models and noise levels in the transform domain. 
 Top-to-bottom are different choices of transformation $\bfM$ and the equivalent matrix setup.  
Left-to-right are the results for various noise levels.  
The relative error $\|\TA \starM \TX_{\rm reg}(\bfM) - \TB\|_F / \|\TB\|_F$ is reported below each image.
(Right): Geometric intuition of $\starM$-optimization for $\eta=10^{-2}$.  
As we iterate, the data points become linearly correlated in the transform domain. \label{fig:toyRegressionHyperplanellustration}]{
\begin{tikzpicture}
\node[draw=yellow, line width=4pt, inner sep=0.0cm] at (0,0) (n0) {\includegraphics[width=0.44\linewidth]{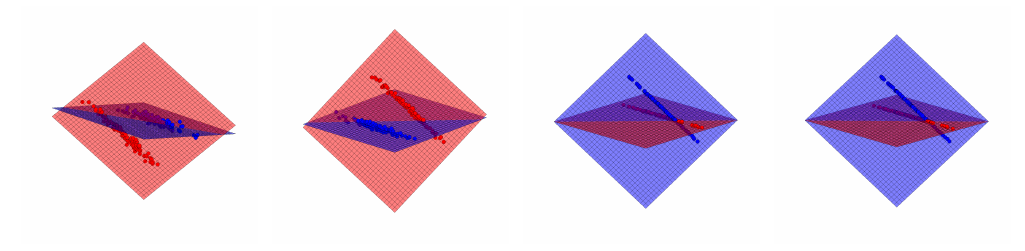}};

\draw[->, line width=4pt, shorten >=0.5cm, shorten <=0.5cm] ($(n0.south west)+(0,-0.3)$) -- node[midway, below] {iteration} ($(n0.south east)+(0,-0.3)$) ;
\node[above=0.0cm of n0.north, anchor=south] {First four iterations for $\eta=10^{-2}$};

\node[left=0.0cm of n0.west, anchor=east] {\includegraphics[width=0.5\linewidth]{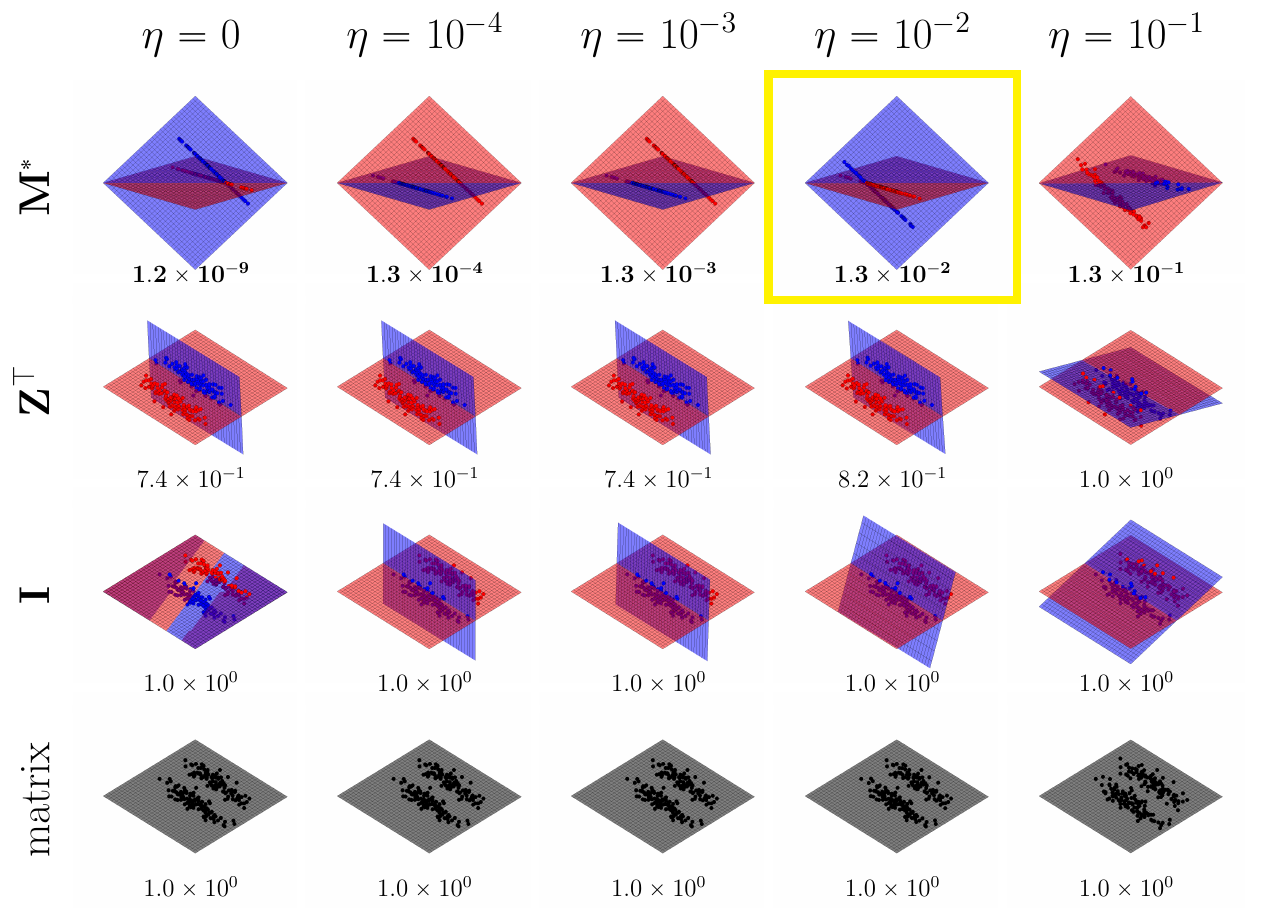}};
\end{tikzpicture}
}

\caption{Convergence of $\starM$-optimization for $t$-linear regression for various noise levels.}
\label{fig:regressionConvergence}
\end{figure}
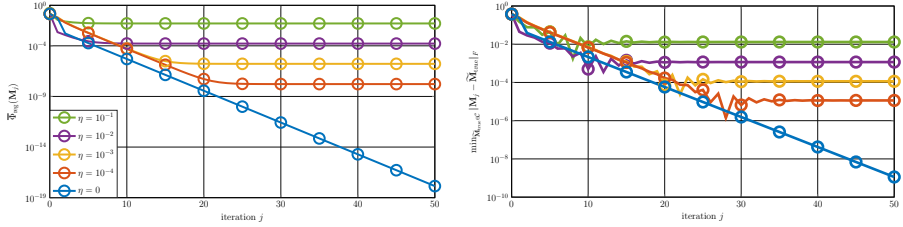
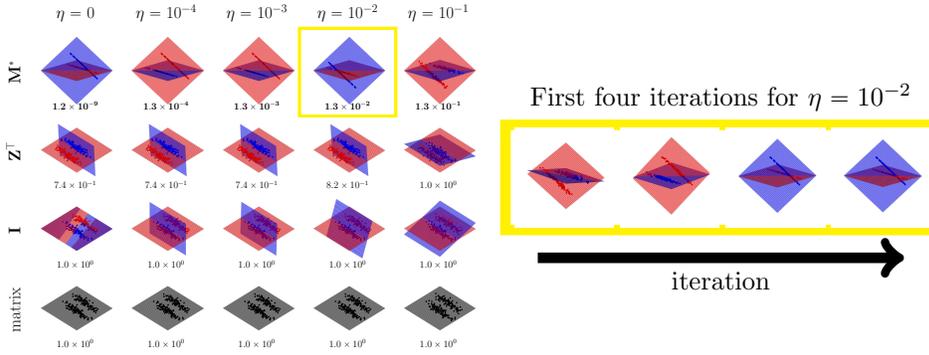

 We present results in~\Cref{fig:regressionConvergence} and discuss the main takeaways.

\paragraph{Effects of Noise} 
In~\Cref{fig:toyRegressionConvergenceIllustration}, we report the convergence of the objective function and the relative error to the nearest true transformation (\Cref{sec:uniqueInvariancePrototype}) for various noise levels.  
As expected, the magnitude of the objective function value and error are comparable to the level of noise. 
For example, the objective function converges to approximately the noise level (squared); e.g., if $\eta = 10^{-3}$, then $\overline{\Phi}(\bfM^*) = \Ocal(10^{-6})$.   
Thus, we fit the data as well as possible and do not overfit the noise because our $t$-linear system is sufficiently overdetermined. 

\paragraph{Comparison to Heuristics} In~\Cref{fig:toyRegressionHyperplanellustration}, we observe that the matrix case does not capture the data with a linear model, as we would expect by construction of the problem. 
For the $t$-linear regression approach, the identity $\bfI$ and data-dependent $\bfZ^\top$ matrices fail to adequately linearly correlate the data and the resulting transformed data to not align well with the learned hyperplanes. 
In comparison, the learned transformation $\bfM^*$ moves data from the spatial domain to the transform domain such that the resulting data is strongly linearly correlated. 
The corresponding hyperplanes thereby capture the behavior of the data well, and we obtain a quality model. 
As we run $\starM$-optimization, we see the data in the transform domain become more linearly correlated, resulting in better hyperplane approximations.

%% file: sections/S_05_varpro_unit_simplex.tex
\section{Variable Projection with Unit Simplex Constraints} 
\label{app:varproConstrained}

In the vector case, suppose we seek a minimizer of the following problem
	\begin{align}\label{eq:varproSimplex}
	\bfx^*(p) = \argmin_{\bfx \in \Delta^n} f(\bfx, p) 
	\end{align}
where $\Delta^n$ is the unit simplex and $p\in \Rbb$ is a scalar parameter for simplicity.  
If $\bfx\in \Delta^n$, then $\bfx$ is a vector of discrete probabilities; that is, $\bfx \ge \bf0$ entrywise and $\bfe^\top \bfx = 1$ where $\bfe$ is the vector of all ones.    
We assume $f$ is convex in $\bfx$, and beause $\Delta^n$ forms a convex set, the resulting optimization problem is convex and admits a global minimizer. 
To solve~\eqref{eq:varproSimplex}, we set up the Lagrangian, following the conventions in~\cite{Beck2014}
	\begin{align}
	\Lcal(\bfx, \bflambda, \mu) \equiv f(\bfx,p)- \bflambda^\top \bfx + \mu(\bfe^\top \bfx - 1)
	\end{align}
where $\bflambda \in \Rbb_+^n$ and $\mu\in \Rbb$ are the dual variables. 
The corresponding Karush-Kuhn-Tucker (KKT) conditions are the following:
	\begin{subequations}
	\begin{align}
	 \nabla_{\bfx} f(\bfx,p) - \bflambda + \mu \bfe &= {\bf0} && \text{stationarity}  \label{eq:kktStationary}\\
	\bfe^\top \bfx &= 1 &&\text{primal feasibility} \label{eq:kktPrimal}\\
	\bflambda &\ge {\bf0} && \text{dual feasibility}  \label{eq:kktDual}\\
	\lambda_i x_i &= 0 \qquad \text{for $i=1,\dots, n$} && \text{complementary slackness}  \label{eq:kktSlack}
	\end{align}
	\end{subequations} 
where stationarity can be equivalently written as $\nabla_{\bfx} \Lcal(\bfx, \bflambda, \mu)$. 
Let $\bfx^*(p)$ be a KKT point.  
We are interested in computing the gradient at a KKT optimal point with respect to the parameters; that is,
	\begin{align}\label{eq:varproSimplexGradient}
	\nabla_p f(\bfx^*(p),p) = \nabla_p \bfx^*(p) \nabla_{\bfx} f(\bfx^*(p),p) + \nabla_p f(\bfx^*(p),p)
	\end{align}
The challenging term to compute is $\nabla_p \bfx^*(p)$, the sensitivity of $\bfx^*(p)$ with respect to the parameters.  
We explore this in two cases.  

\begin{figure}

\centering

\begin{tabular}{cc}
\includegraphics[width=0.4\linewidth]{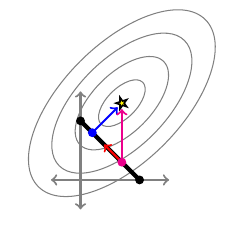}
&
\includegraphics[width=0.4\linewidth]{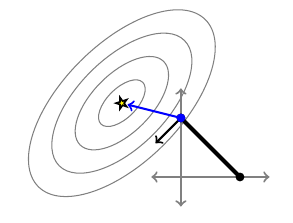}
\end{tabular}
\caption{Illustration of optimizing over the unit simplex. 
The star is the unconstrained minimum of the objective function given by the gray contour lines. 
The black line is the unit simplex, the \blue{\bf blue} and  \magenta{\bf magenta} points are current iterates that satisfy the constraint, and the blue point is optimal. 
The vectors of correspond color are the (negative) gradient directions for the unconstrained problem. 
The \red{\bf red} vector in the left plot shows the direction the magenta point will move when updating with a projected gradient step. 
The black vector in right plot is orthogonal to the simplex, demonstrating that blue point cannot move any closer to the minimum and thus is stationary. 
 }
\label{fig:nonnegativeVarPro}
\end{figure}

\paragraph{Case 1: No Active Inequality Constraints ($\bfx^*(p) > 0$)} 
From complementary slackness, we have $\bflambda^* = \bf0$. 
Substituting $\bflambda^* = \bf0$ into the stationarity condition, we have $\nabla_{\bfx} f(\bfx^*(p),p) = -\mu^* \bfe$. 
As expected, the gradient of the objective function at the stationary point is normal to the equality constraint (\Cref{fig:nonnegativeVarPro}).
The change in the optimal point, $\nabla_p \bfx^*(p)^\top \in \Rbb^n$, must ensure $\bfx^*(p)$ to remains feasible; i.e., 
   \begin{subequations}
    \begin{align}
    \bfe^\top (\bfx^*(p) + \alpha \nabla_p \bfx^*(p)^\top) &= 1\\
    \bfx^*(p) + \alpha \nabla_p \bfx^*(p)^\top &\ge \bf0
    \end{align}
    \end{subequations}
where $\alpha \ge 0$ is a step size.
Because $\bfx^*(p)$ is already feasible, the first equality guarantees that that $\nabla_p \bfx^*(p)^\top$ lies in the null space of $\bfe$. 
Hence, $\nabla_p \bfx^*(p)^\top$ is orthogonal to $\nabla_{\bfx} f(\bfx^*(p),p)$ and the first term in~\eqref{eq:varproSimplexGradient} vanishes. 
Thus, we can ignore the sensitivity of the solution with respect to the parameters when the constraints are inactive. 

\paragraph{Case 2: Active Constraints ($x_i = 0 \Longrightarrow \lambda_i \ge 0$ on $\Acal$)}
Let $\Acal$ be the active set (i.e., the set of indices for which $x_i = 0$) and let $\Acal^C$ be the inactive set (i.e., the set of indices for which $x_i > 0$).    
In this case, $\nabla_{\bfx} f(\bfx^*(p),p) = \bflambda -\mu \bfe$.  
Let $\bfS \in \Rbb^{n\times |\Acal^C|}$ be a subsampling matrix that selects entries in $\Acal^C$. 
Then, we have $\bfS^\top \nabla_{\bfx} f(\bfx^*(p),p) = -\mu \bfS^\top\bfe$.  
Unsurprisingly, on the inactive set, the gradient behaves the same as in Case 1.   

On active constraints, we have $x_i^*(p) = 0$ and the gradient $\nabla_{\bfx} f(\bfx^*(p),p)$ points in the steepest ascent direction, which is no longer guaranteed to be normal to the unit simplex (\Cref{fig:nonnegativeVarPro}). 
This case does not guarantee the first term of~\eqref{eq:varproSimplexGradient} vanishes.  
In practice, we optimize with an interior point method from {\sc Matlab} to avoid hitting the constraints and thus only require Case 1 to hold.  


%% file: sections/S_06_stocks.tex
\section{Portfolio Stocks}
\label{app:stocks}

We provide a list of stocks used for the index tracking experiments in~\Cref{tab:stocks}.  
The tickers can be translated to specific companies through, e.g., Yahoo Finance~\cite{yahoo}. 
The companies are roughly ordered by market cap and we chose among the largest possible companies per sector. 

\begin{table}
\centering
\scriptsize
\caption{Stocks used for the tensor index tracking experiment in~\Cref{sec:numerical}. 
The colors of the sectors are consistent throughout the manuscript.}
\label{tab:stocks}

\begin{tabular}{lcl}
 \cellcolor{jet1} \color{white} 1. Communication Services &: &  \tt CSCO,  TMUS,  VZ, CMCSA, AMX, ORAN, DIS, T, DASH, ZM \\
 \cellcolor{jet2}  \color{white}2.  Consumer Discretionary &: & \tt TGT,  AMZN,  WMT, RCL, HD, LVMHF, TM, MCD,  NKE, SBUX\\
 \cellcolor{jet3}  \color{white} 3. Consumer Staples & : &\tt PG, KO,  PEP, NSRGY, LRLCY, COST, PM, UL, BUD, EL \\
 \cellcolor{jet4}  4. Energy &: & \tt XOM,  CVX,  COP, SHEL, TTE, SLB, BP, EQNR, PBR, EOG\\
 \cellcolor{jet5}  5. Financials &: & \tt JPM,  BAC,  WFC, HSBC, HDB, MS, SCHW, TD, GS, C\\
 \cellcolor{jet6}  6. Healthcare &: & \tt UNH,  JNJ,  LLY,  NVO,  MRK,  RHHBY,  PFE,  TMO,  ABT,  DHR\\
 \cellcolor{jet7}  7. Industrials &: &\tt UPS,  BA,  ACN, CAT, RTX, HON, UNP, GE, DE, ADP\\
 \cellcolor{jet8}  8. Information Technology  &: & \tt AAPL,  MSFT,  GOOGL, IBM, CRM, CSCO, FSLR, ACN, ENPH,  AVGO\\
 \cellcolor{jet9}  9. Materials &: & \tt LIN,  FCX,  SHW, BHP, RIO, APD, SCCO, ECL, GOLD, VALE\\
 \cellcolor{jet10}  10. Real Estate & : &\tt SPG,  PSA,  O, PLD, AMT, EQIX, CCI, WELL, CSGP, DLR\\
 \cellcolor{jet11}  11. Utilities &: & \tt NEE,  D,  ED, SO, DUK, NGG, SRE, D, AEP, XEL \\
\end{tabular}
\end{table}

%% file: sections/S_10_wave_equation.tex
\section{Wave Equation Details}
\label{sec:wave_equation}

\begin{figure}
\centering

\subfloat[Finite element mesh and boundary conditions~\eqref{eq:BC1} and~\eqref{eq:BC2}.]{
\begin{tikzpicture}

\def\b{0.25}
\def\n{4.0}
\def\w{4}
\def\s{0.125}
\node[anchor=south west] at (-\s,-\s) {\includegraphics[width=\w cm]{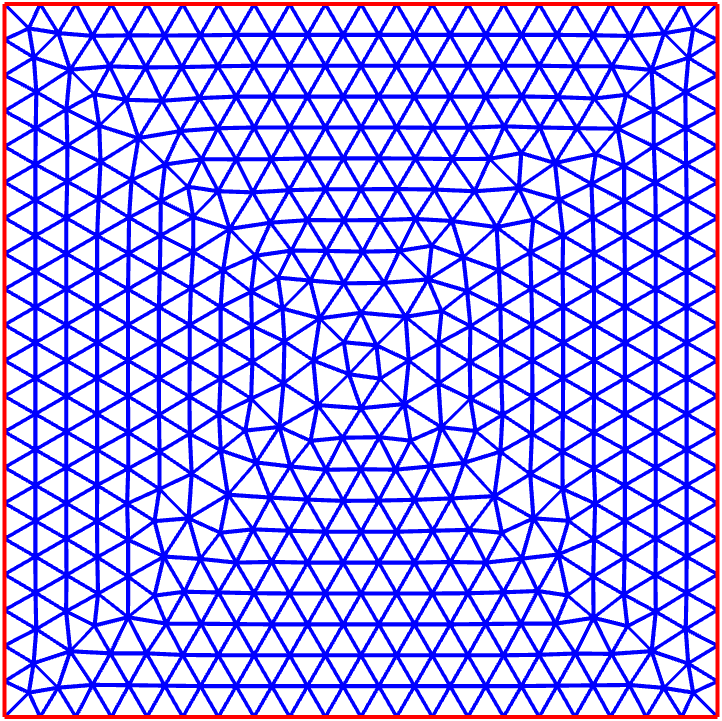}};
\draw[line width=2pt] (0,0) rectangle (\n,\n);
\node at (-\b,-\b) {$-1$};
\node at (\n,-\b) {$1$};
\node at (-\b,\n) {$1$};

\node[rotate=270] at (\n+\b,\n/2) {$x = 0$};
\node[rotate=0] at (\n/2,\n+\b) {$\nabla x = 0$};
\node[rotate=0] at (\n/2,-\b) {$ \nabla x = 0$};
\node[rotate=90] at (-\b, \n/2) {$x = 0$};
\end{tikzpicture}
}
\hspace{1cm}
\subfloat[3D visualization of true solution $\bfx(z_1,z_2,0)$ at initial time point.]{
\includegraphics[width=5cm]{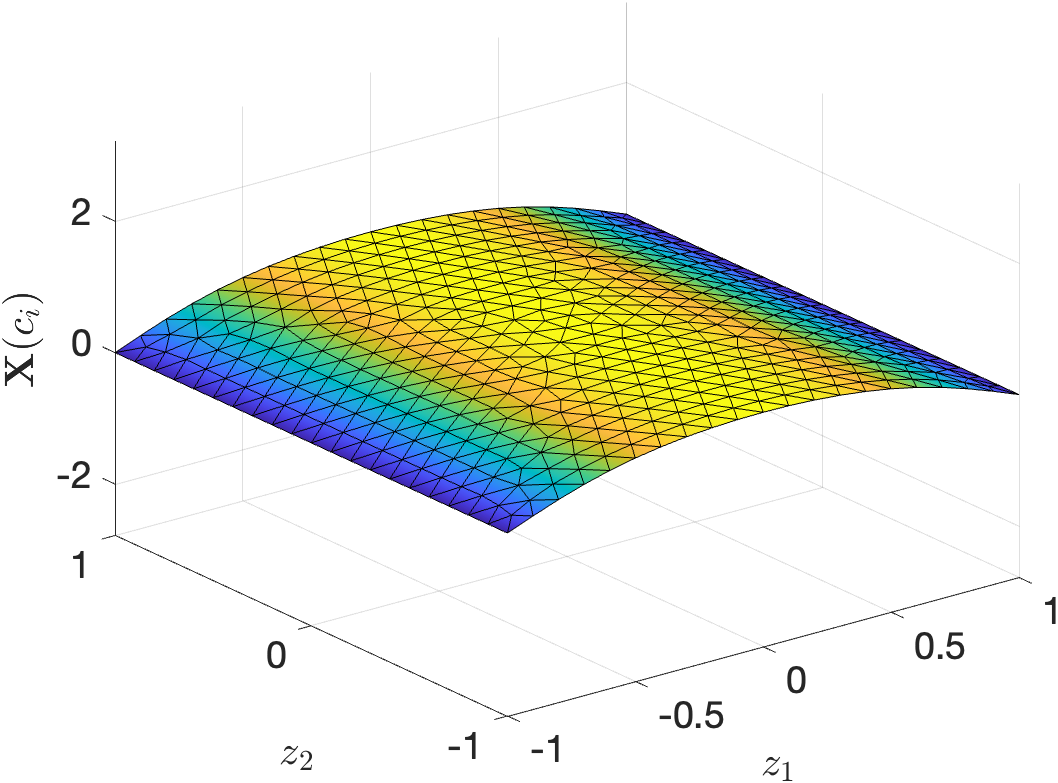}
}

\subfloat[First initial condition~\eqref{eq:IC1}.]{
\begin{tikzpicture}
\begin{axis}[scale only axis,
	width=4cm, 
	height=4cm,
        xmin=-1, xmax=1, 
        ymin=-1, ymax=1, 
        domain=-1:1,  
        xlabel={$z_1$},
        ylabel={$x(z_1,\cdot,0)$},
        xtick={-1,0,1},
        ytick={-1,0,1},
]
\addplot[mark=none, line width=3pt, color=blue, samples=500]    {rad(atan(cos(deg(pi * x / 2))))};
\addplot[mark=none, dashed, line width=1pt, color=gray, samples=500]    {0 * x};
\end{axis}
\end{tikzpicture}
}
\hspace{1cm}
\subfloat[Second initial condition~\eqref{eq:IC1}.]{
\begin{tikzpicture}
\begin{axis}[scale only axis,
	width=4cm, 
	height=4cm,
        xmin=-1, xmax=1, 
        ymin=-2, ymax=5, 
        domain=-1:1,  
        xlabel={$z_1$},
        ylabel={$\dot{x}(z_1,\cdot,0)$},
                xtick={-1,0,1},
        ytick={-2,0,2,5},
]
\addplot[mark=none, line width=3pt, color=red, samples=500]    {2 * sin(deg(pi * x)) * exp(sin(deg(pi * x / 2)))};
\addplot[mark=none, dashed, line width=1pt, color=gray, samples=500]    {0 * x};

\end{axis}
\end{tikzpicture}
}

%

\caption{Illustration of the initial solution to the wave equation, the finite element discretization, and the initial conditions.}
\label{fig:wave_equation}

\end{figure}

For~\Cref{sec:rom}, the  homogeneous two-dimensional wave equation from the {\sc Matlab} PDE Toolbox\footnote{\url{https://www.mathworks.com/help/pde/ug/wave-equation.html}} is given by
	\begin{subequations}\label{eq:wave_pde}
	\begin{align}
	\ddot{x}(z_1,z_2,t)+ c \Delta x(z_1,z_2,t) &= 0 && \text{on $\Omega \times (0,5]$}\label{eq:pde_interior}\\
							x(z_1,z_2,t) &= 0 &&z_1\in \{-1,1\} \label{eq:BC1}\\
							\nabla x(z_1,z_2,t) &= 0 && z_2 \in \{-1,1\} \label{eq:BC2} \\
							x(z_1,z_2,0) &=\arctan(\cos(\tfrac{\pi}{2} z_1)) \label{eq:IC1}\\
							\dot{x}(z_1,z_2,0) &=2 \sin(\pi z_1) e^{\sin(\tfrac{\pi}{2}  z_2)}\label{eq:IC2}
	\end{align}
	\end{subequations}
where $\Omega = (-1,1) \times (-1,1)$ and $c$ is the unknown wave speed.   
We plot the finite element mesh and the initial conditions in~\Cref{fig:wave_equation}.

%% file: Newman_OptimalTensorAlgebras_ArXiv.bbl
\begin{thebibliography}{10}

\bibitem{Absil2008}
{\sc P.-A. Absil, R.~Mahony, and R.~Sepulchre}, {\em Optimization Algorithms on
  Matrix Manifolds}, Princeton University Press, Princeton, NJ, 2008.

\bibitem{ANTIL20121986}
{\sc H.~Antil, M.~Heinkenschloss, R.~H. Hoppe, C.~Linsenmann, and A.~Wixforth},
  {\em Reduced order modeling based shape optimization of surface acoustic wave
  driven microfluidic biochips}, Mathematics and Computers in Simulation, 82
  (2012), pp.~1986--2003,
  \url{https://doi.org/https://doi.org/10.1016/j.matcom.2010.10.027},
  \url{https://www.sciencedirect.com/science/article/pii/S0378475410003496}.
\newblock "The Fourth IMACS Conference : Mathematical Modelling and
  Computational Methods in Applied Sciences and Engineering" Devoted to Owe
  Axelsson in ocassion of his 75th birthday.

\bibitem{Beck2014}
{\sc A.~Beck}, {\em Introduction to Nonlinear Optimization}, Society for
  Industrial and Applied Mathematics, Philadelphia, PA, 2014,
  \url{https://doi.org/10.1137/1.9781611973655},
  \url{https://epubs.siam.org/doi/abs/10.1137/1.9781611973655},
  \url{https://arxiv.org/abs/https://epubs.siam.org/doi/pdf/10.1137/1.9781611973655}.

\bibitem{Benidis2018:indexTracking}
{\sc K.~Benidis, Y.~Feng, and D.~P. Palomar}, {\em Sparse portfolios for
  high-dimensional financial index tracking}, IEEE Transactions on Signal
  Processing, 66 (2018), pp.~155--170,
  \url{https://doi.org/10.1109/TSP.2017.2762286}.

\bibitem{BrennerGugercinWillcox2015:MOR}
{\sc P.~Benner, S.~Gugercin, and K.~Willcox}, {\em A survey of projection-based
  model reduction methods for parametric dynamical systems}, SIAM Review, 57
  (2015), pp.~483--531, \url{https://doi.org/10.1137/130932715},
  \url{https://doi.org/10.1137/130932715},
  \url{https://arxiv.org/abs/https://doi.org/10.1137/130932715}.

\bibitem{boumal2023intromanifolds}
{\sc N.~Boumal}, {\em An introduction to optimization on smooth manifolds},
  Cambridge University Press, 2023,
  \url{https://doi.org/10.1017/9781009166164},
  \url{https://www.nicolasboumal.net/book}.

\bibitem{BoumalAbsilCoralia2018:globalManifold}
{\sc N.~Boumal, P.-A. Absil, and C.~Cartis}, {\em {Global rates of convergence
  for nonconvex optimization on manifolds}}, IMA Journal of Numerical Analysis,
  39 (2018), pp.~1--33, \url{https://doi.org/10.1093/imanum/drx080},
  \url{https://doi.org/10.1093/imanum/drx080},
  \url{https://arxiv.org/abs/https://academic.oup.com/imajna/article-pdf/39/1/1/33283937/drx080.pdf}.

\bibitem{Boyd_Vandenberghe_2018}
{\sc S.~Boyd and L.~Vandenberghe}, {\em Introduction to Applied Linear Algebra:
  Vectors, Matrices, and Least Squares}, Cambridge University Press, 2018.

\bibitem{CarrollChang1970}
{\sc J.~D. Carroll and J.-J. Chang}, {\em Analysis of individual differences in
  multidimensional scaling via an n-way generalization of ``eckart-young''
  decomposition}, Psychometrika, 35 (1970), pp.~283--319.

\bibitem{Cui:2017indianPines}
{\sc B.~Cui, X.~Ma, F.~Zhao, and Y.~Wu}, {\em A novel hyperspectral image
  classification approach based on multiresolution segmentation with a few
  labeled samples}, International Journal of Advanced Robotic Systems, 14
  (2017), p.~1729881417710219, \url{https://doi.org/10.1177/1729881417710219},
  \url{https://doi.org/10.1177/1729881417710219},
  \url{https://arxiv.org/abs/https://doi.org/10.1177/1729881417710219}.

\bibitem{LathauwerMoorVandewalle2000}
{\sc L.~de~Lathauwer, B.~de~Moor, and J.~Vandewalle}, {\em A multilinear
  singular value decomposition}, SIAM Journal on Matrix Analysis and
  Applications, 21 (2000), pp.~1253--1278.

\bibitem{Demmel1997}
{\sc J.~W. Demmel}, {\em Applied Numerical Linear Algebra}, Society for
  Industrial and Applied Mathematics, 1997,
  \url{https://doi.org/10.1137/1.9781611971446},
  \url{https://epubs.siam.org/doi/abs/10.1137/1.9781611971446},
  \url{https://arxiv.org/abs/https://epubs.siam.org/doi/pdf/10.1137/1.9781611971446}.

\bibitem{Eckart1936}
{\sc C.~Eckart and G.~Young}, {\em The approximation of one matrix by another
  of lower rank}, Psychometrika, 1 (1936), pp.~211--218,
  \url{https://doi.org/10.1007/BF02288367},
  \url{https://doi.org/10.1007/BF02288367}.

\bibitem{EdelmanAriasSmith1998}
{\sc A.~Edelman, T.~A. Arias, and S.~T. Smith}, {\em The geometry of algorithms
  with orthogonality constraints}, SIAM Journal on Matrix Analysis and
  Applications, 20 (1998), pp.~303--353,
  \url{https://doi.org/10.1137/S0895479895290954},
  \url{https://doi.org/10.1137/S0895479895290954},
  \url{https://arxiv.org/abs/https://doi.org/10.1137/S0895479895290954}.

\bibitem{FAN20122112}
{\sc Y.-W.~D. Fan and J.~G. Nagy}, {\em An efficient computational approach for
  multiframe blind deconvolution}, Journal of Computational and Applied
  Mathematics, 236 (2012), pp.~2112--2125,
  \url{https://doi.org/https://doi.org/10.1016/j.cam.2011.09.034},
  \url{https://www.sciencedirect.com/science/article/pii/S0377042711005139}.
\newblock Inverse Problems: Computation and Applications.

\bibitem{Giles1948}
{\sc M.~Giles}, {\em An extended collection of matrix derivative results for
  forward and reverse mode algorithmic differentiation matrix product , inverse
  and determinant preliminaries}, University Computing,  (1948).

\bibitem{GleichEtAl2013:algebraOfCirculants}
{\sc D.~F. Gleich, C.~Greif, and J.~M. Varah}, {\em The power and arnoldi
  methods in an algebra of circulants}, Numerical Linear Algebra with
  Applications, 20 (2013), pp.~809--831,
  \url{https://doi.org/https://doi.org/10.1002/nla.1845},
  \url{https://onlinelibrary.wiley.com/doi/abs/10.1002/nla.1845},
  \url{https://arxiv.org/abs/https://onlinelibrary.wiley.com/doi/pdf/10.1002/nla.1845}.

\bibitem{GolubPereyra1973:VarPro}
{\sc G.~H. Golub and V.~Pereyra}, {\em The differentiation of pseudo-inverses
  and nonlinear least squares problems whose variables separate}, SIAM Journal
  on Numerical Analysis, 10 (1973), pp.~413--432,
  \url{https://doi.org/10.1137/0710036}, \url{https://doi.org/10.1137/0710036},
  \url{https://arxiv.org/abs/https://doi.org/10.1137/0710036}.

\bibitem{GoluVanl96}
{\sc G.~H. Golub and C.~F. Van~Loan}, {\em Matrix Computations}, The Johns
  Hopkins University Press, third~ed., 1996.

\bibitem{Harshman1970}
{\sc R.~A. Harshman}, {\em Foundations of the parafac procedure: Models and
  conditions for an "explanatory" multimodal factor analysis}, UCLA Working
  Papers in Phonetics, 16 (1970).

\bibitem{HAY_BORGGAARD_PELLETIER_2009}
{\sc A.~Hat, J.~T. Borggaard, and D.~Pelletier}, {\em Local improvements to
  reduced-order models using sensitivity analysis of the proper orthogonal
  decomposition}, Journal of Fluid Mechanics, 629 (2009), p.~41–72,
  \url{https://doi.org/10.1017/S0022112009006363}.

\bibitem{Higham2008:matrixFunctions}
{\sc N.~J. Higham}, {\em Functions of Matrices}, Society for Industrial and
  Applied Mathematics, 2008, \url{https://doi.org/10.1137/1.9780898717778},
  \url{https://epubs.siam.org/doi/abs/10.1137/1.9780898717778},
  \url{https://arxiv.org/abs/https://epubs.siam.org/doi/pdf/10.1137/1.9780898717778}.

\bibitem{Hitchcock1927}
{\sc F.~L. Hitchcock}, {\em The expression of a tensor or a polyadic as a sum
  of products}, Journal of Mathematics and Physics, 6 (1927), pp.~164--189,
  \url{https://doi.org/https://doi.org/10.1002/sapm192761164},
  \url{https://onlinelibrary.wiley.com/doi/abs/10.1002/sapm192761164},
  \url{https://arxiv.org/abs/https://onlinelibrary.wiley.com/doi/pdf/10.1002/sapm192761164}.

\bibitem{Kaufman1975:VarPro}
{\sc L.~Kaufman}, {\em A variable projection method for solving separable
  nonlinear least squares problems}, BIT Numerical Mathematics, 15 (1975),
  pp.~49--57, \url{https://doi.org/10.1007/BF01932995},
  \url{https://doi.org/10.1007/BF01932995}.

\bibitem{Keegan2022}
{\sc K.~Keegan, T.~Vishwanath, and Y.~Xu}, {\em A tensor svd-based
  classification algorithm applied to fmri data}, SIAM Undergraduate Research
  Online, 15 (2022), pp.~270--294.

\bibitem{KernfeldKilmer2015}
{\sc E.~Kernfeld, M.~Kilmer, and S.~Aeron}, {\em Tensor--tensor products with
  invertible linear transforms}, Linear Algebra and its Applications, 485
  (2015), pp.~545--570,
  \url{https://doi.org/https://doi.org/10.1016/j.laa.2015.07.021},
  \url{https://www.sciencedirect.com/science/article/pii/S0024379515004358}.

\bibitem{KilmerMartinPerrone2008}
{\sc M.~Kilmer, C.~Martin, and L.~Perrone}, {\em A third-order generalization
  of the matrix svd as a product of third-order tensors},  (2008).

\bibitem{KilmerEtAl2013:thirdOrderTensorOperators}
{\sc M.~E. Kilmer, K.~Braman, N.~Hao, and R.~C. Hoover}, {\em Third-order
  tensors as operators on matrices: A theoretical and computational framework
  with applications in imaging}, SIAM Journal on Matrix Analysis and
  Applications, 34 (2013), pp.~148--172,
  \url{https://doi.org/10.1137/110837711},
  \url{https://doi.org/10.1137/110837711},
  \url{https://arxiv.org/abs/https://doi.org/10.1137/110837711}.

\bibitem{Kilmer2021:pnas}
{\sc M.~E. Kilmer, L.~Horesh, H.~Avron, and E.~Newman}, {\em Tensor-tensor
  algebra for optimal representation and compression of multiway data},
  Proceedings of the National Academy of Sciences of the United States of
  America, 118 (2021), \url{https://doi.org/10.1073/pnas.2015851118}.

\bibitem{KILMER2011641}
{\sc M.~E. Kilmer and C.~D. Martin}, {\em Factorization strategies for
  third-order tensors}, Linear Algebra and its Applications, 435 (2011),
  pp.~641--658,
  \url{https://doi.org/https://doi.org/10.1016/j.laa.2010.09.020},
  \url{https://www.sciencedirect.com/science/article/pii/S0024379510004830}.
\newblock Special Issue: Dedication to Pete Stewart on the occasion of his 70th
  birthday.

\bibitem{KoldaBader2009}
{\sc T.~G. Kolda and B.~W. Bader}, {\em Tensor decompositions and
  applications}, SIAM Review, 51 (2009), pp.~455--500.

\bibitem{KongLuLin2021:TensorQRank}
{\sc H.~Kong, C.~Lu, and Z.~Lin}, {\em Tensor q-rank: New data dependent
  definition of tensor rank}, Mach. Learn., 110 (2021), pp.~1867--1900,
  \url{https://doi.org/10.1007/s10994-021-05987-8},
  \url{https://doi.org/10.1007/s10994-021-05987-8}.

\bibitem{Lensky2023:yahooFinance}
{\sc A.~Lensky}, {\em Yahoo finance and quandl data downloader}, 2023,
  \url{https://github.com/Lenskiy/Yahoo-Quandl-Market-Data-Donwloader/releases/tag/v1.131}.

\bibitem{LieuFarhat2007}
{\sc T.~Lieu and C.~Farhat}, {\em Adaptation of aeroelastic reduced-order
  models and application to an f-16 configuration}, AIAA Journal, 45 (2007),
  pp.~1244--1257, \url{https://doi.org/10.2514/1.24512},
  \url{https://doi.org/10.2514/1.24512},
  \url{https://arxiv.org/abs/https://doi.org/10.2514/1.24512}.

\bibitem{LuPengWei2019}
{\sc C.~Lu, X.~Peng, and Y.~Wei}, {\em Low-rank tensor completion with a new
  tensor nuclear norm induced by invertible linear transforms}, in 2019
  IEEE/CVF Conference on Computer Vision and Pattern Recognition (CVPR), 2019,
  pp.~5989--5997, \url{https://doi.org/10.1109/CVPR.2019.00615}.

\bibitem{Luo_2022}
{\sc Y.-S. Luo, X.-L. Zhao, T.-X. Jiang, Y.~Chang, M.~K. Ng, and C.~Li}, {\em
  Self-supervised nonlinear transform-based tensor nuclear norm for
  multi-dimensional image recovery}, IEEE Transactions on Image Processing, 31
  (2022), p.~3793–3808, \url{https://doi.org/10.1109/tip.2022.3176220},
  \url{http://dx.doi.org/10.1109/TIP.2022.3176220}.

\bibitem{ma2020randomized}
{\sc A.~Ma and D.~Molitor}, {\em Randomized kaczmarz for tensor linear
  systems}, 2020, \url{https://arxiv.org/abs/2006.01246}.

\bibitem{minka1997old}
{\sc T.~Minka}, {\em Old and new matrix algebra useful for statistics},
  September 1997,
  \url{https://www.microsoft.com/en-us/research/publication/old-new-matrix-algebra-useful-statistics/}.

\bibitem{murray2023randomized}
{\sc R.~Murray, J.~Demmel, M.~W. Mahoney, N.~B. Erichson, M.~Melnichenko, O.~A.
  Malik, L.~Grigori, P.~Luszczek, M.~Dereziński, M.~E. Lopes, T.~Liang,
  H.~Luo, and J.~Dongarra}, {\em Randomized numerical linear algebra : A
  perspective on the field with an eye to software}, 2023,
  \url{https://arxiv.org/abs/2302.11474}.

\bibitem{Newman2019:thesis}
{\sc E.~Newman}, {\em A Step in the Right Dimension: Tensor Algebra and
  Applications}, PhD thesis, Tufts University, May 2019.

\bibitem{NewmanTNN2018}
{\sc E.~Newman, L.~Horesh, H.~Avron, and M.~Kilmer}, {\em Stable tensor neural
  networks for rapid deep learning}.
\newblock Available at arXiv:1811.06569, November, 2018.

\bibitem{Newman2024:TNN}
{\sc E.~Newman, L.~Horesh, H.~Avron, and M.~E. Kilmer}, {\em Stable tensor
  neural networks for efficient deep learning}, Frontiers in Big Data, 7
  (2024), \url{https://doi.org/10.3389/fdata.2024.1363978},
  \url{https://www.frontiersin.org/articles/10.3389/fdata.2024.1363978}.

\bibitem{Newman2020}
{\sc E.~Newman and M.~E. Kilmer}, {\em Nonnegative tensor patch dictionary
  approaches for image compression and deblurring applications}, SIAM Journal
  on Imaging Sciences, 13 (2020), \url{https://doi.org/10.1137/19M1297026}.

\bibitem{Newman2021:GNvpro}
{\sc E.~Newman, L.~Ruthotto, J.~Hart, and B.~van Bloemen~Waanders}, {\em Train
  like a (var)pro: Efficient training of neural networks with variable
  projection}, SIAM Journal on Mathematics of Data Science, 3 (2021),
  \url{https://doi.org/10.1137/20m1359511}.

\bibitem{OLearyRust2013:VarPro}
{\sc D.~P. O'Leary and B.~W. Rust}, {\em Variable projection for nonlinear
  least squares problems}, Comput. Optim. Appl., 54 (2013), pp.~579--593,
  \url{https://doi.org/10.1007/s10589-012-9492-9},
  \url{https://doi.org/10.1007/s10589-012-9492-9}.

\bibitem{Oseledets2011}
{\sc V.~Oseledets}, {\em Tensor-train decomposition}, SIAM Journal of
  Scientific Computing, 33 (2011), pp.~2295--2317.

\bibitem{10125009}
{\sc K.~Pena-Pena, D.~L. Lau, and G.~R. Arce}, {\em t-hgsp: Hypergraph signal
  processing using t-product tensor decompositions}, IEEE Transactions on
  Signal and Information Processing over Networks, 9 (2023), pp.~329--345,
  \url{https://doi.org/10.1109/TSIPN.2023.3276687}.

\bibitem{yahoo}
{\sc M.~Perlin}, {\em yfR: Downloads and Organizes Financial Data from Yahoo
  Finance}, 2023, \url{https://github.com/ropensci/yfR}.
\newblock R package version 1.1.0.

\bibitem{Petersen07thematrix}
{\sc K.~B. Petersen and M.~S. Pedersen}, {\em The matrix cookbook}, 2007.

\bibitem{Recht2010:nuclearNorm}
{\sc B.~Recht, M.~Fazel, and P.~A. Parrilo}, {\em Guaranteed minimum-rank
  solutions of linear matrix equations via nuclear norm minimization}, SIAM
  Review, 52 (2010), pp.~471--501, \url{https://doi.org/10.1137/070697835},
  \url{https://doi.org/10.1137/070697835},
  \url{https://arxiv.org/abs/https://doi.org/10.1137/070697835}.

\bibitem{Sirovich1987}
{\sc L.~Sirovich}, {\em Turbulence and the dynamics of coherent structures part
  iii: Dynamics and scaling}, Quarterly of Applied Mathematics, 45 (1987),
  pp.~583--590, \url{http://www.jstor.org/stable/43637459} (accessed
  2024-06-07).

\bibitem{SongEtAl2020:tensorNuclearNorm}
{\sc G.~Song, M.~K. Ng, and X.~Zhang}, {\em Robust tensor completion using
  transformed tensor singular value decomposition}, Numerical Linear Algebra
  with Applications, 27 (2020), p.~e2299,
  \url{https://doi.org/https://doi.org/10.1002/nla.2299},
  \url{https://onlinelibrary.wiley.com/doi/abs/10.1002/nla.2299},
  \url{https://arxiv.org/abs/https://onlinelibrary.wiley.com/doi/pdf/10.1002/nla.2299}.

\bibitem{Townsend2016}
{\sc J.~Townsend}, {\em Differentiating the singular value decomposition},
  tech. report, 2016,
  \url{https://j-towns.github.io/papers/svd-derivative.pdf}.

\bibitem{Townsend:22}
{\sc O.~Townsend, S.~Gazzola, S.~Dolgov, and P.~Quinn}, {\em Undersampling
  raster scans in spectromicroscopy for a reduced dose and faster
  measurements}, Opt. Express, 30 (2022), pp.~43237--43254,
  \url{https://doi.org/10.1364/OE.471663},
  \url{https://opg.optica.org/oe/abstract.cfm?URI=oe-30-24-43237}.

\bibitem{Tucker1966}
{\sc L.~R. Tucker}, {\em Some mathematical notes on three-mode factor
  analysis}, Psychometrika, 31 (1966), pp.~279--311.

\bibitem{vanLeeuwenAravkin2021:nonsmoothVarPro}
{\sc T.~van Leeuwen and A.~Y. Aravkin}, {\em Variable projection for nonsmooth
  problems}, SIAM Journal on Scientific Computing, 43 (2021), pp.~S249--S268,
  \url{https://doi.org/10.1137/20M1348650},
  \url{https://arxiv.org/abs/https://doi.org/10.1137/20M1348650}.

\bibitem{Wan2019:complexSVD}
{\sc Z.-Q. Wan and S.-X. Zhang}, {\em Automatic differentiation for complex
  valued svd}, 2019, \url{https://doi.org/10.48550/ARXIV.1909.02659},
  \url{https://arxiv.org/abs/1909.02659}.

\bibitem{WuEtAl2022}
{\sc T.~Wu, J.~Fan, X.~Jize, and W.~L. Woo}, {\em Low-rank tensor completion
  based on self-adaptive learnable transforms}, IEEE Transactions on Neural
  Networks and Learning Systems, PP (2022), pp.~1--13,
  \url{https://doi.org/10.1109/TNNLS.2022.3215974}.

\bibitem{9313596}
{\sc H.-J. Yang, Y.-Y. Zhao, J.-X. Liu, Y.-X. Lei, J.-L. Shang, and X.-Z.
  Kong}, {\em Sparse regularization tensor robust pca based on t-product and
  its application in cancer genomic data}, in 2020 IEEE International
  Conference on Bioinformatics and Biomedicine (BIBM), 2020, pp.~2131--2138,
  \url{https://doi.org/10.1109/BIBM49941.2020.9313596}.

\bibitem{8421595}
{\sc M.~Yin, J.~Gao, S.~Xie, and Y.~Guo}, {\em Multiview subspace clustering
  via tensorial t-product representation}, IEEE Transactions on Neural Networks
  and Learning Systems, 30 (2019), pp.~851--864,
  \url{https://doi.org/10.1109/TNNLS.2018.2851444}.

\bibitem{zhang2016randomized}
{\sc J.~Zhang, A.~K. Saibaba, M.~Kilmer, and S.~Aeron}, {\em A randomized
  tensor singular value decomposition based on the t-product}, 2016,
  \url{https://arxiv.org/abs/1609.07086}.

\bibitem{ZhangEtAl2014:tensorNuclearNorm}
{\sc Z.~Zhang, G.~Ely, S.~Aeron, N.~Hao, and M.~Kilmer}, {\em Novel methods for
  multilinear data completion and de-noising based on tensor-svd}, in 2014 IEEE
  Conference on Computer Vision and Pattern Recognition, 2014, pp.~3842--3849,
  \url{https://doi.org/10.1109/CVPR.2014.485}.

\end{thebibliography}
